\documentclass[reqno,12pt]{amsart}

\usepackage{graphics}
\usepackage{amssymb, amsmath}
\usepackage{version}
\usepackage{fancybox}
\usepackage{bm}
\usepackage{xcolor}
\usepackage{mathrsfs}
\usepackage{mathtools} 

\theoremstyle{plain}
\newtheorem{theorem}{\sc Theorem}[section]
\newtheorem{lemma}[theorem]{\sc Lemma}
\newtheorem{claim}{\it Claim}

\newtheorem{proposition}[theorem]{\sc Proposition}
\newtheorem{corollary}[theorem]{\sc Corollary}

\theoremstyle{definition}

\theoremstyle{remark}
\newtheorem{remark}[theorem]{\sc Remark}

\makeatletter

\@addtoreset{equation}{section}

\renewcommand{\d}{\text{\rm d}}
\newcommand{\D}{\text{\rm D}}
\newcommand{\vep}{\varepsilon}

\newcommand{\e}{{\rm e}}

\newcommand{\R}{\mathbb{R}}
\newcommand{\N}{\mathbb{N}}

\newcommand{\Ro}{\Phi}

\newcommand{\Lphi}{\mathcal L_\phi}
\newcommand{\Ls}{\mathcal L_s}
\renewcommand{\H}{\mathcal{H}}

\newcommand{\M}{\mathcal{M}}

\allowdisplaybreaks[4]

\newcommand{\prf}[1]{}  
\begin{document}
\title[Fast diffusion in thin annuli]{Optimal rate of convergence to nonradial asymptotic profiles for fast diffusion in thin annuli\prf{\\{\tt [Detailed version]}}}
%
%
\author{Goro Akagi}
\address[Goro Akagi]{Mathematical Institute and Graduate School of Science, Tohoku University, Aoba, Sendai 980-8578, Japan}
\email{goro.akagi@tohoku.ac.jp}
\author{Norihisa Ikoma}
\address[Norihisa Ikoma]{Department of Mathematics, Faculty of Science and Technology, Keio University
Yagami Campus, 3-14-1 Hiyoshi, Kohoku-ku, Yokohama, Kanagawa 223-8522, Japan}
\email{ikoma@math.keio.ac.jp}
\author{Yasunori Maekawa}
\address[Yasunori Maekawa]{Department of Mathematics, Kyoto University, Kitashirakawa Oiwakecho, Sakyo-ku, Kyoto 606-8502, Japan}
\email{maekawa.yasunori.3n@kyoto-u.ac.jp}
\thanks{G.A.~is supported by JSPS KAKENHI Grant Numbers JP24H00184, JP21KK0044 and JP25K00920. N.I.~is supported by JSPS KAKENHI Grant Number JP24K06802. Y.M.~is supported by JSPS KAKENHI Grant Numbers 20H00118, 19H05597, 21H04433, 21H00991, 20K03698 and 23H01082. 
}
\date{\today}
\maketitle
\begin{abstract}
Quantitative estimates on the rate of convergence to asymptotic profiles for the fast diffusion equation in bounded domains have so far been obtained almost exclusively for \emph{nondegenerate} profiles. In contrast, least-energy asymptotic profiles in thin annuli are nonradial and hence form a continuum under rotation; such profiles are therefore \emph{degenerate}. This paper establishes the optimal convergence rate toward such profiles in thin annuli, both in two dimensions and in higher dimensions in the nearly Sobolev-critical regime. To this end, the kernel of the linearized operator is characterized and quantitative gradient inequalities are derived. The rate of convergence is also verified for positive asymptotic profiles in general bounded domains complying with the so-called Morse--Bott condition. On the other hand, the proof of optimality of the rate relies crucially on the particular geometry of thin annular domains. Finally, a symmetry-induced jump in the optimal convergence rate under domain perturbations is shown to occur, whereas no such phenomenon occurs for the linear diffusion equation or for the porous medium equation.
\end{abstract}

\section{Introduction}\label{S:Intro}

Berryman and Holland~\cite{BH80} studied the asymptotic behavior of (positive classical) solutions to the Cauchy--Dirichlet problem for the \emph{Fast Diffusion Equation} posed on bounded domains,
\begin{alignat}{3}
\partial_t (|u|^{q-2}u) &= \Delta u \quad &&\mbox{ in } \Omega \times \R_+,\label{u-pde}\\
u &= 0 \quad &&\mbox{ on } \partial \Omega \times \R_+,\label{u-bc}\\
u(\cdot,0) &= u_0 \quad &&\mbox{ in } \Omega,\label{u-ic}
\end{alignat}
where $\Omega$ is a bounded domain in $\R^d$ with boundary $\partial \Omega$ and $\R_+ = (0,\infty)$, under the following assumptions:
\begin{equation}\label{hypo}
2 < q < 2^* := \frac{2d}{(d-2)_+} \quad \mbox{ and } \quad u_0 \in H^1_0(\Omega) \setminus \{0\}.
\end{equation}
Here the condition $q > 2$ corresponds to the fast diffusion regime, whereas the case $1 < q < 2$ corresponds to the \emph{porous medium equation} (see~\cite{Vazquez,Vaz06,GalVaz04}). That study was motivated by the aim of providing a theoretical account for the Wisconsin octupole experiments based on the Okuda--Dawson model (see~\cite{Okuda-Dawson}), which describes anomalous diffusion of hydrogen plasma across a poloidal octupole magnetic field (see also~\cite{BH78,BH82}). The main results originating from~\cite{BH80}, together with a couple of subsequent refinements, can be summarized as follows:

\begin{theorem}[Finite-time extinction and asymptotic profiles of solutions~\cite{BH80,Kwong88-3,SavareVespri,AK13}]\label{T:ap}
Assume that \eqref{hypo} is satisfied and let $u = u(x,t)$ be an energy solution to the Cauchy--Dirichlet problem \eqref{u-pde}--\eqref{u-ic}. Then there exist a positive number $t_*$ and positive constants $c_1, c_2$ such that
$$
c_1 (t_* - t)_+^{1/(q-2)} \leq \|u(\cdot,t)\|_{H^1_0(\Omega)} \leq c_2 (t_* - t)_+^{1/(q-2)} \quad \mbox{ for } \ t \geq 0.
$$
A similar estimate also holds true for the $L^q$-norm of $u(\cdot,t)$. Moreover, for every increasing sequence $(t_n)$ in $\R_+$ with $t_n \to t_*$, there exist a subsequence $(n_k)$ of $(n)$ and $\phi \in H^1_0(\Omega) \setminus \{0\}$ such that
\begin{equation*}
(t_* - t_{n_k})_+^{-1/(q-2)} u(\cdot,t_{n_k}) \to \phi \quad \mbox{ strongly in } H^1_0(\Omega)
\end{equation*}
and $\phi$ is a nontrivial weak solution to the Dirichlet problem,
\begin{equation}\label{elliptic}
-\Delta \phi = \mathsf{c}_q |\phi|^{q-2}\phi \ \mbox{ in } \Omega, \quad \phi = 0 \ \mbox{ on } \partial \Omega,
\end{equation}
where $\mathsf{c}_q$ is the positive constant given by
\begin{equation}\label{cq}
\mathsf{c}_q := \frac{q-1}{q-2} > 0.
\end{equation}
In what follows, $t_*$ and $\phi$ are called the \emph{extinction time} {\rm (}often denoted by $t_*(u_0)${\rm )} and an \emph{asymptotic profile} of $u$, respectively.
\end{theorem}

Berryman and Holland~\cite{BH80} proved the above result for \emph{positive classical solutions} to the Cauchy--Dirichlet problem \eqref{u-pde}--\eqref{u-ic} (cf.~earlier works~\cite{Sabinina62,BC,Diaz88,HerreroVazquez88} on finite-time extinction). It was subsequently extended to energy solutions in~\cite{Kwong88-3,SavareVespri,AK13} (see also~\cite{DKV91,BGV08,BV10,JiXi23,JiXi22+}). On the other hand, while the well-posedness of \eqref{u-pde}--\eqref{u-ic} has been extensively studied (see, e.g.,~\cite{Vazquez}), the classical regularity of nonnegative weak solutions to \eqref{u-pde}--\eqref{u-ic} was mentioned as a conjecture in~\cite{BH80} and remained as an open question for many years. This conjecture was recently settled affirmatively by Jin and Xiong~\cite{JiXi23,JiXi22+} (cf.~\cite{DKV91}). We also refer the reader to~\cite{AS25}, in which the well-posedness result in an energy class is established in a self-contained fashion for classical and fractional porous medium and fast diffusion equations, even allowing for (possibly) sign-changing initial data.

In Theorem \ref{T:ap}, only \emph{quasiconvergence}, that is, convergence along a subsequence of $(t_n)$, is asserted. In~\cite{BH80}, the convergence without extracting any subsequence is also proved for the one-dimensional case $d = 1$, where solutions to \eqref{elliptic} are isolated in $H^1_0(\Omega)$. As for higher-dimensional cases, regardless of isolation, Feireisl and Simondon~\cite{FeiSim00} established the (uniform) convergence,
\begin{equation}\label{FeiSim}
(t_* - t)_+^{-1/(q-2)} u(\cdot,t) \to \phi \quad \mbox{ uniformly on } \overline\Omega \ \mbox{ as } \ t \nearrow t_*,
\end{equation}
provided that $\partial \Omega$ is smooth and that $u_0$ is nonnegative, by developing the so-called \L ojasiewicz--Simon gradient inequality (see also~\cite{AS25}). Furthermore, Bonforte, Grillo and V\'azquez~\cite{BGV12} derived a \emph{relative error convergence} to the positive asymptotic profile $\phi$,
$$
\frac{(t_* - t)_+^{-1/(q-2)} u(\cdot,t) - \phi}{\phi} \to 0 \quad \mbox{ uniformly on } \Omega \ \mbox{ as } \ t \nearrow t_*,
$$
by combining the uniform convergence \eqref{FeiSim} with the \emph{global Harnack principle} established in~\cite{DKV91} and by constructing sub- and supersolutions to \eqref{u-pde}--\eqref{u-ic}.

Rescaling $u = u(x,t)$ as
\begin{equation}\label{v}
v(x,s) := (t_* - t)^{-1/(q-2)} u(x,t), \quad s := \log \frac{t_*}{t_* - t} \quad (0 \leq t < t_*),
\end{equation}
one can rewrite \eqref{u-pde}--\eqref{u-ic} as
\begin{alignat}{3}
\partial_s (|v|^{q-2}v) &= \Delta v + \mathsf{c}_q |v|^{q-2}v \quad &&\mbox{ in } \Omega \times \R_+,\label{v-pde}\\
v &= 0 \quad &&\mbox{ on } \partial \Omega \times \R_+,\label{v-bc}\\
v(\cdot,0) &= v_0 \quad &&\mbox{ in } \Omega,\label{v-ic}
\end{alignat}
where $v_0 \in H^1_0(\Omega) \setminus \{0\}$ is given by
\begin{equation}\label{v0}
v_0 := t_*(u_0)^{-1/(q-2)}u_0.
\end{equation}
We define a functional $J : H^1_0(\Omega) \to \R$ by
$$
J(w) := \frac 12 \int_{\Omega} |\nabla w(x)|^2 \, \d x - \frac{\mathsf{c}_q}{q} \int_{\Omega} |w(x)|^q \, \d x \quad \mbox{ for } \ w \in H^1_0(\Omega).
$$
Then \eqref{v-pde}--\eqref{v-ic} can be represented as the Cauchy problem for the generalized gradient flow equation,
\begin{equation}\label{ggf}
\partial_s (|v|^{q-2}v)(s) = - J'(v(s)) \ \mbox{ in } H^{-1}(\Omega) \ \mbox{ for } \ s > 0, \quad v(0) = v_0,
\end{equation}
where $v(s) := v(\cdot,s)$ for $s \geq 0$. Hence the convergence of $(t_* - t)_+^{-1/(q-2)} u(\cdot,t)$ to an asymptotic profile $\phi$ as $t \to t_*$ can be interpreted as the convergence 
\begin{equation}\label{v-conv} 
 v(\cdot,s) \to \phi \quad \mbox{ strongly in } H^1_0(\Omega) \ \mbox{ as } \ s \to \infty
\end{equation}
of the rescaled solution $v(\cdot,s)$ to an equilibrium $\phi$ of the Cauchy--Dirichlet problem \eqref{v-pde}--\eqref{v-ic}, that is, a critical point of $J$. 

Stability analysis of asymptotic profiles has also been studied in~\cite{AK13,AK14,A16}. Set a phase set $\mathcal{X}$ by
\begin{align}
\mathcal{X} := \left\{ t_*(u_0)^{-1/(q-2)}u_0 \colon u_0 \in H^1_0(\Omega) \setminus \{0\} \right\}\label{X}
= \left\{ v_0 \in H^1_0(\Omega) \colon t_*(v_0) = 1 \right\}
\end{align}
(see~\cite{AK13} for more details). Then the initial datum $v_0$ given by \eqref{v0} belongs to $\mathcal{X}$. Moreover, $\mathcal{X}$ is invariant under the evolution generated by \eqref{v-pde}--\eqref{v-ic}, that is, if $v_0 \in \mathcal{X}$, then $v(s)$ lies in $\mathcal{X}$ for any $s \geq 0$. The notions of stability of asymptotic profiles are defined in the sense of Lyapunov stability of equilibria for the dynamical system generated by \eqref{v-pde}--\eqref{v-ic} on the phase set $\mathcal{X}$. In particular, it was proved in~\cite{A16} that every least-energy solution to \eqref{elliptic} is stable in the sense of asymptotic profiles. Here and henceforth, a nontrivial solution $\phi$ of \eqref{elliptic} is called a \emph{least-energy solution} of \eqref{elliptic}, if $\phi$ attains the infimum of the energy functional $J$ among all nontrivial solutions of \eqref{elliptic}. 

We now turn to quantitative studies of convergence to asymptotic profiles. 
In contrast to the porous medium equation, for which the optimal convergence rate was obtained already in the 1980s by Aronson and Peletier~\cite{ArPel81} (see also~\cite{JiROXi22+} for finer asymptotics), progress in the fast diffusion case remained limited for many years. Bonforte and Figalli~\cite{BF21} singled out a convergence rate that is expected to be sharp, as suggested by the formal linearization around $\phi$ (see~\cite[Section 2]{BF21}), under the assumption that $\phi$ is \emph{nondegenerate}. Here and henceforth, the nondegeneracy of $\phi$ means that the linearized operator $\Lphi : H^1_0(\Omega) \to H^{-1}(\Omega)$ at $\phi$, defined by
\begin{equation}\label{Lphi}
\Lphi w = - \Delta w - \mathsf{c}_q (q-1) |\phi|^{q-2} w \quad \mbox{ for } \ w \in H^1_0(\Omega),
\end{equation}
has trivial kernel. If $\phi$ is nondegenerate, then $\phi$ is isolated among nontrivial solutions to \eqref{elliptic}. Moreover, Bonforte and Figalli~\cite{BF21} were the first to establish this rate via a nonlinear entropy method, proving that
\begin{equation}\label{ent-conv}
\left( \int_\Omega |v(x,s)-\phi(x)|^2 \phi^{q-2}(x) \, \d x \right)^{1/2} \leq C \exp \left( -\frac{\nu_k}{q-1} s\right) \quad \mbox{ for } \ s > 0,
\end{equation}
whenever $\partial \Omega$ is smooth, $v$ is nonnegative (hence $\phi$ is positive) and $\phi$ is nondegenerate. Here $\nu_k$ denotes the smallest positive eigenvalue of the weighted eigenvalue problem,
\begin{equation}\label{ev}
\Lphi e_j = \nu_j |\phi|^{q-2} e_j \ \mbox{ in } \Omega, \quad e_j = 0 \ \mbox{ on } \partial \Omega,
\end{equation}
where the eigenvalues $\{\nu_j\}_{j \in \N}$ are arranged in nondecreasing order, that is, $\nu_1 < \nu_2 \leq \cdots \leq \nu_j \leq \cdots \to \infty$. They also developed a \emph{quantitative relative error estimate}, which plays a crucial role in their analysis. Furthermore, by exploiting the optimal boundary regularity result of~\cite{JiXi23}, the topology of the convergence in \eqref{ent-conv} can be improved when $\partial \Omega$ is sufficiently smooth. See also~\cite{BF23+} for a survey of these developments and~\cite{JiXi20+} for the Sobolev critical case.

In~\cite{A21}, the first author of the present paper established that
\begin{equation}\label{H10-conv}
\|v(\cdot,s) - \phi\|_{H^1_0(\Omega)} \leq C \exp \left( -\frac{\nu_k}{q-1} s\right) \quad \mbox{ for } \ s > 0,
\end{equation}
where $\nu_k$ denotes the smallest positive eigenvalue of \eqref{ev}, provided that $\partial \Omega$ is of class $C^{1,1}$, that $v$ is nonnegative and that $\phi$ is nondegenerate. The method of proof in~\cite{A21} differs from the nonlinear entropy method used in~\cite{BF21} and relies on an energy method combined with a \emph{quantitative gradient inequality}, which can be viewed as a quantitative counterpart of gradient inequalities such as the \L ojasiewicz--Simon inequality. On the other hand, the quantitative relative error estimate developed in~\cite{BF21} is still used. Subsequently, in~\cite{AM23+}, the first and third authors of the present paper extended the estimate \eqref{H10-conv} to (possibly) sign-changing solutions by refining the method developed in~\cite{A21} so as to exclude the use of the quantitative relative error estimate. Furthermore, the optimality of the rates in \eqref{ent-conv} and \eqref{H10-conv} was first established in~\cite{AM23+} by proving an exponential stability result for least-energy asymptotic profiles (cf.~\cite{A21}) and by introducing a change of variables that decomposes \eqref{v-pde}--\eqref{v-ic} into semilinear problems. Finally, the approach of~\cite{A21} has also been developed in a different direction in~\cite{KoYu24}, which will be reviewed later.

On the other hand, all these results, except for~\cite{McCann23} and~\cite{KoYu24}, have been established for \emph{nondegenerate} asymptotic profiles only. Although the nondegeneracy of asymptotic profiles is generic (see~\cite{SaTe,BF21}), it does not necessarily hold for all domains. In fact, such degeneracy can occur, as illustrated by a concrete example given below. Let $A_R$ be the annular domain in $\R^d$ with $d \geq 2$ defined by
\begin{equation}\label{annulus}
A_R := \left\{ x \in \R^d \colon R < |x| < R + 1 \right\}
\end{equation}
for $R > 0$ and consider the Cauchy--Dirichlet problem \eqref{u-pde}--\eqref{u-ic} in $\Omega = A_R$. Although $A_R$ is radially symmetric, all least-energy solutions to \eqref{elliptic} in $\Omega = A_R$ are no longer radially symmetric when $R$ is sufficiently large (see~\cite{Coffman} for $N = 2$ and also~\cite{Li,Byeon} for $N \geq 3$). As a consequence, the set of least-energy solutions forms a continuum under \emph{rotations}. Hence $\phi$ is not isolated in $H^1_0(A_R)$ among asymptotic profiles and is therefore degenerate; that is, the linearized operator $\Lphi$ at $\phi$ has a nontrivial kernel. Consequently, least-energy solutions to \eqref{elliptic} in $\Omega = A_R$ lie beyond the scope of the previous results mentioned above. This paper is concerned with such annular settings.

Rates of convergence have also been discussed for degenerate asymptotic profiles in a few papers. Choi, McCann and Seis~\cite{McCann23} obtained a dichotomy result of the convergence rate when $\partial \Omega$ is sufficiently smooth and $v$ is nonnegative (see also~\cite{CS23+}). Their result shows that exactly one of the following alternatives holds:
\begin{itemize}
 \item the relative entropy decays algebraically or slower:
$$
\left( \int_\Omega |v(x,s)-\phi(x)|^2 \phi^{q-2}(x) \, \d x \right)^{1/2} \geq C s^{-1} \quad \mbox{ for } \ s \gg 1;
$$
 \item the relative error decays exponentially or faster:
$$
\left\|\frac{v(\cdot,s)}{\phi} - 1 \right\|_{L^\infty(\Omega)} \leq C \exp\left(- \frac{\nu_k}{q-1} s\right) \quad \mbox{ for } \ s \gg 1,
$$
where $\nu_k$ is the smallest positive eigenvalue of \eqref{ev}.
\end{itemize}
The proof of the above dichotomy result also differs from those in~\cite{BF21} and~\cite{A21} and crucially relies on a parabolic regularity estimate (see~\cite[Proposition 4.4 and Section 5]{McCann23}). Moreover, some sufficient conditions for the latter case are also provided. However, these results do not appear to be sufficient to draw a conclusion for the annular-domain case mentioned above. Moreover, the optimality of the convergence rate is not addressed in the paper.

K\"onig and Yu~\cite{KoYu24} studied the rate of convergence to asymptotic profiles for the Sobolev-critical fast diffusion equation posed on the whole space $\R^d$,
$$
\partial_t (|u|^{2^*-2}u) = \Delta u \quad \mbox{ in } \R^d \times \R_+,
$$
by refining the energy method developed in~\cite{A21}. It is worth mentioning that, in this setting, positive asymptotic profiles are uniquely determined, up to \emph{translation} and \emph{scaling}, by the so-called Talenti function $\kappa_d (1 + |x|^2)^{-(d-2)/2}$ with some constant $\kappa_d > 0$. Owing to the translation and scaling invariance, each positive asymptotic profile is not isolated and is therefore degenerate. Hence, the method of proof used in~\cite{A21} cannot be directly applied to this setting. K\"onig and Yu~\cite{KoYu24} modified the argument in~\cite{A21} and established a quantitative gradient inequality even for such a degenerate profile. Although the symmetry underlying the degeneracy is different -- translation and scaling invariance in~\cite{KoYu24} versus rotational invariance in the present setting -- the present study is also inspired by~\cite{KoYu24}. The method developed in the present paper to overcome difficulties arising from the degeneracy, however, is slightly different from that used in~\cite{KoYu24} and is closer to the argument originally developed in~\cite{A21}.

The main purpose of the present paper is to determine the optimal rate of the convergence \eqref{v-conv} in thin annuli, in which least-energy asymptotic profiles are degenerate. In what follows, we set
$$
\Omega = A_R
$$
unless otherwise specified. Let $\phi$ be a positive solution to \eqref{elliptic}. For $j \geq 1$, each eigenvalue $\nu_j$ of the linearized operator $\Lphi$ at $\phi$ can be written as 
$$
\nu_j = \mu_j - \mathsf{c}_q(q-1),
$$
where $\mu_j$ is the $j$-th eigenvalue of the weighted eigenvalue problem,
\begin{equation}\label{ep}
- \Delta e_j = \mu_j \phi^{q-2} e_j \ \mbox{ in } A_R, \quad e_j = 0 \ \mbox{ on } \partial A_R.
\end{equation}
The eigenvalues are arranged in nondecreasing order, that is, $0 < \mu_1 < \mu_2 \leq \cdots$, with each eigenvalue repeated according to its multiplicity.

In order to state main results of the present paper, we give an assumption:
\begin{enumerate}
\item[(A)] There exists a constant $R_0 > 0$ such that, for any $R \geq R_0$ and any least-energy solution $\phi > 0$ of \eqref{elliptic} in $\Omega = A_R$, it holds that
\begin{equation}\label{kdim-AR}
\mathrm{dim} \, \ker \Lphi = d-1,
\end{equation}
where $\Lphi : H^1_0(A_R) \to H^{-1}(A_R)$ is the linearized operator at $\phi$ defined by \eqref{Lphi} with $\Omega = A_R$.
\end{enumerate}
We are ready for stating main results of the present paper. 
First of all, the following proposition shows that assumption (A) is satisfied, whenever condition (S) below is fulfilled:
\begin{enumerate}
 \item[(S)] Let $w \in H^1_0(S)$ be a least-energy solution to the Dirichlet problem,
\begin{equation}\label{eq-S}
- \Delta w = |w|^{q-2}w \ \mbox{ in } S, \quad w = 0 \ \mbox{ on } \partial S,
\end{equation}
where $S := \{y = (y',y_d) \in \R^d \colon y' \in \R^{d-1}, \ 0 < y_d < 1\}$, and let $\mathcal{L}_w : H^1_0(S) \to H^{-1}(S)$ be the linearized operator at $w$ given by $\mathcal{L}_w z = - \Delta z - (q-1) |w|^{q-2} z$ for $z \in H^1_0(S)$. Then it holds that 
\begin{equation*}
\dim \ker \mathcal{L}_w = d-1.
\end{equation*} 
\end{enumerate}

\begin{proposition}[Kernels of linearized operators at least-energy profiles in thin annuli]\label{P:nondeg}
Let $d \geq 2$ and assume {\rm (S)} is satisfied. Then there exists a constant $R_0 > 0$ such that, for any $R \geq R_0$ and any least-energy solution $\phi > 0$ of \eqref{elliptic} with $\Omega = A_R$, \eqref{kdim-AR} is satisfied.
\end{proposition}

Assumption (S) holds true for every $q > 2$ when $d = 2$, and for $q \in (2,2^*)$ close enough to $2^*$ when $d \geq 3$. Indeed, we recall the following facts: 

\begin{lemma}[{\cite[Lemma 5]{Dancer89} (see also~\cite[Lemma 2.1]{BKP13})}]\label{L:Dancer}
Assume $d = 2$ and let $S = \{(s,t) \in \R^2 \colon s \in \R, \ 0 < t < 1\}$. For any $2 < q < \infty$, let $w \in H^1_0(S)$ be the unique positive solution to \eqref{eq-S} up to translation. Let $z \in H^1_0(S)$ solve the Dirichlet problem,
\begin{equation}\label{S-lin}
- \Delta z - (q-1) |w|^{q-2} z = 0 \ \mbox{ in } S, \quad z = 0 \ \mbox{ on } \partial S.
\end{equation}
Then $z = c \partial_s w$ for some $c \in \R$, and hence, $\dim \ker \mathcal{L}_w = 1$.
\end{lemma}

\begin{lemma}[{\cite[Propositions B.1 and C.1.]{DaYa02}}\label{L:DaYa}]
Assume $d \geq 3$. For any $2 < q < 2^*$ close enough to $2^*$, let $w \in H^1_0(S)$ be a least-energy solution {\rm (}$=$ mountain-pass solution{\rm )} to \eqref{eq-S}. Let $z \in H^1_0(S)$ solve the Dirichlet problem \eqref{S-lin}. Then $z$ lies on $\mathrm{span} \{ \partial_{y_1} w , \ldots, \partial_{y_{d-1}} w \}$, and therefore, $\dim \ker \mathcal{L}_w = d-1$.
\end{lemma}

We next determine the convergence rate in \eqref{v-conv} under assumption (A).

\begin{theorem}[Rate of convergence in thin annuli]\label{T:rate}
Let $d \geq 2$ and suppose that {\rm (A)} is satisfied. Let $R \geq R_0$, where $R_0$ is given by {\rm (A)}, and let $\Omega = A_R$ be the annular domain in $\R^d$ defined by \eqref{annulus}. Let $v = v(x,s)$ be a nonnegative energy solution to \eqref{v-pde}--\eqref{v-ic} and let $\phi = \phi(x)$ be a least-energy solution to \eqref{elliptic} such that $v(\cdot,s) \to \phi$ strongly in $H^1_0(A_R)$ as $s \to \infty$. Then there exists a constant $C \geq 0$ such that
\begin{equation}\label{rate}
\|v(\cdot,s)-\phi\|_{H^1_0(A_R)} \leq C \exp \left( -\frac{\nu_{d+1}}{q-1} s\right)
\end{equation}
for $s \geq 0$. Here $\nu_{d+1}$ is the least positive eigenvalue of the linearized operator $\Lphi$.
\end{theorem}

The above theorem can be extended to more general settings, see Theorem \ref{T:rate-g} below, where the fast diffusion equation \eqref{u-pde}--\eqref{u-ic} is considered in general bounded domains $\Omega \subset \R^d$ and $\phi$ is assumed to comply with the Morse--Bott condition (MB), which will be stated in Section \ref{S:relerr} (see also Remark \ref{R:MB}). In Sections \ref{S:relerr}--\ref{S:pf}, we shall prove Theorem \ref{T:rate-g}, from which Theorem \ref{T:rate} follows as a corollary.

We also remark that, for general bounded $C^{1,1}$ domains $\Omega$ in $\R^d$, the energy solution $v = v(x,s)$ of \eqref{v-pde}--\eqref{v-ic} converges strongly in $H^1_0(\Omega)$ to a least-energy solution $\phi$ of \eqref{elliptic}, provided that the initial energy $J(v_0)$ is close enough to the \emph{least-energy} $d_J$ of $J$ defined by
$$
d_J := \inf\{J(\psi) \colon \psi \in H^1_0(\Omega) \mbox{ is a nontrivial solution to \eqref{elliptic}}\}
$$
(see Proposition \ref{P:LES-conv} below for more details).

We next discuss the optimality of the convergence rate in \eqref{rate}. To this end, slightly modifying the proof of Theorem \ref{T:rate} in \S \ref{S:relerr}--\ref{S:pf}, we establish a refined exponential convergence estimate near each least-energy asymptotic profile, which plays a crucial role to prove the optimality of the convergence rate. We stress that this estimate holds uniformly in a neighbourhood of each least-energy solution of \eqref{elliptic}, and moreover, the prefactor in the exponential convergence estimate is controlled by the initial discrepancy from the asymptotic profile. Such estimates have never been obtained in the previous works~\cite{BF21,A21,McCann23,KoYu24} except~\cite{AM23+}, which is concerned exclusively with \emph{nondegenerate} asymptotic profiles and hence not available for the present setting.

\begin{corollary}[A refined exponential convergence estimate]\label{C:rate}
Let $d \geq 2$ and suppose that {\rm (A)} is satisfied. Let $R \geq R_0$, where $R_0$ is given by {\rm (A)}, and let $\Omega = A_R$ be the annular domain in $\R^d$ defined by \eqref{annulus}. Let $\phi = \phi(x)$ be a positive least-energy solution to \eqref{elliptic}. Then there exist constants $\delta, C > 0$ such that any energy solution $v = v(x,s)$ to \eqref{v-pde}--\eqref{v-ic} satisfies
\begin{align*}
\|v(\cdot,s)-\phi\|_{H^1_0(A_R)} 
&\leq C \left[ \left( J(v_0) - J(\phi) \right)^{1/2} + \|v_0^{q/2} - \phi^{q/2}\|_{L^2(A_R)} \right] \exp \left( - \frac{\nu_{d+1}}{q-1}s \right)
\end{align*}
for all $s \geq 0$, provided that $v_0 \in \mathcal{X}$, $v_0 \geq 0$, $\|v_0-\phi\|_{H^1_0(A_R)} < \delta$ and $v(\cdot,s) \to \phi$ strongly in $H^1_0(A_R)$ as $s \to \infty$.
\end{corollary}

Then carefully constructing a family of proper initial data which satisfy the assumption of the above corollary and employing the argument developed in~\cite{AM23+}, we establish that

\begin{theorem}[Optimality of the convergence rate]\label{T:opt}
Let $d \geq 2$ and suppose that {\rm (A)} is satisfied. Let $R \geq R_0$, where $R_0$ is given by {\rm (A)}, and let $\Omega = A_R$ be the annular domain in $\R^{d}$ defined by \eqref{annulus}. Let $\phi = \phi(x)$ be a positive least-energy solution to \eqref{elliptic}. Then there exists a nonnegative initial datum $\bar u_0 \in H^1_0(A_R) \setminus \{0\}$ with the following property\/{\rm :} Let $\bar v = \bar v(x,s)$ be the unique energy solution to \eqref{v-pde}--\eqref{v-ic} for the initial datum $v_0 = t_*(\bar u_0)^{-1/(q-2)}\bar u_0$. Then $\bar v(\cdot,s)$ converges to $\phi$ strongly in $H^1_0(A_R)$ as $s \to \infty$, and moreover, there exist constants $c_1,c_2,c_3 > 0$ such that
\begin{align*}
c_1 \exp \left( -\frac{\nu_{d+1}}{q-1} s\right)
&\leq \left( \int_{A_R} |\bar v(x,s)- \phi(x)|^2 \phi^{q-2}(x) \, \d x \right)^{1/2}\\
&\leq c_2 \| \bar v(\cdot,s)- \phi \|_{H^1_0(A_R)} \leq c_3 \exp \left( -\frac{\nu_{d+1}}{q-1} s\right) \quad \mbox{ for } \ s \geq 0.
\end{align*}
In particular, the convergence rate of \eqref{rate} is optimal.
\end{theorem}

In particular, we also have
\begin{corollary}
Assume either $d = 2$ with $q \in (2,\infty)$ or $d \geq 3$ with $q \in (2,2^*)$ close enough to $2^*$. Let $R > 0$ be large enough and set $\Omega = A_R$. Then there exists $\delta_0 > 0$ satisfying the following\/{\rm :} for every nonnegative initial datum $v_0 \in \mathcal X$ with $J(v_0) < d_J + \delta_0$, there exist a least-energy solution $\phi$ of \eqref{elliptic} and a constant $C > 0$ such that the energy solution $v$ to \eqref{v-pde}--\eqref{v-ic} satisfies \eqref{rate} for $s \geq 0$. Moreover, the rate of convergence in \eqref{rate} is optimal in the sense of Theorem {\rm \ref{T:opt}}. 
\end{corollary}

In view of the above results, we observe \emph{symmetry-induced jumps} in the convergence rate under domain perturbations. More precisely, we can construct a sequence of smooth bounded domains $(\Omega_n)$ converging to a thin annular domain $A_R$, to which the above results apply, such that the exponential convergence rate in $\Omega_n$ degenerates to zero as $n \to \infty$ even though $\nu_{d+1}$ is positive (see \eqref{rate}). Such phenomena occur neither in the linear diffusion case ($q = 2$) nor in the porous medium case ($1 < q < 2$). The details are given in Section \ref{S:jump} (see, in particular, Remark \ref{R:jump} below). 

\bigskip
\noindent
{\bf Structure of the paper.} This paper consists of seven sections. Section \ref{S:ker} is devoted to proving Proposition \ref{P:nondeg}. Sections \ref{S:relerr}--\ref{S:pf} are devoted to a proof of Theorem \ref{T:rate-g} below, from which Theorem \ref{T:rate} and Corollary \ref{C:rate} follow as corollaries. More precisely, in Section \ref{S:relerr}, we discuss the Morse--Bott property and develop a relative error estimate, adapted and refined to suit the subsequent proof. In Section \ref{S:QGIa}, we establish a quantitative gradient inequality under the Morse--Bott condition. Our proof of the inequality is based on the method developed in~\cite{A21} along with a suitable modification to handle the degenerate profiles (cf.~\cite{KoYu24}). In contrast to~\cite{KoYu24}, where the target profile function admits an explicit representation (as Talenti functions), least-energy asymptotic profiles in thin annuli do not. Consequently, it is necessary to modify the argument in~\cite{KoYu24} and this issue is also relevant to the relative error estimate established in Section \ref{S:relerr}. Moreover, the proof of the quantitative gradient inequality in the present paper is closer to that of~\cite{A21} and may be simpler than that of~\cite{KoYu24}. In Section \ref{S:pf}, we give a proof of Theorem \ref{T:rate-g}, from which Theorem \ref{T:rate} follows immediately as a corollary, by combining the results in the preceding two sections. Moreover, we establish a refined exponential convergence estimate (see Corollary \ref{C:rate}), which will be needed to prove the optimality of the convergence rate. Owing to the lack of nondegeneracy of the asymptotic profile, the derivation of such a refined exponential convergence estimate is more delicate than the exponential stability results in~\cite{A21,AM23+}, where the asymptotic profile is assumed to be nondegenerate. Section \ref{S:opt} is devoted to a proof of Theorem \ref{T:opt} concerning the optimality of the convergence rate in Theorem \ref{T:rate}. A crucial and particularly delicate step in the proof is the construction of suitable initial data, which deserves special attention (cf.~\cite{AM23+}). In Section \ref{S:jump}, we discuss a symmetry-induced jump in the optimal convergence rate under domain perturbations. Finally, in the appendix, we provide proofs of some auxiliary results for the convenience of the reader. 

\bigskip

\noindent
{\bf Notation.} Let $\Omega$ be a (bounded) domain in $\R^d$. We denote by $H^1_0(\Omega)$ the closure of $C^\infty_c(\Omega)$ in $H^1(\Omega)$, endowed with the inner product
\begin{align*}
(u,v)_{H^1_0(\Omega)} := \int_{\Omega} \nabla u \cdot \nabla v \, \d x \quad \mbox{ for } \ u,v \in H^1_0(\Omega).
\end{align*}
The dual space $H^{-1}(\Omega)$ of $H^1_0(\Omega)$ is a Hilbert space endowed with the inner product
$$
(f,g)_{H^{-1}(\Omega)} := \langle f, (- \Delta)^{-1} g \rangle_{H^1_0(\Omega)} \quad \mbox{ for } \ f,g \in H^{-1}(\Omega),
$$
where $\langle \cdot, \cdot \rangle_{H^1_0(\Omega)}$ stands for the duality pairing. In addition, the operator $-\Delta$ defines a duality (Riesz) map between $H^1_0(\Omega)$ and $H^{-1}(\Omega)$. More precisely, it holds that
\begin{align*}
\|u\|_{H^1_0(\Omega)}^2 &= \|-\Delta u\|_{H^{-1}(\Omega)}^2 = \langle -\Delta u, u \rangle_{H^1_0(\Omega)},\\
\|f\|_{H^{-1}(\Omega)}^2 &= \|(-\Delta)^{-1} f\|_{H^1_0(\Omega)}^2 = \langle f, (-\Delta)^{-1} f \rangle_{H^1_0(\Omega)}
\end{align*}
for $u \in H^1_0(\Omega)$ and $f \in H^{-1}(\Omega)$. For $\phi \in H^1_0(\Omega)$ and $r > 0$, we write $B_{H^1_0(\Omega)}(\phi;r) := \{ w \in H^1_0(\Omega) \colon \|w - \phi\|_{H^1_0(\Omega)} < r\}$. Let $C_\Omega = C(d,\Omega,q)$ stand for the best constant of the Sobolev--Poincar\'e inequality, 
\begin{equation}\label{Sobolev}
\|w\|_{L^q(\Omega)} \leq C_\Omega \|\nabla w\|_{L^2(\Omega)} \quad \mbox{ for } \ w \in H^1_0(\Omega). 
\end{equation}
Throughout this paper, for a function $u = u(x,t)$, we write $u(t) := u(\cdot,t)$. In addition, $C$ denotes a generic nonnegative constant, which may vary from line to line, unless otherwise specified.

\section{Kernels of linearized operators at least-energy profiles}\label{S:ker}

In this section, we prove Proposition {\rm \ref{P:nondeg}}. To this end, by a simple scaling argument, it suffices to consider the Dirichlet problem
\begin{equation}\label{el-AR-d}
-\Delta u = |u|^{q-2}u \ \mbox{ in } A_R, \quad u = 0 \ \mbox{ on } \partial A_R
\end{equation}
posed on the annular domain $A_R \subset \R^d$ defined by \eqref{annulus} for $R > 0$, instead of \eqref{elliptic} with $\Omega = A_R$ (see below).

\begin{remark}[Simple scaling]\label{R:scaling} Let $\Omega$ be a domain in $\R^d$ and let $q \neq 2$. If $\phi$ solves \eqref{elliptic}, then $u := \mathsf{c}_q^{1/(q-2)} \phi$ solves
\begin{equation}\label{elliptic-0}
-\Delta u = |u|^{q-2}u \ \mbox{ in } \Omega, \quad u = 0 \ \mbox{ on } \partial \Omega.
\end{equation}
In particular, if $\phi$ is a least-energy solution of \eqref{elliptic}, then so is $u$ of \eqref{elliptic-0}. Moreover, the corresponding linearized operator $L_u := - \Delta - (q-1) |u|^{q-2}$ coincides with $\mathcal{L}_\phi$ defined by \eqref{Lphi}. 
\end{remark}

It is well known that least-energy solutions to \eqref{el-AR-d} are nonradial, provided that $R$ is large enough (see~\cite{Coffman} for $d = 2$ and also~\cite{Li,Byeon} for $d \geq 3$). We also consider the following auxiliary problem:
\begin{equation}\label{el-S-d}
-\Delta v = |v|^{q-2}v \ \mbox{ in } S, \quad v = 0 \ \mbox{ on } \partial S
\end{equation}
posed on the infinite strip domain $S := \{ y = (y', y_d) \in \R^d \colon 0 < y_d < 1 \}$. 

Define functionals $I_{A_R} : H^1_0(A_R) \to \R$ and $I_S : H^1_0(S) \to \R$ by
\begin{alignat*}{3}
 I_{A_R}(u) &:= \frac 12 \int_{A_R} |\nabla u(x)|^2 \, \d x - \frac 1q \int_{A_R} |u(x)|^q \, \d x \quad &&\mbox{ for } u \in H^1_0(A_R),\\
 I_{S}(v) &:= \frac 12 \int_{S} |\nabla v(y)|^2 \, \d y - \frac 1q \int_{S} |v(y)|^q \, \d x \quad &&\mbox{ for } v \in H^1_0(S)
\end{alignat*}
and define the \emph{least energies} of $I_{A_{R_n}}$ and $I_S$ among nontrivial solutions to \eqref{el-AR-d} and \eqref{el-S-d} by
\begin{align*}
 d_{A_R} &:= \inf \{ I_{A_R}(u) \colon u \in H^1_0(A_R) \setminus \{0\}, \ I_{A_R}'(u) = 0 \},\\
 d_{S} &:= \inf \{ I_{S}(v) \colon v \in H^1_0(S) \setminus \{0\}, \ I_{S}'(v) = 0 \},
\end{align*}
respectively. Moreover, we define the set of all \emph{least-energy solutions} to \eqref{el-AR-d} and \eqref{el-S-d} by
\begin{align*}
 \mathscr{L}_{A_R} &:= \{u \in H^1_0(A_R) \setminus \{0\} \colon I_{A_R}'(u) = 0 \ \mbox{ and } \ I_{A_R}(u) = d_{A_R}\},\\
 \mathscr{L}_{S} &:= \{v \in H^1_0(S) \setminus \{0\} \colon I_{S}'(v) = 0 \ \mbox{ and } \ I_{S}(v) = d_{S}\},
\end{align*}
respectively, and we also set
\begin{align*}
 \mathscr{L}_{A_R}^\circ &:= \{u \in \mathscr{L}_{A_R} \colon |u| \mbox{ takes a maximum on the positive $x_d$-axis} \},\\
 \mathscr{L}_{S}^\circ &:= \{v \in \mathscr{L}_{S} \colon |v| \mbox{ takes a maximum on the positive $x_d$-axis} \}.
\end{align*}
Here we remark that every least-energy solution is sign-definite, and that
\begin{equation}\label{Nehari}
\begin{aligned}
 d_{A_R} &= I_{A_R}(u) = \left(\frac12 - \frac1q\right) \|\nabla u\|_{L^2(A_R)}^2 = \left(\frac12 - \frac1q\right) \|u\|_{L^q(A_R)}^q \quad \mbox{ for } \ u \in \mathscr{L}_{A_R},\\
 d_S &= I_{S}(v) = \left(\frac12 - \frac1q\right) \|\nabla v\|_{L^2(S)}^2 = \left(\frac12 - \frac1q\right) \|v\|_{L^q(S)}^q \quad \mbox{ for } \ v \in \mathscr{L}_{S}.
\end{aligned}
\end{equation}

In what follows, $e_d$ denotes the $d$-th canonical basis vector of $\R^d$, that is, $e_d = (0,\ldots,0,1) \in \R^d$. The zero extensions of $u_R \in H^1_0(A_R)$ and $v \in H^1_0(S)$ belong to $H^1(\R^d)$. By abuse of notation, we use the same symbol to denote a function defined on $A_R$ (or on $S$) and its zero extension to $\R^d$. In order to prove Proposition \ref{P:nondeg}, we start with the following:

\begin{proposition}[Convergence of least-energy solutions]\label{P:les-d}
Let $(u_R)_{R > 0}$ be a family in $H^1(\R^d)$ such that $u_R \in \mathscr{L}_{A_R}^\circ$ for $R > 0$. Then it holds that
$$
\inf_{v \in \mathscr{L}_S^\circ} \|u_R - v(\,\cdot\, - R e_d)\|_{H^1(\R^d)} \to 0 \ \mbox{ as } \ R \to \infty.
$$
\end{proposition}

To prove Proposition \ref{P:les-d}, we first introduce some notation and preliminaries. Let $(\lambda_R,\varphi_R) \in (0,\infty) \times H^1_0(A_R)$ be the first eigenpair of the Dirichlet Laplacian, that is,
$$
-\Delta \varphi_R = \lambda_R \varphi_R \ \mbox{ in } A_R, \quad \varphi_R = 0 \ \mbox{ on } \partial A_R,
$$
such that $\varphi_R > 0$ and $\|\varphi_R\|_{L^\infty(A_R)} = 1$.
Then \prf{$\varphi_R$ is radial, that is, there exists a function $\psi_R \in C^2([0,1])$ such that $\varphi_R(x) = \psi_R(t)$ for $x \in \overline{A_R}$ and $t \in [0,1]$ satisfying $|x| = R + t$. Moreover, }the following fact may be well known (see, e.g.,~\cite[Lemma A.1]{SSLin}), but we also give a complete proof in Appendix \S \ref{apdx:comm} for the convenience of the reader.

\begin{lemma}\label{L:LapEV}
 It holds that $\lambda_R \to \pi^2$ \prf{and $\|\psi_R - \sin (\pi \,\cdot\,) \|_{C^2([0,1])} \to 0$ }as $R \to \infty$.
\end{lemma}

In particular, in view of the Poincar\'e inequality, we can take a constant $c_0 > 0$ independent of $R$ such that
\begin{equation}\label{Poincare-c0}
\|u\|_{L^2(A_R)}^2 \leq c_0 \|\nabla u\|_{L^2(A_R)}^2 \quad \mbox{ for } \ u \in H^1_0(A_R)
\end{equation}
for $R > 0$ large enough. 

Now, we have

\begin{lemma}\label{L:d_AR}
 It holds that
$$
0 < \liminf_{R \to \infty} d_{A_R} \leq \limsup_{R \to \infty} d_{A_R} \leq d_S < \infty.
$$
\end{lemma}

\begin{proof}
 We first recall that the least energy $d_{A_R}$ is characterized as follows:
\begin{align*}
d_{A_R} &= \inf_{\|u\|_{H^1_0(A_R)} = 1} \sup_{t \geq 0} I_{A_R}(tu) = \inf_{u \in H^1_0(A_R) \setminus \{0\}} \sup_{t \geq 0} I_{A_R}(tu) \\
&= \inf_{u \in H^1_0(A_R) \setminus \{0\}} \left( \frac12 - \frac1q \right) \left( \frac{\|\nabla u\|_{L^2(A_R)}}{\|u\|_{L^q(A_R)}} \right)^{2q/(q-2)},
\end{align*}
and moreover, $d_S$ can also be characterized similarly in terms of $I_S$. Fix $\varphi \in C^\infty_c(S) \setminus \{0\}$ arbitrarily. Then we observe that $\mathrm{supp}\,\varphi + R e_d \subset A_R$ for $R > 0$ large enough.\prf{\footnote{Indeed, since $\varphi \in C^\infty_c(S)$, we have $\mathrm{supp}\,\varphi \subset B_L^{d-1}(0) \times (\vep, 1 - \vep)$, where $B_L^{d-1}(0) := \{x' \in \R^{d-1} \colon |x'| < L\}$, for some $\vep,L > 0$. Let $x = (x',x_d) \in \mathrm{supp}\,\varphi + R e_d$. Then since $x_d \in (R+\vep, R+1-\vep)$, we have $|x|^2 = |x'|^2 + x_d^2 < L^2 + (R + 1 - \vep)^2 = (R+1)^2 + L^2 + \vep^2 - 2(R+1)\vep < (R+1)^2$ for $R > 0$ large enough. It is obvious that $|x|^2 > R^2$. Hence we obtain $x \in A_R$.}} Define $\phi_R \in H^1_0(A_R)$ by
$$
\phi_R (x) := \begin{cases}
	      \varphi (x-R e_d) &\mbox{ if } \ x \in \mathrm{supp}\,\varphi + R e_d,\\
	      0 &\mbox{ if } \ x \in A_R \setminus (\mathrm{supp}\,\varphi + R e_d)
	     \end{cases} 
$$
for $R > 0$ large enough. Then it follows that
$$
\|\nabla \phi_R \|_{L^2(A_R)}^2 = \|\nabla \varphi\|_{L^2(S)}^2, \quad 
\|\phi_R \|_{L^q(A_R)}^q = \|\varphi\|_{L^q(S)}^q
$$
for $R > 0$ large enough, and therefore, we obtain
\begin{align*}
\limsup_{R \to \infty} d_{A_R} 
&\leq \left(\frac12 - \frac1q\right) \limsup_{R \to \infty} 
\left( \frac{\|\nabla \phi_R \|_{L^2(A_R)}}{\|\phi_R \|_{L^q(A_R)}} \right)^{2q/(q-2)}\\
&= \left(\frac12 - \frac1q\right) \left( \frac{\|\nabla \varphi\|_{L^2(S)}}{\|\varphi\|_{L^q(S)}} \right)^{2q/(q-2)}.
\end{align*}
From the arbitrariness of $\varphi \in C^\infty_c(S) \setminus \{0\}$, we get
$$
\limsup_{R \to \infty} d_{A_R} \leq \left(\frac12 - \frac1q\right) \inf_{\varphi \in H^1_0(S) \setminus \{0\}} \left( \frac{\|\nabla \varphi\|_{L^2(S)}}{\|\varphi\|_{L^q(S)}} \right)^{2q/(q-2)} = d_S.
$$

In order to prove the positivity of $\liminf_{R \to \infty} d_{A_R}$, suppose to the contrary that $d_{A_{R_n}} \to 0$ for some sequence $R_n \to \infty$. Then for each $u_n \in \mathscr{L}_{A_{R_n}}$, we have $\|\nabla u_n\|_{L^2(A_{R_n})} \to 0$. Since $u_n \in H^1_0(A_{R_n}) \subset H^1(\R^d)$ solves \eqref{el-AR-d}, we derive from \eqref{Nehari}, Sobolev's inequality and \eqref{Poincare-c0} that
$$
\|\nabla u_n\|_{L^2(A_{R_n})}^2 = \|u_n\|_{L^q(A_{R_n})}^q = \|u_n\|_{L^q(\R^d)}^q \leq C \|u_n\|_{H^1(\R^d)}^q \leq C(1 + c_0)^{q/2} \|\nabla u_n\|_{L^2(A_{R_n})}^q,
$$
which along with $q > 2$ contradicts to the fact that $\|\nabla u_n\|_{L^2(A_{R_n})} \to 0$. Thus we obtain $\liminf_{R \to \infty} d_{A_R} > 0$. The proof is completed.
\end{proof}

We also have the following corollary:

\begin{corollary}\label{C:d_AR}
It holds that
\begin{equation}\label{rmk-1}
\limsup_{R \to \infty} \sup_{u \in \mathscr{L}_{A_R}} \|u\|_{C^{2,\alpha}(\overline{A_R})} < \infty
\end{equation}
for any $\alpha \in (0,1)$. Moreover, it also holds that
\begin{equation}\label{rmk-2}
\liminf_{R \to \infty} \inf_{u \in \mathscr{L}_{A_R}} \|u\|_{L^\infty(A_R)} > 0.
\end{equation}
\end{corollary}

\begin{proof}
From Lemma \ref{L:d_AR} together with \eqref{Nehari} and \eqref{Poincare-c0}, we infer that $\|u_R\|_{H^1(\R^d)}$ is uniformly bounded as $R \to \infty$. By Moser's iteration, for any $p \in [2,\infty)$, we observe that $\|u_R\|_{L^p(\R^d)}$ is also uniformly bounded as $R \to \infty$. Then \eqref{rmk-1} can be proved by standard elliptic regularity theory (see, e.g.,~\cite{GT}). Here we note that, for any fixed $R_1>0$, $(A_R)_{R\ge R_1}$ forms a uniformly $C^{2,\alpha}$ (in fact, uniformly $C^\infty$) family of domains in the sense that the boundary can be represented in finitely many charts of a fixed radius with uniformly bounded $C^{2,\alpha}$ norms and uniformly bounded overlap. As a consequence, for any $1<p<\infty$, 
the constants in the $L^p$ elliptic estimates (resp., $C^{2,\alpha}$ estimates for $\alpha \in (0,1)$) for $-\Delta$ on $A_R$ depend only on $d,p$ (resp., $d,\alpha$) and the uniform $C^{1,1}$ (resp., $C^{2,\alpha}$) regularity of the domains and therefore can be chosen independently of $R\ge R_1$. 

Concerning \eqref{rmk-2}, suppose to the contrary that there exist a sequence $(R_n)$ in $(0,\infty)$ and a sequence $(u_n)$ in $H^1(\R^d)$ such that $R_n \to \infty$, $u_n \in \mathscr{L}_{A_{R_n}}$ for $n \in \N$ and $\|u_n\|_{L^\infty(A_{R_n})} \to 0$. Then we can derive from \eqref{Nehari}, \eqref{Poincare-c0} and Lemma \ref{L:d_AR} that
\begin{align*}
0 < \liminf_{n \to \infty} \|\nabla u_n\|_{L^2(A_{R_n})}^2 
&= \liminf_{n \to \infty} \|u_n\|_{L^q(A_{R_n})}^q \\
&\leq \liminf_{n \to \infty} \|u_n\|_{L^\infty(A_{R_n})}^{q-2} \|u_n\|_{L^2(A_{R_n})}^2 \\
&\leq c_0 \Big( \lim_{n \to \infty} \|u_n\|_{L^\infty(A_{R_n})}^{q-2} \Big) \Big( \limsup_{n \to \infty} \|\nabla u_n\|_{L^2(A_{R_n})}^2 \Big)\\
&\leq c_0 \Big( \lim_{n \to \infty} \|u_n\|_{L^\infty(A_{R_n})}^{q-2} \Big) \left( \frac12 - \frac1q \right)^{-1} d_S = 0,
\end{align*}
which yields a contradiction. This completes the proof.
\end{proof}

We further prove that
\begin{lemma}\label{L:d_AR-conv}
It holds that $d_{A_R} \to d_S$ as $R \to \infty$. 
\end{lemma}

\begin{proof}
By Lemma \ref{L:d_AR}, it suffices to show that $d_S \leq \liminf_{R \to \infty} d_{A_R}$. To this end, let $(R_n)$ be a sequence in $(0,\infty)$ such that
\begin{equation}\label{set-Rn}
R_n \to \infty \quad \mbox{ and } \quad \lim_{n \to \infty} d_{A_{R_n}} = \liminf_{R \to \infty} d_{A_R}.
\end{equation}
Hence, for each $n \in \N$, we let $u_n \in \mathscr{L}_{A_{R_n}}^\circ$, that is, $I_{A_{R_n}}(u_n) = d_{R_n}$. We may assume that $u_n \geq 0$ for $n \in \N$. We recall that $(u_n)$ is bounded in $H^1(\R^d)$ by Lemma \ref{L:d_AR}. It then follows immediately that, up to a (not relabeled) subsequence, that
\begin{equation}\label{un-conv}
u_n(\,\cdot\, + R_n e_d) \to v_{0,\infty} \quad \mbox{ weakly in } H^1(\R^d) \ \mbox{ and a.e.~in } \R^d
\end{equation}
for some nonnegative function $v_{0,\infty} \in H^1(\R^d)$. Let $x_{1,n} \in A_{R_n}$ be a maximum point of $u_n \in \mathscr{L}_{A_{R_n}}^\circ$ such that $x_{1,n} \in \{r e_d \colon r \in (R_n,R_n+1)\}$. In view of \eqref{rmk-2} of Corollary \ref{C:d_AR}, we may assume that 
\begin{equation}\label{pos-max}
\lim_{n \to \infty} u_n(x_{1,n}) > 0.
\end{equation}
We next claim that 
\begin{equation*}
v_{0,\infty} \in H^1_0(S).
\end{equation*}
Indeed, for every $\varphi \in C^\infty_c(\R^d \setminus \overline{S})$, there exists $N \in \N$ large enough so that $\mathrm{supp} \, \varphi + R_n e_d \subset \R^d \setminus \overline{A_{R_n}}$ for all $n \geq N$.\prf{\footnote{Indeed, since $\varphi \in C^\infty_c(\R^d \setminus \overline S)$, we have $\mathrm{supp}\,\varphi \subset B_L^{d-1}(0) \times (-\vep, 1 + \vep)^c$ for some $\vep,L > 0$. Let $x = (x',x_d) \in \mathrm{supp}\,\varphi + R e_d$. In case $x_d \leq R-\vep$, we have $|x|^2 = |x'|^2 + x_d^2 < L^2 + (R - \vep)^2 = R^2 + L^2 + \vep^2 - 2R\vep < R^2$ for $R > 0$ large enough. In case $x_d \geq R + 1 + \vep$, we find that $|x|^2 > (R+1)^2$. Thus we obtain $x \in \R^d \setminus \overline{A_R}$.}} Hence it follows that
$$
0 = \int_{\R^d} u_n(y + R_n e_d) \varphi(y) \, \d y \to \int_{\R^d} v_{0,\infty}(y) \varphi(y) \, \d y \quad \mbox{ for } \ \varphi \in C^\infty_c(\R^d \setminus \overline{S}),
$$
which gives $v_{0,\infty} = 0$ a.e.~in $\R^d \setminus \overline{S}$. Hence we obtain $v_{0,\infty} \in H^1_0(S)$ (see, e.g.,~\cite[Proposition 9.18]{B-FA}). Moreover, we claim that
\begin{claim}\label{cll1}
The restriction of $v_{0,\infty}$ onto $S$ is a positive solution to \eqref{el-S-d}.  
\end{claim}

Indeed, for any $\varphi \in C^\infty_c(S)$, we find that $\mathrm{supp}\,\varphi + R_n e_d \subset A_{R_n}$ for $R_n$ large enough. Since $u_n$ is a weak solution to \eqref{el-AR-d} and $\varphi (\,\cdot\, - R_n e_d) \in C^\infty_c(A_{R_n})$, we have
$$
\int_{\R^d} \nabla u_n(y+R_n e_d) \cdot \nabla \varphi(y) \, \d y
= \int_{\R^d} u_n^{q-1}(y+R_n e_d) \varphi(y) \, \d y.
$$
Passing to the limit and using \eqref{un-conv}, we infer that $v_{0,\infty} \in H^1_0(S)$ is a nonnegative weak solution to $-\Delta v_{0,\infty} = v_{0,\infty}^{q-1}$ in $S$. 
Furthermore, recalling \eqref{rmk-1}, we can verify that $u_n (\,\cdot\,+R_n e_d) \to v_{0,\infty}$ uniformly on $[0,1] e_d$, and therefore, by \eqref{pos-max}, $v_{0,\infty}$ attains its positive maximum on the line segment $\{r e_d \colon r \in (0,1)\}$. By virtue of the strong maximum principle, we conclude that $v_{0,\infty}$ is a positive solution to \eqref{el-S-d}. Thus Claim \ref{cll1} has been proved. Then we also have 
\begin{equation}\label{dS-ineq}
d_S \leq I_S(v_{0,\infty}). 
\end{equation}

We next set
$$
u_{1,n}(x) := u_n(x) - v_{0,\infty}(x - R_n e_d) \quad \mbox{ for } \ x \in \R^d.
$$
Then we see from \eqref{un-conv} that $u_{1,n}(\,\cdot\, + R_n e_d) \to 0$ weakly in $H^1(\R^d)$. Hence we can derive that
\begin{align*}
 \|\nabla u_n\|_{L^2(A_{R_n})}^2 = \|\nabla u_n\|_{L^2(\R^d)}^2
 &= \|\nabla v_{0,\infty}\|_{L^2(\R^d)}^2 + \|\nabla u_{1,n}\|_{L^2(\R^d)}^2 + o_n(1)\\
 &= \|\nabla v_{0,\infty}\|_{L^2(S)}^2 + \|\nabla u_{1,n}\|_{L^2(\R^d)}^2 + o_n(1).
\end{align*}
Since $\|u_n(\,\cdot\,+R_n e_d)\|_{L^q(\R^d)} = \|u_n\|_{L^q(\R^d)}$ is bounded and $u_n(\,\cdot\,+R_n e_d) \to v_{0,\infty}$ a.e.~in $\R^d$, the Br\'ezis--Lieb lemma implies that
\begin{align*}
 \|u_n\|_{L^q(A_{R_n})}^q = \|u_n\|_{L^q(\R^d)}^q
 &= \|v_{0,\infty}\|_{L^q(\R^d)}^q + \|u_{1,n}\|_{L^q(\R^d)}^q + o_n(1)\\
 &= \|v_{0,\infty}\|_{L^q(S)}^q + \|u_{1,n}\|_{L^q(\R^d)}^q + o_n(1). 
\end{align*}
Therefore we obtain
\begin{equation}\label{d-decomp}
d_{A_{R_n}} = I_{A_{R_n}}(u_n) = I_S(v_{0,\infty}) + I_{\R^d}(u_{1,n}) + o_n(1),
\end{equation}
where $I_{\R^d} : H^1(\R^d) \to \R$ is given by
$$
I_{\R^d}(u) := \frac12 \|\nabla u\|_{L^2(\R^d)}^2 - \frac 1q \|u\|_{L^q(\R^d)}^2 \quad \mbox{ for } \ u \in H^1(\R^d).
$$
Set
\begin{align*}
\delta_1 &:= \sup \left\{ \limsup_{n \to \infty} \|u_{1,n}\|_{L^q(B_{3}(z_n))} \colon z_n \in A_{R_n} \mbox{ for } n \in \N \right\}.\\
&= \limsup_{n \to \infty} \left( \sup_{z \in A_{R_n}} \|u_{1,n}\|_{L^q(B_3(z))} \right).
\end{align*}
Here and henceforth, $B_r(a)$ denotes the open ball in $\R^d$ centered at $a \in \R^d$ of radius $r > 0$. In case $\delta_1 = 0$, we can verify that $\sup_{z \in \R^d} \|u_{1,n}\|_{L^q(B_1(z))} \to 0$ as $n \to \infty$. Indeed, let $z \in \R^d$ be fixed. In case $B_2(z) \cap A_{R_n} = \emptyset$, since $u_n = 0$ in $B_1(z)$, we see that $\|u_{1,n}\|_{L^q(B_1(z))} = \|v_{0,\infty}(\,\cdot\, - R_n e_d)\|_{L^q(B_1(z))}$. If $S \cap B_1(z-R_n e_d) = \emptyset$, we have $v_{0,\infty}(\,\cdot\, - R_n e_d) = 0$ in $B_1(z)$. If $S \cap B_1(z-R_n e_d) \neq \emptyset$, we find that
$$
z \not\in \{ x \in \R^d \colon R_n - 2 < |x| < R_n + 3 \} \quad \mbox{ and } \quad
z - R_n e_d =: (y', y_d) \in \R^{d-1} \times (-1,2).
$$
Then we claim that $|y'|^2 \geq 2R_n+5$. Indeed, suppose that $|z| \leq R_n - 2$. Then we have $R_n - 2 \geq |z| \geq R_n + y_d > R_n - 1$, which yields a contradiction. Hence we obtain $(R_n + 3)^2 \leq |z|^2 = |y'|^2 + (R_n+y_d)^2$, which yields
$$
|y'|^2 \geq (R_n+3)^2 - (R_n+2)^2 = 2R_n + 5.
$$ 
Therefore from the fact that $v_{0,\infty} \in L^q(S)$, it follows that
$$
\|v_{0,\infty}(\,\cdot\, - R_n e_d)\|_{L^q(B_1(z))} \leq \|v_{0,\infty}\|_{L^q(S_n)} \to 0,
$$
where $S_n := \{y = (y',y_d) \in S \colon |y'| \geq \sqrt{2R_n+5} \ \mbox{ and } \ 0 < y_d < 1\}$. In case $B_2(z) \cap A_{R_n} \neq \emptyset$, we have $\|u_{1,n}\|_{L^q(B_1(z))} \leq \sup_{z' \in A_{R_n}} \|u_{1,n}\|_{L^q(B_3(z'))} \to 0$. 
Hence Lions's lemma (see~\cite{Lions84}) yields
\begin{equation}\label{u1n-conv-Lq}
\|u_{1,n}\|_{L^p(\R^d)} \to 0 \ \mbox{ for } \ p \in (2,2^*),
\end{equation}
which particularly implies that
$$
\liminf_{n \to \infty} I_{\R^d}(u_{1,n}) = \liminf_{n \to \infty} \frac12 \|\nabla u_{1,n}\|_{L^2(\R^d)}^2 \geq 0.
$$
Therefore it follows from \eqref{set-Rn}, \eqref{dS-ineq} and \eqref{d-decomp} that
$$
\liminf_{R \to \infty} d_{A_R} = \lim_{n \to \infty} d_{A_{R_n}} = I_S(v_{0,\infty}) + \liminf_{n \to \infty} I_{\R^d}(u_{1,n}) \geq I_S(v_{0,\infty}) \geq d_S.
$$
Consequently, thanks to Lemma \ref{L:d_AR}, we obtain $d_{A_R} \to d_S$ as $R \to \infty$. 

In case $\delta_1 > 0$, extracting a (not relabeled) subsequence of $(n)$, one can take a sequence $(z_n)$ in $\R^d$ such that $z_n \in A_{R_n}$ for each $n \in \N$ and $\lim_{n \to \infty} \|u_{1,n}\|_{L^q(B_{3}(z_n))} \geq \delta_1/2 > 0$. Set
$$
z_{1,n} := R_n \frac{z_n}{|z_n|} \in \partial A_{R_n} \quad \mbox{ for } \ n \in \N.
$$
Then we observe that $\lim_{n \to \infty} \|u_{1,n}\|_{L^q(B_{4}(z_{1,n}))} \geq \delta_1/2 > 0$ (because $B_{3}(z_n) \subset B_{4}(z_{1,n})$) and also claim that
\begin{equation}\label{ort-1}
|z_{1,n} - R_n e_d| \to \infty.
\end{equation}
Indeed, suppose to the contrary that $\rho_0 := \sup_{k \in \N}|z_{1,n_k} - R_{n_k} e_d|$ is finite for some subsequence $(n_k)$ of $(n)$. Then recalling that $u_{1,n} (\,\cdot\, + R_n e_d) \to 0$ in $L^q_{\rm loc}(\R^d)$ (see \eqref{un-conv}), we have 
\begin{align*}
\|u_{1,n_k}\|_{L^q(B_{4}(z_{1,n_k}))} &\leq \|u_{1,n_k}\|_{L^q(B_{\rho_0+{4}}(R_{n_k}e_d ))} = \|u_{1,n_k}(\,\cdot\,+R_{n_k}e_d)\|_{L^q(B_{\rho_0+{4}}(0))} \to 0,
\end{align*}
which is a contradiction. Moreover, we can take $\tilde v_{1,\infty} \in H^1(\R^d)$ such that
\begin{equation}\label{u1n-conv}
u_{1,n}(\,\cdot\,+z_{1,n}) \to \tilde v_{1,\infty} \quad \mbox{ weakly in } H^1(\R^d) \ \mbox{ and a.e.~in } \R^d
\end{equation}
by extracting a (not relabeled) subsequence if necessary. Noting that 
$$
\|\tilde v_{1,\infty}\|_{L^q(B_{4}(0))} = \lim_{n \to \infty} \|u_{1,n}(\,\cdot\, + z_{1,n})\|_{L^q(B_{4}(0))} \geq \delta_1/2 > 0,
$$
we observe that $\tilde v_{1,\infty} \neq 0$. Extracting a (not relabeled) subsequence, we see that $z_{1,n}/R_n$ converges to a unit vector $p_1 \in \R^d$. Let $\Ro_1 \in \mathrm{SO}(d)$ be such that $\Ro_1 e_d = p_1$. Then we can also prove that 
\begin{equation}\label{v1-zero}
v_{1,\infty} := \tilde v_{1,\infty} \circ \Ro_1 \in H^1_0(S).
\end{equation}
Indeed, let $(\Ro_{1,n})$ be a sequence in $\mathrm{SO}(d)$ such that $\Ro_{1,n}e_d = z_{1,n}/R_n$ for $n \in \N$ and $\Ro_{1,n} \to \Ro_1$. For any $\varphi \in C^\infty_c(\R^d \setminus \overline S)$ and $R_n > 0$ large enough, we can find that
\begin{align*}
\MoveEqLeft{
\int_{\R^d} u_{1,n}(\Ro_{1,n} y + z_{1,n}) \varphi(y) \, \d y
}\\
&= \int_{\R^d} u_n(\Ro_{1,n} y + z_{1,n}) \varphi(y) \, \d y
- \int_{\R^d} v_{0,\infty}(\Ro_{1,n} y + z_{1,n} - R_n e_d) \varphi(y) \, \d y\\
&= \int_{\R^d} u_n(x) \varphi(\Ro_{1,n}^{-1} x - R_n e_d) \, \d x
- \int_{\R^d} v_{0,\infty}(\Ro_{1,n} y + z_{1,n} - R_n e_d) \varphi(y) \, \d y\\
&= - \int_{\mathrm{supp}\,\varphi} v_{0,\infty}(\Ro_{1,n} y + z_{1,n} - R_n e_d) \varphi(y) \, \d y
\to 0.
\end{align*}
Here we used the fact that $\mathrm{supp} \, [\varphi(\Phi_{1,n}^{-1}\,\cdot\, - R_n e_d)] = \Phi_{1,n} (\mathrm{supp}\,\varphi + R_n e_d) \subset A_{R_n}^c$ for $R_n$ large enough and \eqref{ort-1} along with the density of $C^\infty_c(\R^d)$ in $L^2(\R^d)$. Moreover, we can observe that
\begin{align}
\MoveEqLeft{
\int_{\R^d} \left( u_{1,n}(\Ro_{1,n}y + z_{1,n}) - \tilde v_{1,\infty}(\Phi_1 y) \right) \varphi(y) \, \d y
}\nonumber\\
&= \int_{\R^d} \left( u_{1,n}(\eta + z_{1,n}) - \tilde v_{1,\infty}(\Phi_1 \circ \Ro_{1,n}^{-1} \eta)\right) \varphi (\Ro_{1,n}^{-1} \eta) \, \d \eta\nonumber\\
&= \int_{\R^d} \left( u_{1,n}(\eta + z_{1,n}) - \tilde v_{1,\infty}(\eta) \right) \varphi (\Ro_{1,n}^{-1} \eta) \, \d \eta\nonumber\\
&\quad + \int_{\R^d} \left( \tilde v_{1,\infty}(\eta) - \tilde v_{1,\infty}(\Phi_1 \circ \Ro_{1,n}^{-1} \eta)\right) \varphi (\Ro_{1,n}^{-1} \eta) \, \d \eta \to 0\label{Ldc}
\end{align}
by employing \eqref{u1n-conv} and Lebesgue's dominated convergence theorem along with the density again. Hence we conclude that
$$
\int_{\R^d} (\tilde v_{1,\infty} \circ \Ro_1) \varphi \, \d y = 0 \quad \mbox{ for } \ \varphi \in C^\infty_c(\R^d \setminus \overline S),
$$
which implies $v_{1,\infty} := \tilde v_{1,\infty} \circ \Ro_1 = 0$ a.e.~in $\R^d \setminus \overline S$. Therefore we obtain \eqref{v1-zero}. 

We next claim that

\begin{claim}\label{cll2}
The restriction of $v_{1,\infty} = \tilde v_{1,\infty} \circ \Ro_1$ onto $S$ is a nontrivial solution to \eqref{el-S-d}.
\end{claim}

Indeed, repeating the same argument as in \eqref{Ldc} with $\varphi \in C^\infty_c(\R^d)$, we derive from \eqref{ort-1} that
\begin{align*}
u_n(\Ro_{1,n} \,\cdot\, + z_{1,n}) 
&= u_{1,n}(\Ro_{1,n} \, \cdot \, + z_{1,n}) + v_{0,\infty}(\Ro_{1,n}\,\cdot\,+z_{1,n} - R_n e_d)\\
&\to \tilde v_{1,\infty} \circ \Ro_1 = v_{1,\infty} \quad \mbox{ weakly in } H^1(\R^d) \ \mbox{ and a.e.~in } \R^d.
\end{align*}
Let $\varphi \in C^\infty_c(S)$. Since $\mathrm{supp} \, [\varphi(\Phi_{1,n}^{-1}\,\cdot\, - R_n e_d)] = \Phi_{1,n} (\mathrm{supp}\,\varphi + R_n e_d) \subset A_{R_n}$ for all $n \geq 1$ large enough, the function $\psi_n(x) := \varphi(\Ro_{1,n}^{-1}x - R_n e_d)$ belongs to $C^\infty_c(A_{R_n})$, and then, it holds that
$$
\int_{A_{R_n}} \nabla u_n(x) \cdot \nabla \psi_n(x) \, \d x = \int_{A_{R_n}} u_n^{q-1}(x) \psi_n(x) \, \d x.
$$
We find that
\begin{align*}
\int_{A_{R_n}} \nabla u_n(x) \cdot \nabla \psi_n(x) \, \d x 
&= \int_{A_{R_n}} \nabla u_n(x) \cdot {}^t (\Ro_{1,n}^{-1}) \nabla \varphi(\Ro_{1,n}^{-1}x - R_n e_d) \, \d x\\
&= \int_{\R^d} \nabla u_n(\Ro_{1,n}y + z_{1,n}) \cdot \Ro_{1,n} \nabla \varphi(y) \, \d y\\
&= \int_{\R^d} {}^t \Ro_{1,n} \nabla u_n(\Ro_{1,n}y + z_{1,n}) \cdot \nabla \varphi(y) \, \d y\\
&= \int_{\R^d} \nabla \left[ u_n(\Ro_{1,n}\,\cdot\, + z_{1,n}) \right](y) \cdot \nabla \varphi(y) \, \d y\\
&\to \int_{\R^d} \nabla v_{1,\infty} \cdot \nabla \varphi \, \d y.
\end{align*}
On the other hand, we observe that
\begin{align*}
\int_{A_{R_n}} u_n^{q-1}(x) \psi_n(x) \, \d x
&= \int_{\R^d} u_n^{q-1}(x) \varphi(\Ro_{1,n}^{-1}x-R_n e_d) \, \d x\\
&= \int_{\R^d} u_n^{q-1}(\Ro_{1,n}y+z_{1,n}) \varphi(y) \, \d y
\to \int_{\R^d} v_{1,\infty}^{q-1} \varphi \, \d y.
\end{align*}
Thus $v_{1,\infty}$ turns out to be a nontrivial (classical) solution to \eqref{el-S-d}. Hence we find that 
$$
I_{\R^d}(\tilde v_{1,\infty}) = I_{\R^d}(v_{1,\infty}) = I_S(v_{1,\infty}) \geq d_S.
$$

We further set
$$
u_{2,n}(x) := u_{1,n}(x) - \tilde v_{1,\infty}(x - z_{1,n}) \quad \mbox{ for } \ x \in \R^d.
$$
Then we have $u_{2,n}(\,\cdot\, + z_{1,n}) \to 0$ weakly in $H^1(\R^d)$. Using \eqref{ort-1} as well as the Br\'ezis--Lieb lemma, we can similarly derive that
$$
d_{A_{R_n}} = I_{A_{R_n}}(u_n) = I_S(v_{0,\infty}) + I_S(v_{1,\infty}) + I_{\R^d}(u_{2,n}) + o_n(1).
$$
Set
$$
\delta_2 := \sup \left\{ \limsup_{n \to \infty} \|u_{2,n}\|_{L^q(B_{3}(z_n))} \colon z_n \in A_{R_n} \mbox{ for } n \in \N \right\}.
$$
In case $\delta_2 = 0$, then one can derive, as before, that $\|u_{2,n}\|_{L^q(\R^d)} \to 0$, which implies that $\liminf_{n \to \infty} I_{\R^d}(u_{2,n}) \geq 0$. Thus we obtain
$$
\liminf_{R \to \infty} d_{A_R} = \lim_{n \to \infty} d_{A_{R_n}} \geq I_S(v_{0,\infty}) + I_S(v_{1,\infty}) \geq 2 d_S > d_S,
$$
which contradicts Lemma \ref{L:d_AR}. In case $\delta_2 > 0$, extracting a (not relabeled) subsequence, one can take a sequence $(z_{2,n})$ in $\R^d$ such that $|z_{2,n}| = R_n$ for each $n \in \N$ and $\lim_{n \to \infty} \|u_{2,n}\|_{L^q(B_{4}(z_{2,n}))} \geq \delta_2/2 > 0$. Then recalling $u_{2,n}(\,\cdot\, + z_{1,n}) \to 0$ and $u_{1,n}(\,\cdot\, + R_n e_d) \to 0$ in $L^q_{\rm loc}(\R^d)$, we can verify, as before, that $|z_{2,n} - z_{1,n}| \to \infty$ and $|z_{2,n} - R_ne_d| \to \infty$, respectively, from the fact that $\delta_2 > 0$.\prf{\footnote{Indeed, if $\rho_0 := \sup_{k \in \N}|z_{2,n_k}-R_{n_k} e_d|$ is finite for some subsequence $(n_k)$ of $(n)$, then it follows that
\begin{align*}
\|u_{2,{n_k}}\|_{L^q(B_{4}(z_{2,{n_k}}))} & \leq \|u_{2,{n_k}}\|_{L^q(B_{\rho_0 + {4}}(R_{n_k} e_d))}
= \|u_{1,{n_k}} - \tilde v_{1,\infty}(\,\cdot\,-z_{1,{n_k}})\|_{L^q(B_{\rho_0 + {4}}(R_{n_k} e_d))}\\
&\leq \|u_{1,{n_k}}(\,\cdot\,+R_{n_k} e_d)\|_{L^q(B_{\rho_0+{4}}(0))} + \|\tilde v_{1,\infty}(\,\cdot\,-z_{1,{n_k}}+R_{n_k} e_d)\|_{L^q(B_{\rho_0+{4}}(0))} \to 0,
\end{align*}
which yields a contradiction to $\limsup_{n \to \infty} \|u_{2,n}\|_{L^q(B_{4}(z_{2,n}))} > 0$. Thus we obtain $|z_{2,n}-R_n e_d| \to \infty$. Moreover, if $\rho_1 := \sup_{k \in \N}|z_{2,n_k}-z_{1,n_k}|$ is finite for some subsequence $(n_k)$ of $(n)$, then we have
\begin{align*}
\|u_{2,{n_k}}\|_{L^q(B_{4}(z_{2,{n_k}}))} & \leq \|u_{2,{n_k}}\|_{L^q(B_{\rho_1 + {4}}(z_{1,n_k}))}
= \| u_{2,{n_k}}(\,\cdot\,+ z_{1,n_k}) \|_{L^q(B_{\rho_1 + {4}}(0))} \to 0,
\end{align*}
which yields a contradiction to $\limsup_{n \to \infty} \|u_{2,n}\|_{L^q(B_{4}(z_{2,n}))} > 0$. Thus we obtain $|z_{2,n}-z_{1,n}| \to \infty$.}} Moreover, repeating the same argument as in the proof of Claim \ref{cll2}, we can take a limit $\tilde v_{2,\infty} \in H^1(\R^d) \setminus \{0\}$ such that $u_{2,n}(\,\cdot\,+z_{2,n}) \to \tilde v_{2,\infty}$ and $v_{2,\infty} := \tilde v_{2,\infty} \circ \Ro_2 \in H^1_0(S)$ turns to be a nontrivial classical solution to \eqref{el-S-d} for some rotation $\Ro_2 \in \mathrm{SO}(d)$. Then it also follows that
$$
d_{A_{R_n}} = \sum_{k=0}^2 I_S(v_{k,\infty}) + I_{\R^d}(u_{3,n}) + o_n(1),
$$
where $u_{3,n}(x) := u_{2,n}(x) - \tilde v_{2,\infty}(x - z_{2,n})$ for $x \in \R^d$. Set
$$
\delta_3 := \sup \left\{ \limsup_{n \to \infty} \|u_{3,n}\|_{L^q(B_{3}(z_n))} \colon z_n \in A_{R_n} \mbox{ for } n \in \N \right\}.
$$
The argument again leads to the following dichotomy: $\delta_3 = 0$ or $\delta_3 > 0$. The former case leads to a contradiction, as before, whereas in the latter case we proceed to the next step by choosing $\delta_4$ in a similar fashion.

To conclude the proof, it is sufficient to prove that this procedure cannot be repeated infinitely many times. Suppose that $\delta_k > 0$ for $k = 1,\ldots,m$ for some $m \geq 3$. Then
$$
C \geq \|\nabla u_n\|_{L^2(A_{R_n})}^2 = \|\nabla u_n\|_{L^2(\R^d)}^2 = \sum_{k=0}^m \|\nabla v_{k,\infty}\|_{L^2(\R^d)}^2 + \|\nabla u_{m+1,n}\|_{L^2(\R^d)}^2 + o_n(1),
$$
where $v_{k,\infty}$ is a nontrivial solution to \eqref{el-S-d} for $k = 1,\ldots,m$. Thanks to Sobolev--Poincar\'e's inequality, there exists a constant $C_{S} = C_{S}(d,q) > 0$ such that
$$
\|\nabla v_{k,\infty}\|_{L^2(S)}^2 \stackrel{\eqref{el-S-d}}= \|v_{k,\infty}\|_{L^q(S)}^q \leq C_S^q \|\nabla v_{k,\infty}\|_{L^2(S)}^q,
$$
which along with the fact that $v_{k,\infty} \neq 0$ implies that $\|\nabla v_{k,\infty}\|_{L^2(S)} \geq C_S^{-q/(q-2)} > 0$ for any $k = 1,\ldots,m$. Since $\|\nabla u_n\|_{L^2(\R^d)}$ is bounded, we deduce that such an $m$ is bounded from above, and hence we must reach the case $\delta_{m+1} = 0$ after finitely many repetitions. Then we obtain
$$
\liminf_{R \to \infty} d_{A_R} = \lim_{n \to \infty} d_{A_{R_n}} \geq \sum_{k=0}^{m+1} I_S(v_{k,\infty}) \geq (m+2) d_S > d_S,
$$
which contradicts Lemma \ref{L:d_AR}. Therefore, we eventually conclude that $\delta_1 = 0$, and hence $d_{A_R} \to d_S$ as $R \to \infty$. This completes the proof.
\end{proof}

Then we also have:

\begin{corollary}\label{C:conv}
Let $(R_n)$ be a sequence in $(0,\infty)$ such that $R_n \to \infty$. Let $(u_n)$ be a sequence in $H^1(\R^d)$ such that $u_n \in \mathscr{L}^\circ_{A_{R_n}}$ for $n \in \N$. Then there exist a subsequence $(n_k)$ of $(n)$ and $v_{0,\infty} \in \mathscr{L}^\circ_S$ such that
\begin{gather}
u_{n_k}(\,\cdot\,+R_{n_k} e_d) \to v_{0,\infty} \quad \mbox{ weakly in } H^1(\R^d) \ \mbox{ and a.e.~in } \R^d,\label{c:as1}\\
\| u_{n_k} - v_{0,\infty} ( \,\cdot\, - R_{n_k} e_d) \|_{L^p(\R^d)} = \| u_{1,n_k} \|_{L^p(\R^d)} \to 0 \quad \mbox{ for } \ p \in (2,\infty],\label{c:as2}\\ 
I_S (v_{0,\infty}) = d_S \quad \mbox{ and } \quad \| \nabla u_{n_k} - \nabla v_{0,\infty} (\,\cdot\, - R_{n_k} e_d) \|_{L^2(\R^d)} \to 0.\label{c:as3}
\end{gather}
\end{corollary}

\begin{proof}
Recall the proof of Lemma \ref{L:d_AR-conv}, where we concluded that $\delta_1 = 0$. Then \eqref{c:as1} coincides with \eqref{un-conv} for some $v_{0,\infty} \in H^1(\R^d)$, which was also shown to be a positive solution to \eqref{el-S-d}. Furthermore, \eqref{c:as2} for $p \in (2,2^*)$ follows immediately from \eqref{u1n-conv-Lq}. Thanks to \eqref{rmk-1}, we see that $\sup_n \|u_n\|_{W^{1,\infty}(A_{R_n})} < \infty$. By the Gagliardo--Nirenberg inequality (see, e.g.,~\cite[Theorem 12.83]{Leoni}), we deduce that
\begin{align*}
\lefteqn{
\| u_n - v_{0,\infty } (\,\cdot\, - R_n e_d) \|_{L^\infty(\R^d)}
}\\
&\leq c_{d,q} \| u_n - v_{0,\infty } (\,\cdot\, - R_n e_d) \|_{L^q(\R^d)}^{q/(d+q)} \| \nabla u_n - \nabla v_{0,\infty } (\,\cdot\, - R_n e_d) \|_{L^{\infty}(\R^d)}^{d/(d+q)}
\end{align*}
for some constant $c_{d,q} > 0$ depending only on $d,q$. Therefore we can conclude that $\| u_n - v_{0,\infty } (\,\cdot\, - R_n e_d) \|_{L^\infty(\R^d)} \to 0$. Hence \eqref{c:as2} holds for any $p \in (2,\infty]$.

By Lemma \ref{L:d_AR-conv}, recalling \eqref{d-decomp}, we infer that
$$
d_S = \lim_{n \to \infty} d_{A_{R_n}} = I_S(v_{0,\infty}) + \liminf_{n \to \infty} I_{\R^d} (u_{1,n}) \geq I_S(v_{0,\infty}) \geq d_S,
$$
which yields $I_S(v_{0,\infty}) = d_S$. Hence we obtain $v_{0,\infty} \in \mathscr{L}_{S}^\circ$. Furthermore, we observe that
\begin{align*}
\|\nabla u_n\|_{L^2(\R^d)}^2 = \|\nabla u_n\|_{L^2(A_{R_n})}^2 \stackrel{\eqref{Nehari}}= \|u_n\|_{L^q(A_{R_n})}^q = \|u_n\|_{L^q(\R^d)}^q \to \|v_{0,\infty}\|_{L^q(\R^d)}^q \stackrel{\eqref{Nehari}}= \|\nabla v_{0,\infty}\|_{L^2(\R^d)}^2,
\end{align*}
which along with the weak convergence \eqref{un-conv} \prf{and the uniform convexity of $\|\cdot\|_{L^2(\R^d)}$ }implies \eqref{c:as3}. This completes the proof.
\end{proof}

Moreover, we give

\begin{lemma}\label{L:decay}
Let $(R_n)$ be a sequence in $(0,\infty)$ such that $R_n \to \infty$. Let $(u_n)$ be a sequence in $H^1(\R^d)$ such that
\begin{equation}
u_n \in \mathscr{L}^\circ_{A_{R_n}}, \ u_n \geq 0 \ \mbox{ for } \ n \in \N \quad \mbox{ and } \quad 
\lim_{n \to \infty}\left\| u_n - v_\infty(\,\cdot\, - R_n e_d) \right\|_{L^\infty( \R^d )} = 0\label{dec-hypo}
\end{equation}
for some $v_\infty \in \mathscr{L}^\circ_S$. Then there exist constants $\lambda,C_1 > 0$ such that
\begin{equation}\label{barrier}
0 < u_n (x) \leq C_1 \exp \left(- \lambda|x - R_n e_d| \right) \quad \mbox{ for } \ x \in A_{R_n}
\end{equation}
for $n \in \N$ large enough. Furthermore, it also holds that
\begin{equation}\label{barrier-v}
0 < v_\infty(y) \leq C_1 \exp \left( - \lambda |y| \right) \quad \mbox{ for } \ y \in S.
\end{equation}
\end{lemma}

\begin{proof}
Notice that
$$
\sup_{|x-R_n e_d| \geq L} u_n (x) \leq \sup_{|y| \geq L} v_\infty(y) + \left\| u_n - v_\infty(\,\cdot\, - R_n e_d) \right\|_{L^\infty( \R^d \setminus B_L(R_ne_d) )} 
$$
for each $L \geq 1$ and $n \geq 1$. Since $v_\infty \in H^1_0(S)$, one can verify that $\sup_{|y| \geq L} v_\infty (y) \to 0$ as $L \to \infty$ by Moser's iteration\prf{\ (see Appendix \S \ref{apdx:comm:2})}. Therefore since \eqref{dec-hypo} holds and $\lambda_{R_n} \to \pi^2$ (see Lemma \ref{L:LapEV}), there exist $L_0 \in \N$ and $n_0 \in \N$ such that 
\begin{equation}\label{e:sizeuRn}
\sup_{|x-R_n e_d| \geq L_0  } u_n^{q-2}(x) \leq \frac{\pi^2}{27} < \frac{\lambda_{R_n}}{2}
\end{equation}
for any $n \geq n_0$. Next, choose $\lambda \in (0,1)$ such that
\begin{equation}\label{cho-lam}
\frac{2\pi}{3\sqrt 3} \lambda + \lambda^2 \leq \frac{\pi^2}{27}
\end{equation}
and define a function $w_n$ by 
$$
w_n(x) := C_0 w_{n,1} (|x|) w_{n,2} (|x-R_n e_d|) 
:= C_0 \sin \left( \frac{\pi}{3} \left( |x| - \left( R_n - 1 \right) \right) \right) \exp \left( - \lambda \left| x - R_n e_d \right| \right).
$$
Here we choose a positive constant $C_0$ such that
$$
\left(\frac{\pi^2}{27}\right)^{1/(q-2)} \leq C_0 \sin \left( \frac{\pi}3 \right) \e^{-\lambda L_0},
$$
which along with the fact that $1 < |x| - (R_n - 1) < 2$ for $x \in A_{R_n}$ implies
\begin{align*}
u_n (x) \leq C_0 \sin \left( \frac{\pi}{3} \left( |x| - \left( R_n - 1 \right) \right) \right) \e^{ - \lambda L_0 } = w_n(x)
\end{align*}
for each $n \geq n_0$ and $x \in A_{R_n}$ satisfying $|x-R_n e_d|=L_0$. It is easily seen that $w_n > 0$ in $\{x \in \R^d \colon R_n-1 < |x| < R_n+2\}$, and hence, we find that $w_n > 0 = u_n$ on $\partial A_{R_n}$. Moreover, we have
\begin{align*}
-\Delta w_n &= (-\Delta w_{n,1}) w_{n,2} - 2 \nabla w_{n,1} \cdot \nabla w_{n,2} + w_{n,1} (-\Delta w_{n,2})\\
&= C_0 \left[  \frac{\pi^2}{9} \sin \left( \frac{\pi}{3} \left( |x| - \left( R_n - 1 \right) \right) \right) 
- \frac{d-1}{|x|} \frac{\pi}{3} \cos \left( \frac{\pi}{3} \left( |x| - \left( R_n - 1 \right) \right) \right) \right] w_{n,2} \\
& \quad +2 C_0 \frac{\pi}{3} \lambda \cos \left( \frac{\pi}{3} \left( |x| - \left( R_n - 1 \right) \right) \right) \frac{x}{|x|} \cdot \frac{x-R_n e_d}{|x-R_n e_d|}w_{n,2} \\
&\quad + C_0 \sin \left( \frac{\pi}{3} \left( |x| - \left( R_n - 1 \right) \right) \right) \left( -\lambda^2 + \frac{\lambda (d-1) }{|x-R_ne_d|} \right) w_{n,2}
\end{align*}
for $x \in A_{R_n} \setminus B_{L_0}(R_ne_d)$. Noting that
$$
\left| \cos \left( \frac{\pi}{3} \left( |x| - \left( R_n - 1 \right) \right) \right)  \right| \leq \frac1{\sqrt 3} \sin \left( \frac{\pi}{3} \left( |x| - \left( R_n - 1 \right) \right) \right)
$$
for all $x \in A_{R_n}$, i.e., $1 < |x| - (R_n - 1) < 2$, we derive from \eqref{e:sizeuRn} and \eqref{cho-lam} that
\begin{equation*}
-\Delta w_n - u_n^{q-2}(x) w_n 
 \geq \left( \frac{\pi^2}{9} - \frac{\pi}{3\sqrt 3} \frac{d-1}{R_n} - \frac{2\pi}{3\sqrt 3} \lambda - \lambda^2 - u_n^{q-2} (x) \right) w_n
 \geq \left( \frac{\pi^2}{27} - \frac{\pi}{3\sqrt 3} \frac{d-1}{R_n} \right) w_n
\end{equation*}
for every $n \geq n_0$ and $x \in A_{R_n} \setminus B_{L_0}(R_n e_d)$. We may suppose that $\pi^2/27 - \pi (d-1) / (3\sqrt 3R_n) \geq 0$ for every $n \geq n_0$ large enough. Then we obtain 
\begin{equation}\label{wn-sup}
-\Delta w_n - u_n^{q-2} w_n \geq 0 = -\Delta u_n - u_n^{q-2} u_n \quad \text{in \ $A_{R_n} \setminus B_{L_0} (R_n e_d)$}. 
\end{equation}
 Recalling that $w_n \geq u_n$ on $\partial (A_{R_n} \setminus B_{L_0}(R_n e_d))$ and $( u_n - w_n )_+ \in H^1_0(A_{R_n}) \cap C(\overline{A_{R_n}})$, we can deduce that $( u_n - w_n )_+ \bm{1}_{ A_{R_n} \setminus B_{L_0} (R_n e_d) } \in H^1_0(A_{R_n})$ and $\nabla [( u_n - w_n )_+ \bm{1}_{ A_{R_n} \setminus B_{L_0} (R_n e_d) }] = \bm{1}_{ A_{R_n} \setminus B_{L_0} (R_n e_d) } \nabla ( u_n - w_n )_+$, where $\bm{1}_{A}$ stands for the characteristic function of $A \subset \R^d$ and $(s)_+ := \max \{0,s\} \geq 0$ for $s \in \R$.\prf{\footnote{Indeed, we find that $(u_n - w_n)_+|_{A_{R_n} \setminus B_{L_0}(R_n e_d)} \in H^1_0(A_{R_n} \setminus B_{L_0}(R_n e_d))$. Hence its zero extension $h_n := (u_n - w_n)_+ \bm{1}_{ A_{R_n} \setminus B_{L_0} (R_n e_d) }$ onto $A_{R_n}$ belongs to $H^1_0(A_{R_n})$, and moreover, we have $\nabla h_n = \bm{1}_{ A_{R_n} \setminus B_{L_0} (R_n e_d)} \nabla (u_n - w_n)_+ $. We also observe that $\nabla (u_n - w_n) \cdot \nabla h_n = \bm{1}_{ A_{R_n} \setminus B_{L_0}(R_n e_d)} |\nabla (u_n - w_n)_+|^2 = |\nabla h_n|^2$.}} Hence testing \eqref{wn-sup} by $( u_n - w_n )_+ \bm{1}_{ A_{R_n} \setminus B_{L_0} (R_n e_d) } \geq 0$, we see that 
\begin{align*}
\left\| \nabla \left[ (u_n - w_n)_+ \bm{1}_{ A_{R_n} \setminus B_{L_0} (R_n e_d) } \right] \right\|_{L^2(A_{R_n})}^2 - \int_{ A_{R_n} \setminus B_{L_0} (R_n e_d) } u_n^{q-2} \left( u_n - w_n \right)^2_+ \d x \leq 0.
\end{align*}
From the Poincar\'e inequality in $A_{R_n}$ and \eqref{e:sizeuRn}, it follows that
\begin{align*}
\lambda_{R_n} \left\| (u_n - w_n)_+ \bm{1}_{ A_{R_n} \setminus B_{L_0} (R_n e_d) } \right\|_{L^2(A_{R_n})}^2 
&\leq \left\| \nabla \left[ (u_n - w_n)_+ \bm{1}_{ A_{R_n} \setminus B_{L_0} (R_n e_d) } \right]   \right\|_{L^2(A_{R_n})}^2 \\
&\leq  \int_{ A_{R_n} \setminus B_{L_0} (R_n e_d) } u_n^{q-2} \left( u_n - w_n \right)^2_+ \d x \\
&\leq \frac{\lambda_{R_n}}{2} \left\| (u_n - w_n)_+ \bm{1}_{ A_{R_n} \setminus B_{L_0} (R_n e_d) } \right\|_{L^2(A_{R_n})}^2.
\end{align*}
Thus we obtain $(u_n - w_n)_+ \equiv 0$, that is, $u_n \leq w_n$, in $A_{R_n} \setminus B_{L_0} (R_n e_d)$. 

Since $u_n (\,\cdot\,+R_ne_d) \to v_\infty$ in $L^\infty( B_{L_0} (0) )$, we can take a constant $C_1>0$ such that 
$$
u_n(x) \leq C_1 \sin \left( \frac{\pi}{3} \left( |x| - \left( R_n - 1 \right) \right) \right) \exp \left( - \lambda \left| x - R_n e_d \right| \right)
$$
for any $n \geq n_0$ and $x \in A_{R_n}$. Thus \eqref{barrier} follows. Letting $y \in S$, we see that $y + R_n e_d \in A_{R_n}$ for $n \in \N$ large enough. Hence substituting $y + R_n e_d$ to $x$ in \eqref{barrier} and passing to the limit as $n \to \infty$, since $u_n(\,\cdot\, + R_n e_d) \to v_\infty$ uniformly in $\R^d$ as $n \to \infty$ (see \eqref{dec-hypo}), we also obtain \eqref{barrier-v} (cf.~\cite[Proposition A.1.]{DaYa02}). This completes the proof.
\end{proof}

Now, we are ready to prove Proposition \ref{P:les-d}. 

\begin{proof}[Proof of Proposition {\rm \ref{P:les-d}}]
We argue by contradiction. Suppose to the contrary that there exist sequences $(R_n)$ and $(u_n)$ such that 
$$
R_n \to \infty, \quad u_n \in \mathscr{L}_{R_n}^\circ \ \mbox{ for } \ n \in \N, \quad  
c_0 := \liminf_{n \to \infty} \inf_{v \in \mathscr{L}_S^\circ} \left\| u_n -  v(\,\cdot\, - R_n e_d) \right\|_{H^1(\R^d)} > 0.
$$
Thanks to Corollary \ref{C:conv}, there exists $v_\infty \in \mathscr{L}_{S}^\circ$ such that, up to a (not relabeled) subsequence,
\begin{equation}\label{e:conv-1st}
\left\| \nabla u_n - \nabla v_\infty(\,\cdot\, - R_n e_d) \right\|_{L^2(\R^d)} + \left\| u_n - v_\infty(\,\cdot\, - R_n e_d) \right\|_{L^p(\R^d)} \to 0
\end{equation}
for any $p \in (2,\infty]$. Hence we find that 
\begin{equation}\label{hypo-P}
\liminf_{n \to \infty} \| u_n - v_\infty (\,\cdot\, - R_n e_d) \|_{L^2(\R^d)} \geq c_0 > 0.
\end{equation} 
On the other hand, thanks to Lemma \ref{L:decay} along with \eqref{e:conv-1st} with $p = \infty$, it follows that 
$$
\| u_n - v_\infty (\,\cdot\, - R_n e_d) \|_{L^2(\R^d)} \leq C\| u_n - v_\infty (\,\cdot\, - R_n e_d) \|_{L^\infty(\R^d)}^{1/2} \to 0,
$$
which contradicts \eqref{hypo-P}. This completes the proof. 
\end{proof}

Let $p,q \in \R^d$ be such that $|p| = |q| = 1$ and $p \cdot q = 0$. Let $\Ro_{p,q}(\theta) \in \mathrm{SO}(d)$ be the rotation on the $(p,q)$-plane with angle $\theta \in \R$, that is, $\Ro_{p,q}(\theta)x = \e^{\theta (q \otimes p - p \otimes q)} x$ for $x \in \R^d$. Furthermore, the rotational derivative $\Omega_{p,q}w$ of a function $w \in C^1(\R^d)$ on the $(p,q)$-plane is defined by
\begin{equation*}
\Omega_{p,q} w(x) = \lim_{\theta \to 0} \frac{w(\Ro_{p,q}^{-1}(\theta)x) - w(x)}\theta = \left( x_q \partial_{x_p} - x_p \partial_{x_q} \right) w(x) \quad \mbox{ for } \ x \in \R^d,
\end{equation*}
where $x_p := x \cdot p$ and $\partial_{x_p} := p \cdot \nabla$. In what follows, we simply denote $\Ro_i = \Ro_{e_i,e_d}$ and $\Omega_i = \Omega_{e_i,e_d}$ for $i = 1,\ldots,d-1$, where $\{e_j\}_{j=1}^d$ stands for the canonical basis of $\R^d$. For each weak solution $u \in H^1_0(A_R)$ of \eqref{el-AR-d}, denote by $L_{u} = - \Delta - (q-1)|u|^{q-2}$ the linearized operator at $u$, viewed as an operator from $H^1_0(A_R)$ to $H^{-1}(A_R)$. 

We finally prove that

\begin{lemma}[Kernel of the linearized operator at $u_R$]\label{L:les2}
Let $d \geq 2$ and suppose that {\rm (S)} holds. Then there exists $R_0 > 0$ such that, for any $R \geq R_0$ and $u_R \in \mathscr{L}_{A_R}^\circ$, it holds that
$$
\ker L_{u_R} = \mathrm{span} \{\Omega_1 u_R, \ldots, \Omega_{d-1} u_R\},
$$
and hence $\dim \ker L_{u_R} = d-1$. 
%
\end{lemma}

\begin{proof}
We let $u_R \in \mathscr{L}_{A_R}^\circ$ be nonradial for $R > 0$ large enough. We may assume that $u_R > 0$ in $A_R$ for $R > 0$. It is known that any least-energy solution to \eqref{el-AR-d} is foliated Schwarz symmetric (see~\cite[Theorem 3.1]{Pac02},~\cite[Theorem 6.14]{DaPa} and also Lemma \ref{L:sym} in \S \ref{S:opt}). Hence there exists a function $\tilde u_R : [R,R+1] \times [0,\pi] \to \R$ such that $u_R(x) = \tilde{u}_R(r,s)$ with $r = |x|$ and $s = \arccos(e_d \cdot x/|x|)$ for $x \in A_R$ and $\tilde u_R$ is strictly decreasing in $s$, that is, $\partial_s \tilde u_R < 0$ (see~\cite[Remark 6.13]{DaPa}). Since $A_R$ is radial, the rotational derivative $w_i := \Omega_i u_R \in H^1_0(A_R)$ is a weak solution to the Dirichlet problem,
$$
L_{u_R} w_i = 0 \ \mbox{ in } A_R, \quad w_i = 0 \ \mbox{ on } \partial A_R
$$
for $i = 1,\ldots,d-1$ (see~\cite[Lemma 6.11]{DaPa}). We claim that $\{\Omega_1 u_R, \ldots, \Omega_{d-1} u_R\}$ is linearly independent. Indeed, let $c = \{c_i\}_{i=1}^{d-1} \in \R^{d-1}$ be such that $c_1 \Omega_1 u_R(x) + \cdots c_{d-1} \Omega_{d-1} u_R(x) = 0$ for $x \in A_R$. Then setting $x_d = 0$, we have
$$
(c_1 x_1 + \cdots + c_{d-1} x_{d-1}) \partial_{x_d} u_R(x) = 0 \quad \mbox{ for } \ x = (x_1,\ldots,x_{d-1},0) \in A_R. 
$$
Since $\partial_{x_d} u_R(x) = - \partial_s \tilde u_R(r,s)/r > 0$ if $x_d = 0$, we obtain $c = 0$. Thus $\{\Omega_1 u_R, \ldots, \Omega_{d-1} u_R\}$ turns out to be linearly independent, and hence $\dim \ker L_{u_R} \geq d-1$. 

Hence we suppose to the contrary that there exist sequences $(R_n)$ in $(0,\infty)$ and $(u_n)$ in $H^1(\R^d)$ such that $R_n \to \infty$ and $u_n \in \mathscr{L}_{A_{R_n}}^\circ$ and $\dim \ker L_{u_n} \geq d$ for each $n \in \N$. We may assume $u_n \geq 0$. Thanks to Corollary \ref{C:conv}, we may also assume that 
$$
u_n (\cdot+R_n e_d) \to v_\infty \ \mbox{ in } L^p(\R^d) \ \mbox{ for any } p \in (2,\infty]
$$
for some $v_\infty \in \mathscr{L}_S^\circ$. Moreover, for every $n \in \N$, one can take $d$ solutions $e^j_{n} \in H^1_0(A_{R_n}) \setminus \{0\}$ with $j = 1,\ldots,d$ for the linearized problem,
\begin{equation*}
- \Delta e^j_{n} = (q-1) u_n^{q-2}  e^j_{n} \ \mbox{ in } A_{R_n}, \quad e^j_{n} = 0 \ \mbox{ on } \partial A_{R_n}, \quad \|e^j_{n}\|_{L^q(A_{R_n})} = 1
\end{equation*}
such that
\begin{equation}\label{orth-R}
\int_{A_{R_n}} \nabla e^j_{n} \cdot \nabla e^k_{n} \, \d x = (q-1)\int_{A_{R_n}} e^j_{n} e^k_{n} u_n^{q-2} \, \d x = 0 \quad \mbox{ if } \ j \neq k
\end{equation}
for $j,k \in \{1,\ldots,d\}$. Since $(u_n)$ is bounded in $L^q(\R^d)$, it follows from \eqref{Poincare-c0} that
\begin{equation}\label{eR-est}
\|e^j_{n}\|_{H^1(\R^d)} \leq C \quad \mbox{ for } \ n \in \N \ \mbox{ and } \ j \in \{1,\ldots,d\}.
\end{equation}
For each $j \in \{1,\ldots, d\}$, we set
$$
\delta_j := \sup \left\{ \limsup_{n \to \infty} \|e^j_{n}\|_{L^q(B_2(z_n))} \colon z_n \in A_{R_n} \mbox{ for } n \in \N \right\}.
$$
If $\delta_j = 0$, we have $\sup_{z \in \R^d} \|e^j_{n}\|_{L^q(B_1(z))} \to 0$. By virtue of Lions's lemma, we obtain $e^j_{n} \to 0$ in $L^q(\R^d)$. This is a contradiction to $\|e^j_{n}\|_{L^q(A_{R_n})} = 1$. Hence we deduce that $\delta_j > 0$. Therefore there exists a sequence $(x^j_{n})$ in $\R^d$ such that $x^j_{n} \in A_{R_n}$ for $n \in \N$ and $\lim_{n \to \infty} \|e^j_{n}\|_{L^q(B_2(x^j_{n}))} \geq \delta_j/2 > 0$ after extracting a subsequence of $(n)$, which we do not relabel.

We claim that there exists $L_0 > 0$ such that 
\begin{equation}\label{2}
|x^j_{n} - R_n e_d| < L_0
\quad \mbox{ for all } \ j = 1,\ldots,d
\end{equation}
for $n \in \N$ large enough. Suppose to the contrary that there exist a subsequence $(n_k)$ of $(n)$ and a sequence $(j_k)$ in $\{1,\ldots,d\}$ such that
$$
|x^{j_k}_{n_k} - R_{n_k} e_d| \geq 2k \quad \mbox{ for } \ k \in \N.
$$
We can assume $j_k \equiv j \in \{1,\ldots,d\}$ for all $k \in \N$ by taking a further subsequence of $(n_k)$ if necessary. For each $k \in \N$, let $\Ro^j_{n_k} \in \mathrm{SO}(d)$ be such that $\Ro^j_{n_k} e_d = x^j_{n_k}/|x^j_{n_k}|$. Then we observe that
\begin{equation}\label{x-apart}
|\Ro^j_{n_k} x - R_{n_k} e_d| \geq |x^j_{n_k} - R_{n_k} e_d| - |\Ro^j_{n_k} x - x^j_{n_k}| \geq 2k - \left| x - |x^j_{n_k}| e_d \right| \geq k - 1
\end{equation}
for all $x \in B_{k}(R_{n_k} e_d)$. Here we used the fact that $|x - |x^j_{n_k}| e_d| \leq |x - R_{n_k}e_d| + |R_{n_k} - |x^j_{n_k}|| \leq k + 1$ for $x \in B_{k}(R_{n_k} e_d)$. Setting $\hat e^j_{n_k} := e^j_{n_k} \circ \Ro^j_{n_k}$, we have
$$
- \Delta \hat e^j_{R_{n_k}} = (q-1)(u_{n_k} \circ \Ro^j_{n_k})^{q-2} \hat e^j_{R_{n_k}} \ \mbox{ in } A_{R_{n_k}}, \quad \hat e^j_{R_{n_k}} = 0 \ \mbox{ on } \partial A_{R_{n_k}}.
$$
We derive from \eqref{barrier} and \eqref{x-apart} that
\begin{equation}\label{uR-decay-k}
\sup \{ (u_{n_k} \circ \Ro^j_{n_k})(x) \colon x \in A_{n_k} \cap B_k(R_{n_k} e_d) \} \leq C_1 \e^{-\lambda (k-1)}.
\end{equation}
Using \eqref{eR-est}, we can deduce that
$$
\hat f^j_{n_k} := \hat e^j_{{n_k}} (\,\cdot\,+R_{n_k} e_d) \to \hat f^j \ \mbox{ weakly in } H^1(\R^d) \ \mbox{ and a.e.~in } \R^d,
$$
and moreover, thanks to \eqref{uR-decay-k}, arguing as in Claim \ref{cll1} in the proof of Lemma \ref{L:d_AR-conv}, we can prove that $\hat f^j \in H^1_0(S)$ and $\hat f^j$ is a weak solution to
$$
- \Delta \hat f^j = 0 \ \mbox{ in } S, \quad \hat f^j = 0 \ \mbox{ on } \partial S,
$$
which yields $\hat f^j \equiv 0$. On the other hand, noting that 
\begin{align*}
\|\hat f^j_{{n_k}}\|_{L^q(B_3(0))}
&= \|\hat e^j_{n_k}(\,\cdot\,+R_{n_k} e_d)\|_{L^q(B_3(0))}\\ 
&= \|e^j_{n_k} \circ \Ro^j_{n_k} (\,\cdot\,+R_{n_k} e_d) \|_{L^q(B_3(0))}\\ 
&= \|e^j_{n_k} (\,\cdot\,+R_{n_k} x^j_{n_k}/|x^j_{n_k}|) \|_{L^q(B_3(0))}\\ 
&= \|e^j_{n_k} \|_{L^q(B_3(R_{n_k} x^j_{n_k}/|x^j_{n_k}|))} \geq \|e^j_{n_k} \|_{L^q(B_2(x^j_{n_k}))}
\end{align*}
and $\hat f^j_{{n_k}} \to \hat f^j$ in $L^q(B_3(0))$ and recalling that $\lim_{n \to \infty} \|e^j_{n}\|_{L^q(B_2(x^j_{n}))} \geq \delta_j/2 > 0$, we conclude that $\|\hat f^j\|_{L^q(B_3(0))} \geq \delta_j /2 > 0$, which yields a contradiction. Thus \eqref{2} follows.

Let $j \in \{1,\ldots,d\}$ be fixed. Using \eqref{eR-est}, we have
\begin{equation}\label{c:e2}
f^j_{n} := e^j_{{n}}(\,\cdot\,+{R_n} e_d) \to f^j \quad \mbox{ weakly in } H^1(\R^d) \ \mbox{ and a.e.~in } \R^d.
\end{equation}
Arguing as in Claim \ref{cll1} in the proof of Lemma \ref{L:d_AR-conv}, we can verify that $f^j \in H^1_0(S)$ and $f^j$ is a weak solution to
$$
- \Delta f^j = (q-1) v^{q-2}_\infty f^j \ \mbox{ in } S, \quad f^j = 0 \ \mbox{ on } \partial S.
$$
Furthermore, setting $y_{n}^j := x^j_{n} - {R_n} e_d \in B_{L_0}(0)$ (see \eqref{2}), we observe that 
$$
\limsup_{n \to \infty} \|f^j_{n}\|_{L^q(B_{L_0+2}(0))} \geq \limsup_{n \to \infty} \|f^j_{n}\|_{L^q(B_2(y^j_{n}))} = \limsup_{n \to \infty} \|e^j_{n}\|_{L^q(B_2(x^j_{n}))} \geq \delta_j/2 > 0,
$$
which yields $\|f^j\|_{L^q(B_{L_0+2}(0))} \geq \delta_j/2 > 0$. Thus we obtain $f^j \not\equiv 0$.

One finds from the assumption \eqref{orth-R} that
\begin{align*}
0 &= \int_{A_{R_n} \cap B_L({R_n} e_d)} e^j_{n} e^k_{n} u_n^{q-2} \, \d x + \int_{A_{R_n} \setminus B_L({R_n} e_d)} e^j_{n} e^k_{n} u_n^{q-2} \, \d x.
\end{align*}
Thanks to Corollary \ref{C:conv}, passing to the limits as $n\to \infty$ and then as $L \to \infty$, we observe from \eqref{barrier-v} and \eqref{c:e2} that
\begin{align*}
\int_{A_{R_n} \cap B_L({R_n} e_d)} e^j_{n} e^k_{n} u_n^{q-2} \, \d x 
&= \int_{(A_{R_n} - {R_n}e_d) \cap B_L(0)} f^j_{n} f^k_{n} v^{q-2}_{n} \, \d y\\
&\stackrel{n\to \infty}\to \int_{S \cap B_L(0)} f^j f^k v^{q-2}_\infty \, \d y
\stackrel{L \to \infty}\to \int_S f^j f^k v^{q-2}_\infty \, \d y,
\end{align*}
where $v_{n} := u_n(\,\cdot\,+{R_n} e_d)$, and moreover, since $|x - {R_n}e_d| \geq L$ for $x \in A_{R_n} \setminus B_L({R_n} e_d)$, we derive from \eqref{barrier} that
\begin{align*}
\left| \int_{A_{R_n} \setminus B_L({R_n}e_d)} e^j_{n} e^k_{n} u_n^{q-2} \, \d x \right|
&\leq \|u_n \|_{L^q(A_{R_n} \setminus B_L({R_n}e_d))}^{q-2} \\
&\leq C_1^{q-2} \left( \int_{A_{R_n} \setminus B_L({R_n}e_d)} \e^{-\lambda q|x - {R_n} e_d|} \,\d x \right)^{(q-2)/q}\\
&\leq C_1^{q-2} \left( \int_{\R^d \setminus B_L(0)} \e^{-\lambda q|x|} \, \d x\right)^{(q-2)/q} 
\stackrel{L\to\infty}\to 0.
\end{align*}
Thus it follows that
\begin{align*}
\int_S \nabla f^j \cdot \nabla f^k \, \d y = (q-1)\int_S f^j f^k v^{q-2}_\infty \, \d y = 0 \quad \mbox{ if } \ j \neq k
\end{align*}
for $j,k \in \{1,\ldots,d\}$. However, it is a contradiction to the assumption (S). Therefore the dimension of $\ker L_{u_R}$ is exactly $d-1$, and hence $\ker L_{u_R} = \mathrm{span}\{\Omega_1 u_R, \ldots, \Omega_{d-1} u_R\}$, for $R > 0$ large enough. This completes the proof.
\end{proof}

In view of Remark \ref{R:scaling}, Proposition \ref{P:nondeg} follows immediately from Lemma \ref{L:les2}.

\section{Morse--Bott condition and relative error convergence}\label{S:relerr}

\subsection{Some preliminaries} Let $\Omega$ be a bounded domain in $\R^d$ with boundary $\partial \Omega$. We define a functional $J : H^1_0(\Omega) \to \R$ by
\begin{equation*}
J(w) := \frac 12 \int_{\Omega} |\nabla w(x)|^2 \, \d x - \frac{\mathsf{c}_q}{q} \int_{\Omega} |w(x)|^q \, \d x \quad \mbox{ for } \ w \in H^1_0(\Omega),
\end{equation*}
where $2 < q < 2^*$ and $\mathsf{c}_q = (q-1)/(q-2)$ (see \eqref{cq}). We also define the Rayleigh quotient $\mathcal{Q} : H^1_0(\Omega) \setminus \{0\} \to [C_\Omega^{-1},\infty)$ associated with \eqref{Sobolev} by
$$
\mathcal{Q}(w) := \frac{\|\nabla w\|_{L^2(\Omega)}}{\|w\|_{L^q(\Omega)}} \quad \mbox{ for } \ w \in H^1_0(\Omega) \setminus \{0\},
$$
whose minimizers are exactly least-energy solutions to \eqref{elliptic} up to multiplying constants and conversely. We recall that $\mathcal Q(v(\cdot))$ is nonincreasing (see~\cite{BH80} and also~\cite{AK13,AS25}) and hence, if $v(s)$ converges to a limit $\phi \neq 0$ in $H^1_0(\Omega)$ as $s \to \infty$, we have $\mathcal{Q}(\phi) \leq \mathcal{Q}(v(s)) \leq \mathcal{Q}(v_0)$. 

Let $\phi$ be a positive solution to \eqref{elliptic}. We have $\|\phi\|_{H^1_0(\Omega)} = \mathsf{c}_q^{-1/(q-2)}\mathcal{Q}(\phi)^{q/(q-2)}$ immediately. Using Moser's iteration, one has
$$
\|\phi\|_{L^\infty(\Omega)} \leq B
$$
with some constant $B = B(d,\Omega,q,\|\phi\|_{H^1_0(\Omega)})$ depending only on $d,\Omega,q,\|\phi\|_{H^1_0(\Omega)}$. In particular, if $\phi$ is a least-energy solution to \eqref{elliptic}, it follows from the minimality that
$$
\|\phi\|^2_{H^1_0(\Omega)} = \mathsf{c}_q \|\phi\|^q_{L^q(\Omega)} = \left(\mathsf{c}_q C^q_{\Omega} \right)^{-2/(q-2)},
$$
and hence, $B$ depends only on $d,\Omega,q$.

\subsection{Morse--Bott condition} We shall prove Theorem \ref{T:rate} in a more general form (see Theorem \ref{T:rate-g} below). To this end, we introduce the following assumption on $J$ at a positive solution $\phi$ to \eqref{elliptic}:
\begin{enumerate}
 \item[(MB)] There exists a connected closed $C^1$ submanifold $\M_\phi$ of $H^1_0(\Omega)$ such that the following {\rm (i)--(iii)} are satisfied\/:
\begin{enumerate}
 \item[(i)] $\phi \in \M_\phi$,
 \item[(ii)] all elements of $\M_\phi$ are (positive) solutions to \eqref{elliptic},
 \item[(iii)] $\ker \Lphi = T_\phi \M_\phi$, where $\Lphi = J''(\phi)$ and $T_\phi \M_\phi$ stands for the tangent space to $\M_\phi$ at $\phi$.
\end{enumerate}
\end{enumerate}
Since $\dim \ker \Lphi < \infty$ due to Fredholm alternatives, we note that $\M_\phi$ is finite-dimensional.

\begin{remark}[Morse--Bott property and integrable kernel]\label{R:MB}
Assumption (MB) is also referred to as the \emph{Morse--Bott property} of $J$ at $\phi$ (see, e.g.,~\cite{Fe20,FeMa}). In~\cite{AdSi88,AlAl}, the linearized operator $\Lphi = J''(\phi)$ at $\phi$ is said to have an \emph{integrable kernel} under a similar condition, which is (in general) slightly weaker than the Morse--Bott property but still equivalent to that under some additional condition (see~\cite{Fe20}). In these literatures, the \L ojasiewicz--Simon gradient inequality is often derived with the best exponent $\theta = 1/2$ under the Morse--Bott property; however, the constant of the gradient inequality is never specified. In order to derive the optimal convergence rate as in Theorem \ref{T:rate}, it is also necessary to determine the optimal constant of the gradient inequality, which will be achieved in Section \ref{S:QGIa}. This point constitutes one of the main novelties of the present paper.
\end{remark}

As for the thin annular domain case $\Omega = A_R$ as in Proposition \ref{P:nondeg}, we can observe that

\begin{proposition}\label{P:int-ker}
Assume that {\rm (A)} is satisfied and let $R \geq R_0$. Let $\Omega = A_R$ and let $\phi > 0$ be a least-energy solution of \eqref{elliptic}. Then $J$ satisfies {\rm (MB)} at $\phi$.
\end{proposition}

\begin{proof}
Let $\phi > 0$ be a nontrivial solution to \eqref{elliptic}. We may assume that $\phi$ is foliated Schwarz symmetric with respect to the $x_d$-axis (see~\cite{Pac02,DaPa}). We define
\begin{equation}\label{Mphi}
\M_\phi := \left\{ \phi(\Ro^{-1}\,\cdot\,) \colon \ \Ro \in \mathrm{SO}(d) \right\}.
\end{equation}
Then $\M_\phi$ is a connected closed smooth submanifold of $H^1_0(A_R)$, and by rotational invariance, all elements of $\M_\phi$ are positive least-energy solutions to \eqref{elliptic}. Moreover, we observe that $\mathrm{span}\,\{\Omega_j \phi\}_{j=1,\ldots,d-1} \subset T_\phi \M_\phi \subset \ker \Lphi$. Thus it follows immediately from (A) that $T_\phi \M_\phi \subset \ker \Lphi = \mathrm{span}\, \{\Omega_j \phi\}_{j=1,\ldots,d-1} \subset T_\phi \M_\phi$. This completes the proof.
\end{proof} 

In what follows, assume that $\Omega$ is a bounded $C^{1,1}$ domain in $\R^d$. Let $\phi$ be a positive solution to \eqref{elliptic} such that
\begin{equation}\label{H-int}
J \mbox{ satisfies {\rm (MB)} at } \phi.
\end{equation}
Then we have
\begin{equation}\label{nu}
\underbrace{\nu_1 = - \mathsf{c}_q(q-2) \leq \cdots \leq \nu_{\ell-1}}_{\text{negative eigenvalues}} < \underbrace{0 = \nu_{\ell} = \cdots = \nu_{k-1}}_{\text{zero eigenvalues}} < \underbrace{\nu_k \leq \nu_{k+1} \leq \cdots}_{\text{positive eigenvalues}}
\end{equation}
for some $1 < \ell < k < \infty$. We remark that $\nu_1$ is simple, that is, $\nu_1 < \nu_2$. 
\begin{lemma}[Morse--Bott manifold]\label{L:MB}
Let $\M_\phi$ be a connected closed $C^1$ submanifold of $H^1_0(\Omega)$ satisfying {\rm (i)--(iii)} of Condition {\rm (MB)}. Then the following conditions hold true\/{\rm :}
\begin{enumerate}
 \item[(i)] $J$ is constant on $\M_\phi$.
 \item[(ii)] For every $w \in H^1_0(\Omega)$, there exists an element $m_w \in \M_\phi$ such that
\begin{equation}\label{min}
\|w - m_w\|_{H^1_0(\Omega)} = \inf_{m \in \M_\phi} \|w - m\|_{H^1_0(\Omega)}.
\end{equation}
We denote by $\mathbb{P}_{\M_\phi}(w)$ the set of such $m_w \in \M_\phi$ satisfying \eqref{min}. 
 \item[(iii)] For every $w \in {H^1_0(\Omega)}$, it holds that
\begin{equation}\label{o1}
\left( w - m_w, z \right)_{H^1_0(\Omega)} = 0 \quad \mbox{ for } \ m_w \in \mathbb{P}_{\M_\phi}(w), \ z \in T_{m_w} \M_\phi.
\end{equation}
\end{enumerate}
\end{lemma}

\begin{proof}
We first prove (i). 
For any $C^1$ curves $\gamma : [0,1] \to \M_\phi$ satisfying $\gamma(0) = \phi$ (see (i) of (MB)), it follows from (ii) of (MB) that
$$
\dfrac{\d}{\d t} J(\gamma(t)) = \left\langle J'(\gamma(t)), \frac{\d \gamma}{\d t}(t) \right\rangle_{H^1_0(\Omega)} = 0 \quad \mbox{ for all } \ t \in (0,1),
$$
which along with the connectedness of $\M_\phi$ yields (i). We next prove (ii). Let $w \in H^1_0(\Omega)$ and let $(m_n)$ be a sequence in $\M_\phi$ such that $\|w - m_n\|_{H^1_0(\Omega)} \to \inf_{m \in \M_\phi} \|w - m\|_{H^1_0(\Omega)}$. We then observe that $J'(m_n) = 0$ and $J(m_n)$ is constant for $n \in \N$. Since $J$ satisfies the Palais--Smale condition, we obtain, up to a subsequence, $m_n \to m_w$ in $H^1_0(\Omega)$ for some $m_w \in H^1_0(\Omega)$. Since $\M_\phi$ is closed in $H^1_0(\Omega)$, we conclude that $m_w \in \M_\phi$ and $\|w - m_w\|_{H^1_0(\Omega)} = \inf_{m \in \M_\phi} \|w - m\|_{H^1_0(\Omega)}$. Finally, (iii) follows immediately from the minimality in (ii) and the $C^1$ structure of $\M_\phi$. This completes the proof.
\end{proof}


\subsection{Relative error estimate} We next derive a relative error estimate, which was originally established for the limit $\phi$ of each solution $v(s)$ to \eqref{v-pde}--\eqref{v-ic} in~\cite[Theorem 4.1]{BF21} (cf.~\cite[Lemma 6.1]{A21}). Here we are concerned with the relative error between $v(s)$ and $\phi_s \in \mathbb{P}_{\M_\phi}(v(s))$ instead of $\phi$.

\begin{lemma}[Relative error estimates (cf.~\cite{BF21,A21})]\label{L:re-exp}
Let $\Omega$ be a bounded $C^{1,1}$ domain in $\R^d$ with boundary $\partial \Omega$. Let $\phi = \phi(x)$ be a positive solution to \eqref{elliptic} such that \eqref{H-int} is satisfied. Let $v = v(x,s)$ be a nonnegative energy solution to \eqref{v-pde}--\eqref{v-ic} and let $(\phi_s)$ be a family of positive solutions to \eqref{elliptic} such that $\phi_s \in \mathbb{P}_{\M_\phi}(v(s))$ for $s \geq 0$ {\rm (}see {\rm (ii)} of Lemma {\rm \ref{L:MB}}{\rm )}. Then there exists a constant $C_0$ depending only on $d$, $\Omega$, $q$, $\mathcal{Q}(v_0)$, $\|\phi\|_{H^1_0(\Omega)}$ such that
\begin{equation}\label{re-bdd}
 \left\| \frac{v(s)}{\phi_s} - 1\right\|_{L^\infty({\Omega})} \leq C_0 \quad \mbox{ for } \ s \geq s_*
\end{equation}
with $s_* := \log 2 + 1 > 0$. Furthermore, if $v(s) \to \phi$ strongly in $H^1_0(\Omega)$ as $s \to \infty$, then 
\begin{equation*}
 \left\| \frac{v(s)}{\phi_s} - 1\right\|_{L^\infty({\Omega})} \to 0 \quad \mbox{ as } \ s \to \infty.
\end{equation*}
In addition, if there exist constants $c_1, \mu > 0$ such that $\|v(s)-\phi\|_{H^1_0({\Omega})} \leq c_1 \e^{-\mu s}$ for any $s \geq 0$, then there exists a constant $C_1$ depending only on $d$, $\Omega$, $q$, $\mathcal{Q}(v_0)$, $c_1$, $\mu$ such that
\begin{equation}\label{re-exp}
\left\| \frac{v(s)}{\phi_s} - 1\right\|_{L^\infty({\Omega})} 
\leq C_1 \e^{-b s} \quad \mbox{ for } \ s \geq s_*
\end{equation}
with $b := \mu/(2d) > 0$.
\end{lemma}

\begin{proof}
Throughout this proof, we denote by $C$ general constants which depend only on $d$, $\Omega$, $q$, $\mathcal{Q}(v_0)$, $\|\phi\|_{H^1_0(\Omega)}$ and may vary from line to line. Thanks to Theorem 4.1 of~\cite{BF21}, we deduce that
\begin{align}
\left\|\frac{v(s)}{\phi}-1\right\|_{L^\infty({\Omega})}
&\leq C \frac{\e^{L (s-s_1)}}{s-s_1} (1 + s - s_1) \sup_{\sigma \geq s_1} \left( \int_{\Omega} |v^{q-1}(\sigma)-\phi^{q-1}| \, \d x \right)^{1/d}\nonumber\\
 &\quad + C(s-s_1) \e^{L(s-s_1)},\label{quan-re}
\end{align}
where $L := 2/(q-2) > 0$, for any $s > s_1 \geq \log 2$. Just for the convenience of the reader, we also give a proof of \eqref{quan-re} in Appendix \S \ref{A:S:smoothing}. Moreover, we observe that
\begin{align}
\int_{\Omega} |v^{q-1}(\sigma)-\phi^{q-1}| \, \d x
&\leq C \|v(\sigma)-\phi\|_{H^1_0({\Omega})}
\quad \mbox{ for } \ \sigma \geq 0.\label{quan-re3}
\end{align}
Indeed, according to~\cite[Proposition 2]{AK13}, we have 
\begin{equation}\label{bound}
\sup_{\sigma \geq 0} \|v(\sigma)\|_{L^q({\Omega})} \leq \mathsf{c}_q^{-1/(q-2)} \mathcal{Q}(v_0)^{2/(q-2)}, \quad \sup_{\sigma \geq 0} \|v(\sigma)\|_{H^1_0({\Omega})} \leq \mathsf{c}_q^{-1/(q-2)} \mathcal{Q}(v_0)^{q/(q-2)}.
\end{equation}

First, substituting $s_1 = s-1$ and $\phi = \phi_s$ to \eqref{quan-re} along with \eqref{quan-re3}, we observe that
\begin{align*}
\left\|\frac{v(s)}{\phi_s}-1\right\|_{L^\infty({\Omega})}
 &\leq 2 C \e^{L} \sup_{\sigma \geq s-1} \|v(\sigma)-\phi_s\|_{H^1_0({\Omega})}^{1/d} + C\e^{L}
\end{align*}
for $s \geq s_* := \log 2 + 1$. Thus \eqref{re-bdd} follows from \eqref{bound} and the fact that $\|\phi_s\|_{H^1_0(\Omega)} = \|\phi\|_{H^1_0(\Omega)}$ (see (i) of Lemma \ref{L:MB}) for $s \geq 0$. 

Next, we let $\delta : [0,\infty) \to (0,\infty)$ be a function satisfying $\lim_{s \to \infty} \delta(s) = 0$. For any $s \geq s_*$ large enough, setting $s_1 = s - \delta(s)$, we have $0 < \delta(s) \leq 1$ and $s_1 \geq \log 2$. 
Moreover, putting $\phi = \phi_s$ in~\eqref{quan-re} and \eqref{quan-re3}, we see that
\begin{align*}
\left\|\frac{v(s)}{\phi_s}-1\right\|_{L^\infty({\Omega})}
\leq C \e^L (\delta(s)^{-1} + 1) \sup_{\sigma \geq s-1} \|v(\sigma)-\phi_s\|_{H^1_0({\Omega})}^{1/d} + C \delta(s) \e^L,
\end{align*} 
provided that $s - \delta(s) \geq \log 2$. We note that
\begin{align}
\|v(\sigma)-\phi_s\|_{H^1_0({\Omega})}
&\leq \|v(\sigma)-\phi\|_{H^1_0({\Omega})} + \|\phi-v(s)\|_{H^1_0({\Omega})} + \|v(s)-\phi_s\|_{H^1_0({\Omega})} \nonumber\\
&\leq \|v(\sigma)-\phi\|_{H^1_0({\Omega})} + 2\|\phi-v(s)\|_{H^1_0({\Omega})} \label{phisphi}
\end{align}
for $s,\sigma \geq 0$. Therefore using \eqref{phisphi}, we obtain
\begin{align*}
\left\|\frac{v(s)}{\phi_s}-1\right\|_{L^\infty({\Omega})}
&\leq C \e^L (\delta(s)^{-1} + 1) \left[ \sup_{\sigma \geq s-1} \|v(\sigma)-\phi\|_{H^1_0({\Omega})} + 2\|\phi-v(s)\|_{H^1_0({\Omega})} \right]^{1/d}\\
&\quad + C \delta(s) \e^L,
\end{align*}
provided that $s - \delta(s) \geq \log 2$. If $v(s) \to \phi$ strongly in $H^1_0(\Omega)$ as $s \to \infty$, putting
$$
\delta(s) := \left[ \sup_{\sigma \geq s-1} \|v(\sigma)-\phi\|_{H^1_0({\Omega})} + 2\|\phi-v(s)\|_{H^1_0({\Omega})} \right]^{1/(2d)} \to 0 \quad \mbox{ as } \ s \to \infty,
$$
we see that
$$
\left\|\frac{v(s)}{\phi_s}-1\right\|_{L^\infty({\Omega})} \to 0 \quad \mbox{ as } \ s \to \infty.
$$
Here we also note that $\|\phi_s\|_{H^1_0(A_R)} = \|\phi\|_{H^1_0(A_R)}$ can be estimated in terms of $\mathcal{Q}(v_0)$, since this property is inherited from $v(s)$ as its limit as $s \to \infty$.

In addition, if $\|v(s)-\phi\|_{H^1_0({\Omega})} \leq c_1 \e^{-\mu s}$ for $s \geq 0$, putting $\delta(s) = \e^{-as}$, where $a > 0$ will be determined later, we can deduce that
\begin{align*}
\left\|\frac{v(s)}{\phi_s}-1\right\|_{L^\infty({\Omega})}
\leq C \e^L (\e^{as} + 1) \left( c_1 \e^{-\mu (s-1)} + 2 c_1 \e^{-\mu s} \right)^{1/d} + C \e^{-a s} \e^L
\end{align*}
for $s \geq s_* = \log 2 +1$ (hence, $s - \e^{-as} \geq \log 2$). Hence choosing $a = \mu/(2d)$, we obtain \eqref{re-exp} with $b = \mu/(2d) > 0$. This completes the proof.
\end{proof}

\section{Quantitative gradient inequality under the Morse--Bott condition}\label{S:QGIa}

Let $\phi > 0$ be a positive solution to \eqref{elliptic} and let $v$ be a nonnegative energy solution to \eqref{v-pde}--\eqref{v-ic} such that $v(s) \to \phi$ strongly in $H^1_0(\Omega)$ as $s \to \infty$. Under Assumption \eqref{H-int}, we can verify that

\begin{lemma}[Morse--Bott]\label{L:MB2}
Assume that \eqref{H-int} is satisfied. Then it holds that
\begin{equation}\label{H-int+}
\mbox{there exists } r_0 > 0 \mbox{ such that } \ker \mathcal{L}_m = T_m \M_\phi \mbox{ for all } m \in \M_\phi \cap B_{H^1_0(\Omega)}(\phi;r_0),
\end{equation}
where $\M_\phi$ is a connected closed $C^1$ submanifold of $H^1_0(\Omega)$ satisfying {\rm (i)--(iii)} of Condition {\rm (MB)} and $\mathcal{L}_m : H^1_0(\Omega) \to H^{-1}(\Omega)$ is the linearized operator defined by $\mathcal L_m w = - \Delta w - \mathsf{c}_q (q-1) m^{q-2} w$ for $w \in H^1_0(\Omega)$.  
\end{lemma}

The above fact is known as a conclusion of the Morse--Bott lemma (see~\cite[\S 1.1]{Fe20}), but for the convenience of the reader, we give a proof of \eqref{H-int+} in Appendix \S \ref{A:MB}. It also implies that
$$
\dim \ker \mathcal{L}_m = \dim \ker \Lphi \quad \mbox{ for all } \ m \in \M_\phi \cap B_{H^1_0(\Omega)}(\phi;r_0).
$$
Then for each $s \geq 0$ (large enough) so that $2 \|v(s) - \phi\|_{H^1_0(\Omega)} < r_0$, we let $\phi_s \in \mathbb{P}_{\M_\phi}(v(s))$, that is, 
\begin{equation}\label{phis}
\|v(s) - \phi_s\|_{H^1_0(\Omega)} = \inf_{m \in \M_\phi} \|v(s) - m\|_{H^1_0(\Omega)}
\ \mbox{ and } \ v(s)-\phi_s \in (\ker \mathcal{L}_{\phi_s})^\perp,
\end{equation}
the latter of which follows from \eqref{o1} along with \eqref{H-int+} (indeed, we note that $\|\phi_s-\phi\|_{H^1_0(\Omega)} \leq \|\phi_s-v(s)\|_{H^1_0(\Omega)} + \|v(s)-\phi\|_{H^1_0(\Omega)} \leq 2 \|v(s)-\phi\|_{H^1_0(\Omega)} < r_0$) and will play a crucial role in this section. Concerning the thin annular domain case $\Omega = A_R$ for $R \geq R_0$, Condition \eqref{H-int+} with $r_0 = \infty$ follows directly from Proposition \ref{P:int-ker} by rotation. 

Next, for each $s > 0$, let us consider the weighted eigenvalue problem,
\begin{equation}\label{evs}
- \Delta e^s_j = \mu_j^{s} \phi_s^{q-2} e^s_j \ \mbox{ in } {\Omega}, \quad e^s_j = 0 \ \mbox{ on } \partial \Omega,
\end{equation}
where $0 < \mu_1^s = \mathsf{c}_q < \mu_2^{s} \leq \cdots$ stand for the eigenvalues and $e^s_j \in H^1_0(\Omega)$ denotes an eigenfunction (orthonormalized in ${H^1_0(\Omega)}$) corresponding to $\mu_j^{s}$ for each $j \in \N$. We similarly denote by $\{(\mu_j,e_j)\}_{j\in \N}$ the eigenpairs of
\begin{equation}\label{ev-}
- \Delta e_j = \mu_j \phi^{q-2} e_j \ \mbox{ in } \Omega, \quad e_j = 0 \ \mbox{ on } \partial \Omega.
\end{equation}
Then both $(e^s_j)$ and $(e_j)$ are complete orthonormal systems of ${H^1_0(\Omega)}$, and hence, $(- \Delta e^s_j)$ and $(- \Delta e_j)$ are complete orthonormal systems of the dual space $H^{-1}(\Omega)$ of $H^1_0(\Omega)$. 

\begin{remark}[Annular domain case]
In particular, when $\Omega = A_R$, we note that  $\mu^s_j = \mu_j$, where $\mu_j$ is the $j$-th eigenvalue of \eqref{ep}, for all $j \in \N$ by rotation. Hence we may not need Lemmas \ref{L:ev-var} and \ref{L:ev-pert} below for proving Theorem \ref{T:rate}; however, they will be used to prove more general Theorem \ref{T:rate-g} below, which yields Theorem \ref{T:rate} as a corollary. On the other hand, $e^s_j$ may differ from the eigenfunction $e_j$ of \eqref{ep} even for the case $\Omega = A_R$.  
\end{remark}

Then it is well known that

\begin{lemma}[Variational representation of eigenvalues]\label{L:ev-var}
For each $s \geq 0$ and $j \geq 1$, the eigenvalue $\mu^s_j > 0$ of \eqref{evs} can be characterized as the following max-min value\/{\rm :}
\begin{align*}
\frac 1 {\mu^s_j} &= \sup_{\substack{Y \subset H^1_0(\Omega)\\[0mm]\dim Y = j}} \inf_{\substack{w \in Y\\[0mm]\|w\|_{H^1_0(\Omega)}=1}} \int_\Omega \phi_s^{q-2} w^2 \, \d x.
\end{align*}
Here $Y$ denotes a subspace of $H^1_0(\Omega)$. Similarly, each eigenvalue $\mu_j > 0$ of \eqref{ev-} can be written as
\begin{equation*}
\frac 1 {\mu_j} = \sup_{\substack{Y \subset H^1_0(\Omega)\\[0mm]\dim Y = j}} \inf_{\substack{w \in Y\\[0mm]\|w\|_{H^1_0(\Omega)}=1}} \int_\Omega \phi^{q-2} w^2 \, \d x.
\end{equation*}
\end{lemma}

Therefore as in~\cite[Lemma 3.4]{AM23+}, recalling that $\|\phi_s\|_{H^1_0(\Omega)} = \|\phi\|_{H^1_0(\Omega)}$ (see (i) of Lemma \ref{L:MB}) and $\|\phi_s - \phi\|_{H^1_0(\Omega)} \leq 2 \|v(s)-\phi\|_{H^1_0(\Omega)}$ for $s \geq 0$, we can deduce that 

\begin{lemma}[Convergence of eigenvalues]\label{L:ev-pert}
There exists a positive constant $C$ which depends only on $q$, $C_\Omega$ and $\|\phi\|_{H^1_0(\Omega)}$ such that
\begin{equation*}
\left|\frac 1 {\mu^s_j} - \frac 1 {\mu_j}\right| \leq C \| \phi_s - \phi \|_{L^q(\Omega)}^{\rho} \quad \mbox{ for } \ s \geq 0 \ \mbox{ and } \ j \in \N,
\end{equation*}
where $\rho := \min\{q-2,1\} \in (0,1]$. Moreover, for each $j \in \N$, it holds that
\begin{equation*}
 \mu^s_j \to \mu_j \quad \mbox{ as } \ s \to \infty,
\end{equation*}
provided that $v(s) \to \phi$ in $H^1_0(\Omega)$ as $s \to \infty$.
\end{lemma}
Set
$$
\H_s := L^2({\Omega};\phi_s^{q-2}\d x), \quad \H_s' := L^2({\Omega};\phi_s^{2-q}\d x)
$$
equipped with norms defined by
\begin{align*}
\|w\|_{\H_s}^2 &:= \int_{\Omega} |w(x)|^2 \phi_s^{q-2}(x) \, \d x \quad \mbox{ for } \ w \in \H_s,\\
\|w\|_{\H_s'}^2 &:= \int_{\Omega} |w(x)|^2 \phi_s^{2-q}(x) \, \d x \quad \mbox{ for } \ w \in \H_s',
\end{align*}
respectively. Then $H^1_0(\Omega)$ is continuously embedded in $\H_s$, and therefore, $\H_s'$ is continuously embedded in $H^{-1}(\Omega)$. For later use, we also remark that $( (\mu^s_j)^{-1/2}(-\Delta e^s_j))$ forms a (complete) orthonormal system in $\H_s'$. Furthermore, as we shall see in the proof of Theorem \ref{T:rate-g} below (see, in particular, \eqref{U2}), we have
$$
J'(v(s)) \stackrel{\eqref{ggf}}= - \partial_s \left(v^{q-1}\right)(s) \in \H_s' \quad \mbox{ for a.e. } \ s > 0.
$$

Now, we are ready to state a key lemma, which provides a quantitative gradient inequality, to prove Theorem \ref{T:rate} in a more general form (see Theorem \ref{T:rate-g} below).

\begin{lemma}[Quantitative gradient inequality under the Morse--Bott condition]\label{L:QGIa}
Let $v = v(x,s)$ be a nonnegative energy solution to \eqref{v-pde}--\eqref{v-ic} and let $\phi = \phi(x)$ be a positive solution to \eqref{elliptic} such that $v(s) \to \phi$ strongly in $H^1_0(\Omega)$ as $s \to \infty$. Assume that \eqref{H-int} is satisfied, {\rm (}and then \eqref{H-int+} also holds{\rm )}. In addition, let $\phi_s \in \mathbb{P}_{\M_\phi}(v(s))$ for each $s \geq 0$. Then there exist positive constants $c_0 \geq 0$ and $s_0 \geq 0$ such that
\begin{equation}\label{QGI}
 J(v(s)) - J(\phi_s)
\leq \frac1{2 \nu_{{k}}} \left( 1 + c_0 \|v(s)-\phi\|_{H^1_0({\Omega})}^{\rho} \right) \int_{{\Omega}} |J'(v(s))|^2 \phi_s^{2-q} \, \d x
\end{equation}
for a.e.~$s \geq s_0$. Here $\nu_{{k}}$ is the least positive eigenvalue for the linearized operator $\Lphi = - \Delta - \mathsf{c}_q(q-1)\phi^{q-2}$ at $\phi$ {\rm (}see \eqref{nu}{\rm )} and $\rho = \min\{q-2,1\} \in (0,1]$. Furthermore, there exists $\delta_0 > 0$ such that the above assertion holds with $s_0 = 0$, provided that $\sup_{s \geq 0} \|v(s)-\phi\|_{H^1_0({\Omega})} < \delta_0$. The constants $c_0$ and $\delta_0$ depend only on $d,\Omega,q,\mathcal{Q}(v_0)$, $r_0$, $\nu_{k}$ and the largest negative eigenvalue $\nu_{\ell - 1}$ of $\Lphi$ {\rm (}see also Remark {\rm \ref{R:const}} below{\rm )}.
\end{lemma}

The following proof is closer to the original proof of a quantitative gradient inequality developed in~\cite{A21}, which concerns nondegenerate profiles only (see also~\cite{AM23+}), than that of~\cite{KoYu24}, which modified the proof of~\cite{A21} in order to treat degenerate profiles (in a different setting from the present paper). Moreover, our proof may be simpler than that of~\cite{KoYu24}.

\begin{proof} Throughout this proof, we denote by $C$ general constants which may depend at most on $d,\Omega,q,\mathcal{Q}(v_0), \nu_k = \mu_k - \mathsf{c}_q(q-1) > 0$ and $\nu_{\ell-1} = \mu_{\ell-1} - \mathsf{c}_q(q-1) < 0$ and may vary from line to line. We expand $J$ and $J'$ in a Taylor series around $\phi_s$. Since $\phi_s$ is also a critical point of $J$, we see that
\begin{align}\label{Jvs}
 J(v(s)) = J(\phi_s) + \frac12 \langle J''(\phi_s)(v(s)-\phi_s), v(s)-\phi_s \rangle_{H^1_0(\Omega)} + \mathscr{R}_s \quad \mbox{ for } \ s \geq 0,
\end{align}
where $\mathscr{R}_s$ denotes the remainder term given by
$$
\mathscr{R}_s := - \mathsf{c}_q \int_\Omega \left( \frac1q |v(s)|^q - \frac1q |\phi_s|^q - \phi_s^{q-1}(v(s)-\phi_s) - \frac{q-1}2 \phi_s^{q-2}(v(s)-\phi_s)^2 \right) \, \d x.
$$
Here, by Sobolev's inequality, we note that
\begin{equation}\label{Rs-est}
|\mathscr{R}_s| 
\leq C \|v(s)-\phi_s\|_{H^1_0(\Omega)}^{2+\rho}  \quad \mbox{ for } \ s \geq 0
\end{equation}
with $\rho := \min \{q-2,1\} \in (0,1]$ (see~\cite[Proof of Lemma 3.10]{AM23+}). Furthermore, it follows that
\begin{equation}\label{J'vs}
J'(v(s)) = J''(\phi_s)(v(s)-\phi_s) + r_s  \quad \mbox{ for } \ s \geq 0,
\end{equation}
where $r_s \in L^{q'}(\Omega)$ is given by
$$
r_s := - \mathsf{c}_q \left( v^{q-1}(s) - \phi_s^{q-1} - (q-1) \phi_s^{q-2} (v(s) - \phi_s) \right).
$$
Then we also note that
\begin{equation}\label{rs-est}
\|r_s\|_{H^{-1}(\Omega)} \lesssim \|r_s\|_{L^{q'}(\Omega)} \leq C \|v(s)-\phi_s\|_{H^1_0(\Omega)}^{1+\rho}  \quad \mbox{ for } \ s \geq 0.
\end{equation}

Setting
$$
\alpha_j^s := (v(s)-\phi_s, e^s_j)_{H^1_0(\Omega)} \quad \mbox{ for } \ j \in \N,
$$
we see that
$$
v(s)-\phi_s = \sum_{j=1}^\infty \alpha^s_j e^s_j \ \mbox{ in } {H^1_0(\Omega)}.
$$
Moreover, let $\Ls := J''(\phi_s) = -\Delta - \mathsf{c}_q(q-1) \phi_s^{q-2}$ and set $\nu_j^s := \mu_j^s - \mathsf{c}_q(q-1)$ for $j \in \N$. Then
$$
\Ls e^s_j = \nu_j^s \phi_s^{q-2} e^s_j \ \mbox{ in } \Omega, \quad e^s_j = 0 \ \mbox{ on } \partial \Omega
$$
for $j \in \N$. In view of \eqref{nu} and \eqref{H-int+}, we see that
$$
\underbrace{\nu_1^s = \nu_1 = - \mathsf{c}_q(q-2) \leq \cdots \leq \nu_{\ell-1}^s}_{\text{negative eigenvalues}} < \underbrace{0 = \nu_{\ell}^s = \cdots = \nu_{k-1}^s}_{\text{zero eigenvalues}} < \underbrace{\nu_k^s \leq \nu_{k+1}^s \leq \cdots}_{\text{positive eigenvalues}}
$$
when $\|v(s)-\phi\|_{H^1_0(\Omega)}$ is small enough. Indeed, by virtue of Lemma \ref{L:ev-pert}, we have $\nu^s_j \to \nu_j$ as $\|v(s)-\phi\|_{H^1_0(\Omega)} \to 0$. Moreover, it follows from \eqref{H-int+} that the number of zero eigenvalues is same as \eqref{nu} if $\|\phi_s - \phi\|_{H^1_0(\Omega)} \leq 2 \|v(s)-\phi\|_{H^1_0(\Omega)} < r_0$. Therefore $\ell$ and $k$ are independent of $s$ for which $\|v(s)-\phi\|_{H^1_0(\Omega)}$ is small enough, and moreover, $\nu^s_{\ell-1}$ and $\nu^s_k$ are uniformly away from zero for such $s$.

Since $v(s) - \phi_s \in (\ker \Ls)^\perp$ (see \eqref{phis}) and $e^s_j \in \ker \Ls$ for $j \in \{\ell,\ldots,k-1\}$, we obtain $\alpha^s_j = 0$ for $j \in \{\ell,\ldots, k-1 \}$. Setting 
\begin{align*}
\beta^s_j &:= (J'(v(s)), - \Delta e^s_j)_{H^{-1}(\Omega)} = \langle J'(v(s)),e^s_j \rangle_{H^1_0(\Omega)},\\
\gamma^s_j &:= (r_s, - \Delta e^s_j)_{H^{-1}(\Omega)} = \langle r_s,e^s_j \rangle_{H^1_0(\Omega)}
\end{align*}
for $j \in \N$, we deduce from \eqref{J'vs} that
$$
\sum_{j=1}^\infty \frac{\nu^s_j}{\mu^s_j} \alpha^s_j (-\Delta e^s_j) = \Ls \left( v(s)-\phi_s \right) = J'(v(s)) - r_s = \sum_{j=1}^\infty \left(\beta^s_j - \gamma^s_j\right)(-\Delta e^s_j)
$$
for a.e.~$s > 0$. Therefore it follows that
$$
\beta^s_j - \gamma^s_j = 0 \ \mbox{ for } \ j \in \{\ell, \ldots,k-1\},
\quad \alpha^s_j = \frac{\mu^s_j}{\nu^s_j} \left(\beta^s_j - \gamma^s_j\right) \ \mbox{ for } \ j \not\in \{\ell,\ldots, k-1\}.
$$
Furthermore, we also infer that
\begin{align*}
\|v(s)-\phi_s\|_{H^1_0(\Omega)}^2 
&= \sum_{j \not\in \{\ell,\ldots,k-1\}} \left( \frac{\mu^s_j}{\nu^s_j} \right)^2 \left(\beta^s_j - \gamma^s_j \right)^2 \\
&\leq C \left( \|J'(v(s))\|_{H^{-1}(\Omega)}^2 + \|r^s\|_{H^{-1}(\Omega)}^2 \right)\\
&\stackrel{\eqref{rs-est}}\leq C \left( \|J'(v(s))\|_{H^{-1}(\Omega)}^2 + \|v(s)-\phi_s\|_{H^1_0(\Omega)}^{2(1+\rho)} \right),
\end{align*}
which implies that
\begin{equation}\label{absorb}
 \|v(s)-\phi_s\|_{H^1_0(\Omega)}^2 \leq C \|J'(v(s))\|_{H^{-1}(\Omega)}^2,
\end{equation}
provided that $\|v(s)-\phi_s\|_{H^1_0(\Omega)}$ is small enough. Here we used the fact that
$$
0 > \frac{\mu^s_j}{\nu^s_j} \geq \frac{\mu^s_{\ell-1}}{\nu^s_{\ell-1}} \ \mbox{ for } \ j \leq \ell-1 \quad \mbox{ and } \quad 1 < \frac{\mu^s_j}{\nu^s_j} \leq \frac{\mu^s_k}{\nu^s_k} \ \mbox{ for } j \geq k,
$$
and moreover, by Lemma \ref{L:ev-pert}, we see that $|\mu^s_{\ell-1}/{\nu^s_{\ell-1}}| + \mu^s_k/\nu^s_k \leq C$ for $s \geq 0$ at which $\|v(s)-\phi\|_{H^1_0(\Omega)}$ is small enough. Combining these facts, since $\nu^s_j < 0$ for $j \leq \ell-1$ and $\nu^s_j \geq \nu^s_k > 0$ for $j \geq k$, we derive that
\begin{align}
\MoveEqLeft{
 J(v(s)) - J(\phi_s)
}\nonumber\\
&= \frac12 \langle \Ls(v(s)-\phi_s), v(s)-\phi_s \rangle_{H^1_0(\Omega)} + \mathscr{R}_s \nonumber\\
&= \frac12 \langle J'(v(s))-r_s, v(s)-\phi_s \rangle_{H^1_0(\Omega)} + \mathscr{R}_s \nonumber\\
&= \frac12 \sum_{j=1}^\infty \left(\beta^s_j-\gamma^s_j\right) \alpha^s_j + \mathscr{R}_s \nonumber\\
&= \frac12 \sum_{j=1}^{\ell-1} \frac{\mu_j^s}{\nu_j^s} \left( \beta^s_j - \gamma^s_j \right)^2 + \frac 12 \sum_{j = k}^\infty \frac{\mu_j^s}{\nu_j^s} \left( \beta^s_j - \gamma^s_j \right)^2 + \mathscr{R}_s\nonumber\\
&\stackrel{\eqref{Rs-est}}\leq \frac 12 \sum_{j=k}^\infty \frac{\mu_j^s}{\nu_j^s} \left(\beta^s_j-\gamma^s_j\right)^2 + C \|v(s)-\phi_s\|_{H^1_0(\Omega)}^{2+\rho} \nonumber\\
&= \frac 12 \sum_{j=k}^\infty \frac{\mu_j^s}{\nu_j^s} \left(\beta^s_j\right)^2 - \sum_{j=k}^\infty \frac{\mu_j^s}{\nu_j^s} \beta^s_j \gamma^s_j + \frac 12 \sum_{j=k}^\infty \frac{\mu_j^s}{\nu_j^s} \left(\gamma^s_j\right)^2 + C \|v(s)-\phi_s\|_{H^1_0(\Omega)}^{2+\rho} \nonumber\\
&\leq \frac 1{2\nu^s_k} \|J'(v(s))\|_{\mathcal{H}_s'}^2 + C \left( \|J'(v(s))\|_{H^{-1}(\Omega)} \|r_s\|_{H^{-1}(\Omega)} + \|r_s\|_{H^{-1}(\Omega)}^2 \right) \nonumber\\ 
&\quad + C \|v(s)-\phi_s\|_{H^1_0(\Omega)}^{2+\rho} \nonumber\\
&\leq \frac 1{2\nu^s_k} \|J'(v(s))\|_{\mathcal{H}_s'}^2 + C \left( \|J'(v(s))\|_{H^{-1}(\Omega)} \|v(s)-\phi_s\|_{H^1_0(\Omega)}^{1+\rho} + \|v(s)-\phi_s\|_{H^1_0(\Omega)}^{2(1+\rho)} \right) \nonumber\\ 
&\quad + C \|v(s)-\phi_s\|_{H^1_0(\Omega)}^{2+\rho} \nonumber\\
&\stackrel{\eqref{absorb}}\leq \frac 1{2\nu^s_k} \left( 1 + C\|v(s)-\phi_s\|_{H^1_0(\Omega)}^{\rho} \right) \|J'(v(s))\|_{\mathcal{H}_s'}^2.
\nonumber
\end{align}
Here we used the facts that 
$$
\sum_{j=k}^\infty \mu^s_j (\beta^s_j)^2 \leq \sum_{j=1}^\infty \mu^s_j (\beta^s_j)^2 = \|J'(v(s))\|_{\H_s'}^2
$$
and that $1 < \mu_j^s/\nu_j^s \leq \mu_k^s/\nu_k^s \leq C$ for $j \geq k$ and $s \geq 0$ at which $\|v(s)-\phi\|_{H^1_0(\Omega)}$ is small enough (see Lemma \ref{L:ev-pert}). We recall that $\|v(s)-\phi_s\|_{H^1_0(\Omega)} \leq \|v(s)-\phi\|_{H^1_0(\Omega)}$. Exploiting Lemma \ref{L:ev-pert} again, we deduce that
$$
\frac 1{2\nu^s_k} \leq \frac 1 {2\nu_k} + C \|\phi_s - \phi\|_{H^1_0(\Omega)}^\rho \leq \frac 1 {2\nu_k} \left(1  + 2^{\rho+1} \nu_k C \|v(s) - \phi\|_{H^1_0(\Omega)}^\rho \right) \ \mbox{ for } \ s \geq 0,
$$
and therefore, taking $s_0 > 0$ large enough if needed, we obtain \eqref{QGI} for $s \geq s_0$. We can take $s_0 = 0$, provided that $v(s)$ lies on the open ball in $H^1_0(\Omega)$ centered at $\phi$ of radius $\delta_0  \in (0,r_0/2)$ small enough for all $s \geq 0$. Finally, we remark that $\delta_0 \in (0,r_0/2)$ can be chosen so as to depend only on $d,\Omega,q,\mathcal{Q}(v_0), r_0,\nu_k$ and $\nu_{\ell-1}$. The proof is completed.
 \end{proof}

\begin{remark}\label{R:const}
According to the proof above, if $\mathcal{Q}(v_0)$ is bounded from above for any $v_0$ in a subset $B$ of $\mathcal{X}$, e.g., a small neighborhood of $\phi$ on $\mathcal{X}$, one can choose the constants $c_0$ and $\delta_0$ uniformly for $v_0\in B$ in Lemma \ref{L:QGIa}. Moreover, the constants $C_0$ and $C_1$ in Lemma \ref{L:re-exp} can also be chosen uniformly for $v_0\in B$. These observations are not important at all to prove Theorem \ref{T:rate}, but they will be taken into account to prove Corollary \ref{C:rate}.
\end{remark}

\section{Proofs of Theorem \ref{T:rate} and Corollary \ref{C:rate}}\label{S:pf}

We are now in a position to prove Theorem \ref{T:rate} and Corollary \ref{C:rate}. We first give a proof for a more general form of Theorem \ref{T:rate}.

\begin{theorem}[Rate of convergence under the Morse--Bott property]\label{T:rate-g}
Let $\Omega$ be a bounded $C^{1,1}$ domain in $\R^d$ with boundary $\partial \Omega$. Let $v = v(x,s)$ be a nonnegative energy solution to \eqref{v-pde}--\eqref{v-ic} and let $\phi = \phi(x)$ be a positive solution to \eqref{elliptic} such that $v(\cdot,s) \to \phi$ strongly in $H^1_0(\Omega)$ as $s \to \infty$. Assume that \eqref{H-int} holds, {\rm (}and then \eqref{H-int+} also holds{\rm )}. Then there exists a constant $C \geq 0$ such that
\begin{equation}\label{rate-g}
\|v(\cdot,s)-\phi\|_{H^1_0(\Omega)} \leq C \exp \left( -\frac{\nu_k}{q-1} s\right)
\end{equation}
for $s \geq 0$. Here $\nu_k$ stands for the least positive eigenvalue of the linearized operator $\Lphi$.
\end{theorem}

\begin{proof}
As in~\cite{A21}, testing \eqref{v-pde} with $\partial_s v(s)$, we can derive the energy inequality,
\begin{equation}\label{eneq}
\frac{4}{qq'} \left\| \partial_s (v^{q/2})(s) \right\|_{L^2({{\Omega}})}^2 + \frac{\d}{\d s} J(v(s)) \leq 0 \quad \mbox{ for a.e.~} s > 0
\end{equation}
(in fact, even equality holds; see~\cite{JiXi23}. See also~\cite{AS25}). Noting that
\begin{align*}
\partial_s (v^{q-1})(s) = \frac{2(q-1)}{q} v^{(q-2)/2}(s) \partial_s (v^{q/2})(s), 
\end{align*}
we derive from Lemma \ref{L:re-exp} that
\begin{align}
\int_{\Omega} \left| \partial_s (v^{q-1})(s) \right|^2 \phi_s^{2-q} \, \d x\nonumber
&= \frac{4(q-1)^2}{q^2} \int_\Omega \left( \frac{v(s)}{\phi_s} \right)^{q-2} \left| \partial_s (v^{q/2})(s) \right|^2 \, \d x\\
&\leq \frac{4(q-1)^2}{q^2} \left\| \frac{v(s)}{\phi_s} \right\|_{L^\infty(\Omega)}^{q-2} \left\| \partial_s (v^{q/2})(s) \right\|_{L^2({\Omega})}^2 \label{U2}
\end{align}
for a.e.~$s > s_*$. Thus we obtain
\begin{align*}
\frac{1}{q-1} \left\| \frac{v(s)}{\phi_s} \right\|_{L^\infty({{\Omega}})}^{2-q} \left\| \partial_s (v^{q-1})(s) \right\|_{\H_s'}^2 + \frac{\d}{\d s} J(v(s)) \leq 0 \quad \mbox{ for a.e.~} \ s > s_*.
\end{align*}
By virtue of Lemma \ref{L:QGIa}, since $J(\phi_s) = J(\phi)$ (see (i) of Lemma \ref{L:MB}), we have
\begin{align}
\MoveEqLeft{
\frac{2\nu_{k}}{q-1} \left\| \frac{v(s)}{\phi_s} \right\|_{L^\infty({{\Omega}})}^{2-q} 
{\left( 1 + c_0\|v(s)-\phi\|_{H^1_0({\Omega})}^{\rho} \right)}^{-1} \left( J(v(s)) - J(\phi) \right)
}\nonumber\\
&+ \frac{\d}{\d s} \left( J(v(s)) - J(\phi) \right) \leq 0
\quad \mbox{ for a.e.~} \ s > \tilde s_0\label{DI}
\end{align}
for some constant $\tilde s_0 \geq s_*$. By virtue of \eqref{re-bdd} of Lemma \ref{L:re-exp}, we find that $\|v(s)/\phi_s\|_{L^\infty({\Omega})} \leq C_0 + 1$ for $s \geq s_* = \log 2 + 1$. Hence
\begin{align*}
\MoveEqLeft{
\frac{2\nu_{k}}{q-1} (C_0+1)^{2-q} 
\left( 1 + c_0M^{\rho} \right)^{-1} 
\left( J(v(s)) - J(\phi) \right)
}\nonumber\\
&+ \frac{\d}{\d s} \left( J(v(s)) - J(\phi) \right) \leq 0
\quad \mbox{ for a.e.~} \ s > \tilde s_0,
\end{align*}
where $M := \sup_{s \geq 0} \|v(s)-\phi\|_{H^1_0(\Omega)} < \infty$. Therefore there exists a constant $\mu \in (0,\nu_{k}/(q-1))$ depending only on $d, \Omega, q, \mathcal{Q}(v_0), \nu_{k},\nu_{\ell-1}$ such that
$$
0 \leq J(v(s)) - J(\phi) \leq \big( J(v(\tilde s_0)) - J(\phi) \big) \e^{-2\mu s} \quad \mbox{ for all } \ s \geq \tilde s_0,
$$
which along with the fact that $J(v(s)) \leq J(v_0)$ for $s \geq 0$ yields
\begin{equation}\label{J-exp0}
0 \leq J(v(s)) - J(\phi) \leq \e^{2\mu \tilde s_0} \left( J(v_0) - J(\phi) \right) \e^{-2\mu s}
\quad \mbox{ for all } \ s \geq 0.
\end{equation}

In what follows, we may assume that $J(v(s))-J(\phi) > 0$ for $s \geq 0$. Indeed, if $J(v(s_0)) = J(\phi)$ for some $s_0 \geq 0$, then in view of \eqref{eneq} one can verify that $v(s) = v(s_0)$ for all $s \geq s_0$. From the assumption that $v(s) \to \phi$ in $H^1_0(\Omega)$ as $s \to \infty$, we obtain $v(s) = \phi$ for $s \geq s_0$, which gives the desired conclusion (see Theorem \ref{T:rate-g}). 

The following lemma implies the convergence rate of $v(s)$ in $H^1_0(\Omega)$ in view of \eqref{J-exp0}.

\begin{lemma}[Energy convergence implies Sobolev convergence]\label{L:H10-conv}
Let $\mu > 0$ and suppose that there exists a positive constant $C_\mu$ such that 
\begin{equation}\label{p:ass}
0 \leq J(v(s)) - J(\phi) \leq C_\mu^2 \e^{-2\mu s} \quad \mbox{ for } \ s \geq 0.
\end{equation}
Then the following hold true\/{\rm :} 
\begin{enumerate}
 \item[\rm (i)] There exist constants $K,s_1 > 0$ such that
$$
\|v(s)-\phi\|_{H^1_0(\Omega)} \leq K C_\mu \e^{-\mu s} \quad \mbox{ for } \ s \geq s_1.
$$
Hence $K$ depends only on $d, \Omega, q, \mathcal{Q}(v_0), \nu_{k},\nu_{\ell-1}$.
 \item[\rm (ii)] Furthermore, there exists a constant $\delta_1 \in (0,r_0/2)$ depending only on $d, \Omega, q, r_0,\nu_{k},\nu_{\ell-1}$ such that, if $\sup_{s \geq 0} \|v(s)-\phi\|_{H^1_0(\Omega)} < \delta_1$, then 
$$
\|v(s)-\phi\|_{H^1_0(\Omega)} \leq \tilde K \e^{\mu s_*}\left( C_\mu + \|v_0^{q/2}-\phi^{q/2}\|_{L^2(\Omega)} \right) \e^{-\mu s} \quad \mbox{ for } \ s \geq 0,
$$
where $s_* = \log 2 + 1$, for some constant $\tilde K$ depending only on $d, \Omega, q, \delta_1,\nu_{k},\nu_{\ell-1}$.
\end{enumerate}
\end{lemma}

\begin{proof}
We first prove (i). Recalling \eqref{eneq} and applying \eqref{U2} in a slightly different manner now than before, we see that
\begin{align*}
- \dfrac{\d}{\d s} J(v(s))
&\geq 
 \frac{4(q-1)}{q^2} \left\|\partial_s (v^{q/2})(s)\right\|_{L^2({\Omega})}^2 \\
&\geq \frac 2 q \left\| \frac{v(s)}{\phi_s} \right\|_{L^\infty({{\Omega}})}^{(2-q)/2} \left\|\partial_s (v^{q-1})(s)\right\|_{\H_s'} \left\|\partial_s (v^{q/2})(s)\right\|_{L^2(\Omega)}
\end{align*}
for a.e.~$s > s_*$. 
Thus Lemmas \ref{L:re-exp} and \ref{L:QGIa} again yield
\begin{align*}
\MoveEqLeft{
\frac2q(C_0+1)^{(2-q)/2} \sqrt{2 \nu_{k}} \left(1 + c_0 M^{\rho}\right)^{-1/2} \left(J(v(s))-J(\phi)\right)^{1/2} \left\|\partial_s (v^{q/2})(s)\right\|_{L^2(\Omega)}}\\
&\leq - \dfrac{\d}{\d s} \left( J(v(s)) - J(\phi) \right)
\quad \mbox{ for a.e. } s > \tilde s_0,
\end{align*}
where $M = \sup_{s \geq 0} \|v(s)-\phi\|_{H^1_0(\Omega)} < \infty$ and $\tilde s_0$ is given in \eqref{DI}, and hence, we have
\begin{align*}
\left\|\partial_s (v^{q/2})(s)\right\|_{L^2(\Omega)} 
&\leq - K_0 \dfrac{\d}{\d s} \left( J(v(s)) - J(\phi) \right)^{1/2}
\quad \mbox{ for a.e. } s > \tilde s_0
\end{align*}
for some positive constant $K_0$ depending only on $d,\Omega,q,\mathcal{Q}(v_0),\nu_{k},\nu_{\ell-1}$. Thus employing \eqref{p:ass}, we obtain
\begin{align}
\left\| v^{q/2}(s) - \phi^{q/2} \right\|_{L^2(\Omega)} 
&\leq 
\int^\infty_s \left\| \partial_s (v^{q/2})(s)\right\|_{L^2(\Omega)} \, \d s \nonumber\\
&\leq K_0 \left( J(v(s)) - J(\phi) \right)^{1/2}
\stackrel{\eqref{p:ass}}\leq K_0 C_\mu \e^{-\mu s} \label{AA}
\end{align}
for $s \geq \tilde s_0$. Furthermore, using an elementary inequality,
$$
0 < \frac{a^p-b^p}{a-b} \leq \min\{a^{p-1},b^{p-1}\} \quad \mbox{ for } \ a,b > 0, \ a\neq b, \ p \in (0,1),
$$
and putting $p = 2/q \in (0,1)$, we observe that
\begin{align*}
\left| v(x,s)-\phi(x) \right|
= \left| \beta\left( v^{q/2}(x,s) \right) - \beta\left( \phi^{q/2}(x) \right)\right| 
\leq \phi^{(2-q)/2}(x) \left| v^{q/2}(x,s) - \phi^{q/2}(x) \right|,
\end{align*}
where $\beta(r) := r^{2/q}$ for $r \geq 0$, for $x \in \Omega$ and $s > 0$. Thus we obtain
\begin{align}
 \int_{\Omega} |v(s)-\phi|^2 \phi^{q-2} \, \d x
\leq \int_{\Omega} \left| v^{q/2}(s) - \phi^{q/2} \right|^2 \, \d x
\stackrel{\eqref{AA}}\leq K_0^2 C_\mu^2 \e^{-2\mu s} \label{E-ec}
\end{align}
for $s \geq \tilde s_0$. In addition, as in \eqref{Jvs} and \eqref{Rs-est}, one can derive that
\begin{align}
J(v(s))-J(\phi)
&= \frac 1 2 \|\nabla v(s) - \nabla \phi\|_{L^2({\Omega})}^2 + \mathsf{c}_q \int_{\Omega} \phi^{q-1} \left( v(s)-\phi \right) \, \d x\nonumber\\
&\quad - \frac{\mathsf{c}_q}{q} \|v(s)\|_{L^q({\Omega})}^q + \frac{\mathsf{c}_q}{q} \|\phi\|_{L^q({\Omega})}^q\nonumber\\
&\geq \frac 1 2 \|\nabla v(s) - \nabla \phi\|_{L^2({\Omega})}^2 - \frac{\mathsf{c}_q}2 (q-1) \int_{\Omega} \phi^{q-2} | v(s)-\phi |^2 \, \d x\nonumber\\
&\quad - C_* \|\nabla v(s)-\nabla \phi\|_{L^2({\Omega})}^{2+\rho} \quad \mbox{ for all } \ s \geq 0,\label{AAA}
\end{align}
where $\rho = \min\{q-2,1\} > 0$. Here the constant $C_*$ above depends on $d, \Omega, q, \mathcal{Q}(v_0)$ only. Therefore combining the above with \eqref{p:ass} and \eqref{E-ec}, we conclude that
\begin{equation}\label{BBB}
\frac14 \|\nabla v(s) - \nabla \phi\|_{L^2({\Omega})}^2 \leq \left[ 1 + \frac{\mathsf{c}_q}2 (q-1) K_0^2 \right] C_\mu^2 \e^{-2\mu s} 
=: \frac14 K^2 C_\mu^2 \e^{-2\mu s} \quad \mbox{ for } \ s \geq s_1
\end{equation}
for some $s_1 \geq \tilde s_0$, where $s_1$ is large enough so that $C_* \|\nabla v(s)-\nabla \phi\|_{L^2({\Omega})}^\rho \leq 1/4$ for $s \geq s_1$. Hence $K$ depends only on $d, \Omega, q, \mathcal{Q}(v_0), \nu_{k},\nu_{\ell-1}$.

We next prove (ii). According to Lemma \ref{L:QGIa}, if $\sup_{s \geq 0} \|v(s)-\phi\|_{H^1_0(\Omega)}$ is small enough, then we can take $s_0 = 0$ in \eqref{QGI}, and hence, \eqref{AA} and \eqref{E-ec} hold for $s \geq \tilde s_0 = s_* = \log 2 + 1$. Hence we can also take $s_1 = s_*$ in the proof of (i) above. Furthermore, it follows from \eqref{eneq} that
\begin{align*}
\|v^{q/2}(s)-v_0^{q/2}\|_{L^2(\Omega)} 
&\leq \int^s_0 \|\partial_s (v^{q/2})(\sigma)\|_{L^2(\Omega)} \, \d \sigma\\
&\leq \sqrt{s} \left( \int^s_0 \|\partial_s (v^{q/2})(\sigma)\|_{L^2({\Omega})}^2 \, \d \sigma \right)^{1/2}\\
&\leq \sqrt{s} \left( - \frac{qq'}4 \int^s_0 \dfrac{\d}{\d s} J(v(\sigma)) \, \d \sigma \right)^{1/2}\\
&\leq \frac{\sqrt{s q q'}}2 \left(J(v_0)-J(v(s))\right)^{1/2}\\
&\leq \frac{\sqrt{s_* q q'}}2 \left(J(v_0)-J(\phi)\right)^{1/2}
\stackrel{\eqref{p:ass}}\leq q C_\mu
\end{align*}
for $s \in [0,s_*]$. Hence we see that
\begin{align*}
\|v^{q/2}(s)-\phi^{q/2}\|_{L^2(\Omega)} 
&\leq \|v^{q/2}(s)-v_0^{q/2}\|_{L^2(\Omega)} + \|v_0^{q/2}-\phi^{q/2}\|_{L^2(\Omega)}\\
&\leq q C_\mu + \|v_0^{q/2}-\phi^{q/2}\|_{L^2(\Omega)}
= : \tilde C_\mu \leq \tilde C_\mu \e^{\mu s_*} \e^{-\mu s}
\end{align*}
for $s \in [0,s_*]$. Therefore by \eqref{AA} we have
$$
\|v^{q/2}(s)-\phi^{q/2}\|_{L^2(\Omega)} \leq \max\{K_0,\e^{\mu s_*}\} \tilde C_\mu \e^{-\mu s} \leq (K_0 + 1) \e^{\mu s_*} \tilde C_\mu \e^{-\mu s} =: \tilde K_0 \tilde C_\mu \e^{-\mu s}
$$
for all $s \geq 0$. Accordingly, \eqref{E-ec} holds for all $s \geq 0$ with $C_\mu$ and $K_0$ replaced by $\tilde C_\mu$ and $\tilde K_0$, respectively. Therefore, by virtue of \eqref{AAA}, one can also derive \eqref{BBB} for all $s \geq 0$ with $C_\mu$ and $K_0$ replaced by $\tilde C_\mu$ and $\tilde K_0$, respectively, since we have already ensured that $\sup_{s \geq 0} \|v(s) - \phi\|_{H^1_0(\Omega)}$ is sufficiently small. We finally remark that 
$$
1 + \frac{\mathsf{c}_q}2 (q-1) \tilde K_0^2 
\leq \left[ 1 + \frac{\mathsf{c}_q}2 (q-1) (K_0+1)^2 \right] \e^{2\mu s_*}
=: \frac1{4q^2} \tilde K^2 \e^{2\mu s_*}
$$
(see \eqref{BBB}) and recall the definition of $\tilde{C}_\mu$ to obtain the conclusion. The proof is completed.
\end{proof}

Thus recalling \eqref{J-exp0} and exploiting (i) of Lemma \ref{L:H10-conv} above, we obtain
\begin{align*}
\|v(s)-\phi\|_{H^1_0(\Omega)} 
\leq K \e^{\mu s_0} \left( J(v_0) - J(\phi) \right)^{1/2} \e^{-\mu s} \quad \mbox{ for all } \ s \geq s_1
\end{align*}
for some constants $K,s_1 > 0$. Then recalling \eqref{DI} along with \eqref{re-exp} of Lemma \ref{L:re-exp}, as in~\cite{A21}, we can improve \eqref{J-exp0} as
$$
0 \leq J(v(s)) - J(\phi) \leq C \left( J(v_0) - J(\phi) \right) \exp \left(-\frac{2\nu_{k}}{q-1} s\right)
\quad \mbox{ for all } \ s \geq 0,
$$
which along with (i) of Lemma \ref{L:H10-conv} further implies
\begin{equation*}
\|v(s)-\phi\|_{H^1_0(\Omega)} \leq C \left( J(v_0) - J(\phi) \right)^{1/2} \exp \left(-\frac{\nu_{k}}{q-1} s\right)
\quad \mbox{ for all } \ s \geq s_1
\end{equation*}
for some constant $s_1 > 0$. \prf{Indeed, (even for energy solutions, not just classical ones) since the singular part $(\D H)_{\rm sing}$ of the Radon measure $\D H$ is nonpositive, as in~\cite[Section 5]{AS25}, we can deduce that
\begin{align*}
\e^{\frac{2\nu_k}{q-1} s} H(s) - H(0) &= \D \left( \e^{\frac{2\nu_k}{q-1} s} H(s) \right)((0,s])\\
&= \int^s_0 \dfrac{\d}{\d \sigma} \left( \e^{\frac{2\nu_k}{q-1} \sigma} H (\sigma) \right) \, \d \sigma + \e^{\frac{2\nu_k}{q-1} s} (\D H)_{\rm sing}((0,s])\\
&\leq \int^s_0 C \e^{-\mu \sigma} \e^{\frac{2\nu_k}{q-1} \sigma} H(\sigma) \, \d \sigma,
\end{align*}
which along with Gronwall's inequality gives
$$
\e^{\frac{2\nu_k}{q-1} s} H(s) \leq H(0) \exp \left( C \int^s_0 \e^{-\mu \sigma} \, \d \sigma \right)
$$
for all $s \geq 0$. Thus it follows that
$$
H(s) \leq H(0) \e^{\frac{C}{\mu}} \e^{-\frac{2\nu_k}{q-1} s}
$$
for all $s \geq 0$. }Thus \eqref{rate-g} follows (in particular for all $s \geq 0$) by means of \eqref{bound}.\footnote{Here we do not retain the dependence of the constant $C$ above on the initial discrepancy from $\phi$ when extending the estimate up to the initial time $s = 0$. On the other hand, it is discussed in Corollary \ref{C:rate}.} This completes the proof of Theorem \ref{T:rate-g}.
\end{proof}

\begin{proof}[Proof of Theorem {\rm \ref{T:rate}}]
Thanks to Proposition \ref{P:int-ker}, assumption \eqref{H-int} holds true for any least-energy solutions to \eqref{elliptic} in $\Omega = A_R$ for any $R \geq R_0$, where $R_0$ is given by (A). Hence, since the Morse index of any least-energy solution to \eqref{elliptic} is one, \eqref{nu} holds with $\ell = 2$ and $k = d+1$. Furthermore, \eqref{H-int+} also holds with $r_0 = \infty$ by rotation. Thus the assertion follows immediately from Theorem \ref{T:rate-g}.
\end{proof}

We next give a proof of Corollary \ref{C:rate} in a more general form.

\begin{proof}[Proof of Corollary {\rm \ref{C:rate}}]
Let $\phi > 0$ be a \emph{least-energy} solution to \eqref{elliptic} and let $\delta_0$ be the positive constant given in Lemma \ref{L:QGIa}. Thanks to~\cite[Theorem 1]{A16} (see also Lemma \ref{L:stbl} below), taking a constant $\delta> 0$ small enough, we have 
$$
\sup_{s \geq 0} \|v(s) - \phi\|_{H^1_0(\Omega)} < \delta_0,
$$
provided that $v_0 \in \mathcal X$ and $\|v_0 - \phi\|_{H^1_0(\Omega)} < \delta$. Hence one can take $s_0 = 0$ in Lemma \ref{L:QGIa}. We also recall Remark \ref{R:const}, which asserts that the constants appeared in Lemmas \ref{L:re-exp} and \ref{L:QGIa} can be chosen uniformly for such $v_0$. Therefore the assertion of the corollary can be proved from \eqref{J-exp0} by repeating the same argument as in the proof of Theorem \ref{T:rate-g} and by using (ii) of Lemma \ref{L:H10-conv} instead of (i).
\end{proof}

We close this section with the following proposition:

\begin{proposition}[Convergence to least-energy profiles]\label{P:LES-conv}
Let $\Omega$ be a bounded $C^{1,1}$ domain in $\R^d$. Then there exists $\delta_0 > 0$ such that for any $v_0 \in \mathcal{X}$ satisfying $J(v_0) < d_J + \delta_0$, where $d_J$ denotes the \emph{least-energy} of $J$ given by
$$
d_J = \inf\{J(\psi) \colon \psi \in H^1_0(\Omega) \mbox{ is a nontrivial solution to \eqref{elliptic}}\},
$$
the energy solution $v = v(x,s)$ of the Cauchy--Dirichlet problem \eqref{v-pde}--\eqref{v-ic} converges in $H^1_0(\Omega)$ to a least-energy solution $\phi$ of \eqref{elliptic} as $s \to \infty$.
\end{proposition}

\begin{proof}
We first recall that the minimum of $J$ over the phase set $\mathcal{X}$ is $d_J$ (see~\cite[Propositions 8 and 10]{AK13}). We first prove quasiconvergence, i.e., convergence of $v(s_k)$ to a least-energy solution of \eqref{elliptic} along some sequence $s_k \to \infty$. Suppose to the contrary that there exists a sequence $(v_{0,n})$ in $\mathcal{X}$ such that, for each $n \in \N$, $J(v_{0,n}) < d_J + 1/n$ and the energy solution $v_n = v_n(x,s)$ of \eqref{v-pde}--\eqref{v-ic} with $v_0 = v_{0,n}$ never converges in $H^1_0(\Omega)$ to any least-energy solution of \eqref{elliptic} along any sequence in time diverging to $\infty$. Then by Theorem \ref{T:ap}, for every $n \in \N$, there exist some sequence $s^n_k \to \infty$ and a nontrivial solution $\psi_n$ of \eqref{elliptic} such that $v_n(s^n_k) \to \psi_n$ in $H^1_0(\Omega)$ as $k \to \infty$. From the decrease of the energy $s \mapsto J(v_n(s))$, it follows that 
$$
d_J < J(\psi_n) \leq J(v_{0,n}) < d_J + 1/n,
$$
which also implies that $(\psi_n)$ is a Palais--Smale sequence of $J$, and therefore, $\psi_n$ converges in $H^1_0(\Omega)$ to a least-energy solution $\phi$ of \eqref{elliptic} as $n \to \infty$. However, due to the \L ojasiewicz--Simon gradient inequality (see~\cite[Theorem 6.2]{AS25},~\cite{AGG19}), we observe that
$$
0 = \|J'(\psi_n)\|_{H^{-1}(\Omega)} \geq \omega \left( J(\psi_n) - J(\phi) \right)^{1-\theta} > 0 \quad \mbox{ for } \ n \in \N \ \mbox{ large enough}
$$
for some $\omega > 0$ and $\theta \in (0,1/2]$ (cf.~\cite{FeiSim00}). This yields a contradiction. Moreover, the convergence (as $s \to \infty$) follows from the quasiconvergence due to~\cite{FeiSim00} (see also~\cite{AS25}). We have completed the proof.
\end{proof}

\section{Optimality of the convergence rate}\label{S:opt}

In this section, we prove Theorem \ref{T:opt} concerning the optimality of the convergence rate obtained in Theorem \ref{T:rate}. To this end, we first recall the stability result in~\cite{A16} for (possibly non-isolated) least-energy asymptotic profiles for fast diffusion. See \eqref{X} for the definition of the phase set $\mathcal X$.

\begin{lemma}[Stability of least-energy asymptotic profiles~\cite{A16}]\label{L:stbl}
Let $\Omega$ be a smooth bounded domain in $\R^d$. Let $\phi = \phi(x)$ be a least-energy solution to \eqref{elliptic}. Then for every $\vep > 0$ there exists a constant $\delta_\vep > 0$ such that, for every $v_0 \in \mathcal X$, the unique energy solution $v = v(x,s)$ to \eqref{v-pde}--\eqref{v-ic} for the initial datum $v_0$ satisfies
$$
\sup_{s \geq 0} \|v(s)-\phi\|_{H^1_0(\Omega)} < \vep,
$$
provided that $\|v_0 - \phi\|_{H^1_0(\Omega)} < \delta_\vep.$
\end{lemma}

The next lemma is concerned with the foliated Schwarz symmetry of least-energy solutions in annuli to \eqref{elliptic}. 

\begin{lemma}[Symmetry of least-energy solutions in annuli{~\cite[Theorem 3.1]{Pac02}}]\label{L:sym}
Let $\Omega = A_R$ be the annular domain in $\R^d$ defined by \eqref{annulus} for $R > 0$.
All positive least-energy solutions to \eqref{elliptic} are foliated Schwarz symmetric.
\end{lemma}

We recall that, when $\phi$ is \emph{foliated Schwarz symmetric} in $A_R$ with respect to the $x_d$-axis, there exists a function $\tilde \phi : [R,R+1] \times [0,\pi] \to \R$ such that $\phi(x) = \tilde{\phi}(r,\theta)$ with $r = |x|$ and $\theta = \arccos(e_d \cdot x/|x|)$ for $x \in A_R$ and $\tilde \phi$ is nonincreasing in $\theta$. We also remark that $\partial_\theta \tilde \phi \not\equiv 0$ if, in addition, $\phi$ is nonradial.

Throughout this section, $\M_\phi$ denotes the submanifold of $H^1_0(A_R)$ defined by \eqref{Mphi}. With the aid of Proposition \ref{P:nondeg}, we have the following

\begin{lemma}[$\M_\phi$ is isolated from the other profiles]\label{L:isola}
Suppose that {\rm (A)} is satisfied and let $R \geq R_0$, where $R_0$ is given by {\rm (A)}. Let $\phi > 0$ be a nonradial least-energy solution to \eqref{elliptic} in $\Omega = A_R$. Then the set $\M_\phi$ is isolated in $H^1_0(A_R)$ from the other solutions to \eqref{elliptic}.
\end{lemma}

\begin{proof}
Applying the \L ojasiewicz--Simon gradient inequality (see~\cite[Theorem 6.2]{AS25}), we can observe that $\M_\phi$ is isolated in $H^1_0(A_R)$ from all solutions having higher energies. Indeed, suppose to the contrary that there exists a sequence $(z_n)$ of nontrivial solutions to \eqref{elliptic} (i.e., critical points of $J$) such that $J(z_n) > J(\phi)$ for $n \in \N$ and $z_n \to \phi$ strongly in $H^1_0(A_R)$. Then the \L ojasiewicz--Simon gradient inequality (see~\cite[Theorem 6.2]{AS25}) yields
$$
0 = \|J'(z_n)\|_{H^{-1}(A_R)} \geq \omega \left(J(z_n)-J(\phi)\right)^{1-\theta}> 0 \quad \mbox{ for } \ n \in \N \ \mbox{ large enough}
$$
for some constants $\omega > 0$ and $\theta \in (0,1/2]$. This yields a contradiction. Hence it suffices to show that $\M_\phi$ is isolated in $H^1_0(A_R)$ from all the other least-energy solutions to \eqref{elliptic}. Suppose to the contrary that there exists a sequence $(\psi_n)$ of least-energy solutions to \eqref{elliptic} such that $\psi_n \not\in \M_\phi$ for $n \in \N$ and $\psi_n \to \phi$ strongly in $H^1_0(A_R)$. Setting $\phi_n \in \mathbb{P}_{\M_\phi}(\psi_n)$, one finds that $(\psi_n - \phi_n, \Omega_{p,q} \phi_n)_{H^1_0(A_R)} = 0$ for any $p,q \in \R^d$ be such that $|p| = |q| = 1$ and $p \cdot q = 0$. Moreover, setting $z_n := (\phi_n-\psi_n)/\|\phi_n-\psi_n\|_{H^1_0(A_R)} \in H^1_0(A_R)$, we have
\begin{align*}
 - \Delta z_n = \mathsf{c}_q \frac{\phi_n^{q-1}-\psi_n^{q-1}}{\phi_n-\psi_n} z_n \ \mbox{ in } A_R, \quad z_n = 0 \ \mbox{ on } \partial A_R.
\end{align*}
We note that $\|\phi_n-\phi\|_{H^1_0(A_R)} \leq \|\phi_n-\psi_n\|_{H^1_0(A_R)} + \|\psi_n-\phi\|_{H^1_0(A_R)} \leq 2 \|\psi_n-\phi\|_{H^1_0(A_R)} \to 0$, which also yields $\Omega_{p,q} \phi_n \to \Omega_{p,q} \phi$ in $L^2(A_R)$ for any $p,q \in \R^d$ with $|p| = |q| = 1$ and $p \cdot q = 0$. Passing to the limit as $n \to \infty$ and employing the elliptic regularity theory, we can deduce that
$$
z_n \to z \quad \mbox{ strongly in } C^2(\overline{A_R})
$$
and
$$
- \Delta z = \mathsf{c}_q (q-1) \phi^{q-2} z \ \mbox{ in } A_R, \quad z = 0 \ \mbox{ on } \partial A_R.
$$
Furthermore, it follows that $\|z\|_{H^1_0(A_R)} = 1$ and
\begin{align*}
0 = (z_n, \Omega_{p,q} \phi_n)_{H^1_0(A_R)}
 = \int_{A_R} (-\Delta z_n) \Omega_{p,q} \phi_n \, \d x 
 \to \int_{A_R} (-\Delta z) \Omega_{p,q} \phi \, \d x = (z,\Omega_{p,q} \phi)_{H^1_0(A_R)},
\end{align*}
which yields $\dim \ker \Lphi \geq d$ from the foliated Schwarz symmetry of $\phi$ (i.e., $\{\Omega_j \phi\}_{j=1,\ldots,d-1}$ are linearly independent). However, this contradicts the assertion of Proposition \ref{P:nondeg}. This completes the proof.
\end{proof}

We finally establish that

\begin{lemma}[Symmetric eigenfunctions]\label{L:e3}
Suppose that {\rm (A)} is satisfied. Let $R \geq R_0$, where $R_0$ is given by {\rm (A)}, and let $\phi > 0$ be a nonradial least-energy solution to \eqref{elliptic} in $\Omega = A_R$ such that $\phi$ is foliated Schwarz symmetric with respect to the $x_d$-axis. Then the following holds\/{\rm :}
\begin{enumerate}
 \item[(i)] In case $d \geq 3$, there exists an eigenfunction $e_{d+1}$ corresponding to the smallest positive eigenvalue $\nu_{d+1}$ of $\Lphi$ such that $e_{d+1}$ is symmetric with respect to the $(x_j,x_d)$-plane for any $j = 1,2$.
 \item[(ii)] In case $d = 2$, there exists an eigenfunction $e_3$ corresponding to the smallest positive eigenvalue $\nu_3$ of $\Lphi$ such that $e_3$ is even in $x_1$. 
\end{enumerate}
\end{lemma}

\begin{proof}
We first prove (i) for $d \geq 3$. Since $\phi$ is foliated Schwarz symmetric with respect to the $x_d$-axis, that is, $\phi = \tilde\phi(r,\theta)$ for $r = |x|$ and $\theta = \arccos(e_d \cdot x/|x|)$, the eigenvalue problem \eqref{ep} can be reduced into
\begin{equation*}
\left\{
\begin{aligned}
\left[ - \partial_r^2 - \frac{d-1}r \partial_r - \frac 1 {r^2} \left( \partial_\theta^2 + \frac{d-2}{\tan \theta} \partial_\theta - \frac{\ell(\ell + d - 3)}{\sin^2 \theta} \right) \right] \psi_{\ell,k} &= \lambda_{\ell,k} \tilde \phi^{q-2} \psi_{\ell,k} \ &&\mbox{ in } A_R,\\
\psi_{\ell,k} &= 0 \ &&\mbox{ on } \partial A_R
\end{aligned}
\right.
\end{equation*}
for $\ell \in \mathbb{Z}_{\geq 0}$ and $k \in \N$. For each $\ell \in \mathbb{Z}_{\geq 0}$, the eigenvalues $\{\lambda_{\ell,k}\}_{k\in\N}$ are arranged in nondecreasing order with respect to $k$, that is,
\begin{equation}\label{order-k}
\lambda_{\ell,k} \leq \lambda_{\ell,k+1} \quad \mbox{ for } \ \ell \in \mathbb{Z}_{\geq 0} \ \mbox{ and } \ k \in \N.
\end{equation}
Let $\{Y_{\ell,m}\}_{m = 1, \ldots, N(d-1,\ell)}$, where $N(d-1,\ell) := (2\ell+d-3) (\ell + d - 4)! /\{\ell! (d-3)!\}$, denote an orthonormal basis of spherical harmonics of degree $\ell$ on the $(d-2)$-dimensional unit sphere $\mathbb{S}^{d-2}$, that is, eigenfunctions corresponding to the eigenvalue $\ell(\ell + d-3)$ of the Laplace--Beltrami operator $- \Delta_{\mathbb{S}^{d-2}}$ on $\mathbb{S}^{d-2}$ (see, e.g.,~\cite{DaXu,SH}). Then each eigenvalue $\mu_j$ of \eqref{ep} coincides with $\lambda_{\ell,k}$ for some $\ell \in \mathbb{Z}_{\geq 0}$ and $k \in \N$, and moreover, $\psi_{\ell,k} Y_{\ell,m}$ becomes an eigenfunction corresponding to the eigenvalue $\mu_j = \lambda_{\ell,k}$ for $m \in \{ 1,\ldots, N(d-1,\ell) \}$. Noting variational representations of the eigenvalues,
$$
\lambda_{\ell,1} = \inf_{\psi \in X \setminus \{0\}} Q_\ell(\psi) \quad \mbox{ and } \quad \lambda_{\ell,2} = \inf_{\substack{E \subset X\\[1mm]\dim E = 2}} \sup_{\psi \in E \setminus \{0\}} \ Q_\ell(\psi),
$$
where $X := C^\infty_c((R,R+1)\times[0,\pi])$ and $Q_\ell : X \to [0,\infty]$ is given 
by
$$
Q_\ell(\psi) := \dfrac{\int^{R+1}_R \int^\pi_0 \left[ |\partial_r \psi|^2 + \frac 1 {r^2} |\partial_\theta \psi|^2 + \frac{\ell(\ell + d - 3)}{r^2 \sin^2 \theta} |\psi|^2 \right] r^{d-1} (\sin \theta)^{d-2} \, \d \theta \, \d r} 
{\int^{R+1}_R \int^\pi_0 \tilde\phi^{q-2} |\psi|^2 r^{d-1} (\sin \theta)^{d-2} \, \d \theta \, \d r}
$$
for $\ell \in \mathbb{Z}_{\geq 0}$, we find that 
$$
\lambda_{\ell,1} \leq \lambda_{\ell+1,1}, \quad \lambda_{\ell,2} \leq \lambda_{\ell+1,2} \quad \mbox{ for all } \ \ell \in \mathbb{Z}_{\geq 0}.
$$
Thanks to Proposition \ref{P:nondeg}, since eigenfunctions corresponding to the eigenvalue $\mu_2$ of \eqref{ep} are not axially symmetric with respect to the $x_d$-axis, we observe that
$$
\mu_1 = \lambda_{0,1} < \mu_2 = \cdots = \mu_d = \lambda_{1,1}.
$$
Furthermore, noting that
$$
\lambda_{0,2} \leq \lambda_{1,2} \quad \mbox{ and } \quad \lambda_{2,1} \leq \lambda_{3,1}
$$
(see also \eqref{order-k}), we obtain
$$
\mu_2 < \mu_{d+1} = \min\{\lambda_{0,2},\lambda_{2,1}\}. 
$$
In case $\mu_{d+1} = \lambda_{0,2}$, since $Y_{0,1}$ is constant, we observe that $e_{d+1} := \psi_{0,2}/\|\psi_{0,2}\|_{H^1_0(A_R)}$ is an eigenfunction corresponding to the eigenvalue $\mu_{d+1}$, and $e_{d+1}$ is then axially symmetric with respect to the $x_d$-axis. In case $\mu_{d+1} = \lambda_{2,1}$, there exists a spherical harmonic function $Y_2 = \sum_{m=1}^{N(d-1,2)} c_m Y_{2,m}$ of degree $2$ on $\mathbb{S}^{d-2}$ such that $Y_2$ is symmetric with respect to the $(x_j,x_d)$-plane for any $j = 1,2$ (e.g., $Y_2(\omega) = \omega_1^2 - \omega_2^2$ for $\omega \in \mathbb{S}^{d-2}$). Hence $e_{d+1} := \psi_{2,1} Y_2/\|\psi_{2,1} Y_2\|_{H^1_0(A_R)}$ turns to be an eigenfunction corresponding to the eigenvalue $\mu_{d+1}$ of \eqref{ep} such that $e_{d+1}$ is symmetric with respect to the same planes. Thus (i) has been proved.

We next prove (ii) for $d = 2$. 
In view of Proposition \ref{P:nondeg} with $d = 2$, we have $\ker \Lphi = \mathrm{span} \{\partial_\theta \phi\}$. Let $e_1$ and $e_2$ be a first and a second eigenfunctions of \eqref{ep}. Then $e_1$ and $e_2$ are multiples of $\phi$ and $\partial_\theta \phi$, respectively. Since $\phi$ is even in $x_1$, all eigenfunctions of \eqref{ep} are either even or odd in $x_1$. \prf{ Indeed, $H^1_0(A_R)$ can be expressed as the direct sum $L^+ \oplus L^-$ of closed subspaces $L^\pm = \{w \in H^1_0(A_R) \colon w(x_1,x_2) = \pm w(-x_1,x_2) \mbox{ for } (x_1,x_2) \in A_R \}$ of $H^1_0(A_R)$. Moreover, we see that $L^- = (L^+)^\perp$. For every $w \in H^1_0(A_R)$, there exist unique $w^+ \in L^+$ and $w^- \in L^-$ such that $w = w^+ + w^-$, and hence, $\|w\|_{H^1_0(A_R)}^2 = \|w^+\|_{H^1_0(A_R)}^2 + \|w^-\|_{H^1_0(A_R)}^2$. Furthermore, since $\phi$ is even, we find that $A(w) = A(w^+) + A(w^-) \in L^+ \oplus L^-$, where $A:H^1_0(A_R) \to H^1_0(A_R)$ is given by $A(w) = (-\Delta)^{-1}(\phi^{q-1}w)$ for $w \in H^1_0(A_R)$ (then eigenpairs of $A$ are $(1/\mu_j,e_j)$ for $j \in \N$). Then a Rayleigh quotient associated with $A$ reads,
$$
\frac{(A(w),w)_{H^1_0(A_R)}}{\|w\|_{H^1_0(A_R)}^2}
= \frac{\|w^+\|_{H^1_0(A_R)}^2}{\|w\|_{H^1_0(A_R)}^2} \frac{(A(w^+),w^+)_{H^1_0(A_R)}}{\|w^+\|_{H^1_0(A_R)}^2} + \frac{\|w^-\|_{H^1_0(A_R)}^2}{\|w\|_{H^1_0(A_R)}^2}\frac{(A(w^-),w^-)_{H^1_0(A_R)}}{\|w^-\|_{H^1_0(A_R)}^2}
$$
for $w = w^+ + w^- \in L^+ \oplus L^-$. Hence exploiting a minimax principle for eigenvalues of $A$, we conclude that each eigenvalue can be characterized in terms of the minimax value of the Rayleigh quotient over either $L^+$ or $L^-$. }Set $A_R^+ := \{x \in A_R \colon x_1 > 0\}$ and denote by $\mu^{\mathrm{D}}_j$ the $j$-th eigenvalue of \eqref{ep} with $A_R$ replaced by $A_R^+$. Moreover, we denote by $\mu^{\mathrm{N}}_j$ the $j$-th eigenvalue for
$$
-\Delta e_j = \mu^{\mathrm{N}}_j \phi^{q-2} e_j \ \mbox{ in } A_R^+, \quad e_j = 0 \ \mbox{ on } \partial A_R \cap [x_1 \geq 0], \quad \partial_{x_1} e_j = 0 \ \mbox{ on } A_R \cap [x_1 = 0].
$$
Since $e_1$ is even and $e_2$ is odd in $x_1$, we have $\mu_1 = \mu^{\mathrm{N}}_1$ and $\mu_2 = \mu^{\mathrm{D}}_1$. 
In view of the inequality $\mu^{\mathrm{N}}_2 \leq \mu^{\mathrm{D}}_2$, it follows that $\mu_3$ coincides with $\mu^{\mathrm{N}}_2$. Here we can check $\mu^{\mathrm{N}}_2 \leq \mu^{\mathrm{D}}_2$ by noting that $H^1_0(A_R^+) \subset X := \{w|_{A_R^+} \colon w \in H^1_0(A_R) \}$ and by using Rayleigh--Ritz's variational characterization of the eigenvalues:
\begin{align*}
\frac 1 {\mu^{\mathrm{D}}_2} &= \sup_{\substack{Y \subset H^1_0(A_R^+)\\[1mm]\dim Y = 2}} \inf \left\{ \int_{A_R^+} \phi^{q-2} w^2 \, \d x \colon w \in Y, \ \|\nabla w\|_{L^2(A_R^+)}=1 \right\}\\
&\leq \sup_{\substack{Y \subset X\\[1mm]\dim Y = 2}} \inf \left\{ \int_{A_R^+} \phi^{q-2} w^2 \, \d x \colon w \in Y, \ \|\nabla w\|_{L^2(A_R^+)}=1 \right\} = \frac 1 {\mu^{\mathrm{N}}_2},
\end{align*}
where $Y$ denotes a subspace of $H^1_0(A_R^+)$ or $X$ (see Lemma \ref{L:ev-var}). Therefore there exists an eigenfunction $e_3$ corresponding to the smallest positive eigenvalue $\nu_3 > 0$ of $\Lphi$ such that $e_3$ is even in $x_1$. Thus (ii) has been proved.
\end{proof}

We are now ready for proving Theorem \ref{T:opt}.

\begin{proof}[Proof of Theorem {\rm \ref{T:opt}}]
Let $\phi > 0$ be a nonradial least-energy solution to \eqref{elliptic}. Thanks to Lemma \ref{L:sym}, we may assume that $\phi$ is foliated Schwarz symmetric with respect to the $x_{d}$-axis. Furthermore, thanks to Lemma \ref{L:e3}, we can take an eigenfunction $e_{d+1}$ corresponding to the smallest positive eigenvalue $\nu_{d+1}$ of $\Lphi$ such that $\|e_{d+1}\|_{H^1_0(A_R)} = 1$, and for $d \geq 3$, $e_{d+1}$ is symmetric with respect to the $(x_j,x_d)$-plane for any $j = 1,2$ (for $d = 2$, $e_3$ is even in $x_1$). We set
$$
u_{0,\vep} := \phi + \vep e_{d+1} \quad \mbox{ for } \ \vep > 0.
$$
By virtue of Hopf's lemma along with the fact that $e_{d+1} \in C^1(\overline{A_R})$, we note that $u_{0,\vep} > 0$ in $A_R$ for $\vep > 0$ small enough. Moreover, we set
$$
v_{0,\vep} := t_*(u_{0,\vep})^{-1/(q-2)} u_{0,\vep} \in \mathcal{X} \quad \mbox{ for } \ \vep > 0.
$$
Putting $c_\vep := t_*(u_{0,\vep})^{-1/(q-2)}$, one observes that
\begin{equation}\label{v0e}
v_{0,\vep} = \phi + (c_\vep - 1) \phi + \vep c_\vep e_{d+1} \quad \mbox{ for } \ \vep > 0.
\end{equation}
Here we note that 
\begin{equation}\label{c:cep}
c_\vep = 1 + \mathrm{O}(\vep) \quad \mbox{ as } \ \vep \to 0_+
\end{equation}
 (see~\cite[Lemma 5.1]{AM23+} based on~\cite[Proposition 4]{AK13}), and hence, it follows that $v_{0,\vep} \to \phi$ strongly in $H^1_0(A_R)$ as $\vep \to 0_+$. Let $v_\vep = v_\vep(x,s)$ be the unique energy solution to \eqref{v-pde}--\eqref{v-ic} with $v_0 = v_{0,\vep}$. Then it follows from Lemma \ref{L:stbl} that $v_\vep(s)$ lies on a small neighbourhood of $\phi$ in $H^1_0(A_R)$ for any $s \geq 0$, provided that $\vep > 0$ is small enough. Moreover, thanks to Lemma \ref{L:isola}, there exists $\phi_\vep \in \M_\phi$ such that
$$
v_\vep(s) \to \phi_\vep \quad \mbox{ strongly in } H^1_0(A_R) \quad \mbox{ as } \ s \to \infty,
$$
provided that $\vep > 0$ is small enough. Since $v_{0,\vep}$ is symmetric with respect to the $(x_j,x_d)$-plane for any $j = 1,2$ when $d \geq 3$, and it is even in $x_1$ when $d = 2$, one can verify that $v_\vep(s)$ also inherits the same symmetry as $v_{0,\vep}$ for $s \geq 0$ (see~\cite[Lemma 2 in \S 3.3]{INdAM13}), and therefore $\phi_\vep$ also has the same symmetry as $v_{0,\vep}$. Now we claim that
\begin{claim}
The limit $\phi_\vep \in \M_\phi$ coincides with $\phi$, provided that $\vep > 0$ is small enough.
\end{claim}

Indeed, we first prove that \emph{only} $\phi$ and $\phi(-\,\cdot\,)$ among $\M_\phi$ enjoy the same symmetry as $\phi_\vep$. In case $d = 2$, it is obvious. In case $d \geq 3$, let $\Ro \in \mathrm{SO}(d)$ be such that $\psi := \phi(\Ro^{-1}\,\cdot\,)$ is symmetric with respect to the $(x_1,x_d)$-plane. Then denoting by $S_1$ the reflection with respect to the $(x_1,x_d)$-plane, i.e., $S_1 x = (x_1,-x_2,\ldots,-x_{d-1},x_d)$ for $x = (x_1,\ldots,x_d)$, we have $\psi(x) = \psi(S_1 x)$ for $x \in A_R$. Recall that $\phi(x) = \tilde \phi(|x|,\arccos(e_d \cdot x/|x|))$ for some function $\tilde \phi = \tilde \phi(r,\theta)$, which is strictly decreasing in $\theta$ (see~\cite[Remark 6.13]{DaPa}). Hence we note that 
$$
\psi(x) = \phi(\Ro^{-1}x) = \tilde \phi \left(|x|, \arccos (e_d\cdot \Ro^{-1} x/|x|) \right) = \tilde \phi \left(|x|, \arccos (\Ro e_d \cdot x/|x|) \right),
$$
which along with the fact that $\psi(x) = \psi(S_1 x)$ for $x \in A_R$ and the strict monotonicity of $\tilde \phi(r,\theta)$ in $\theta$ yields $\Ro e_d \cdot x = \Ro e_d \cdot S_1 x$ for $x \in A_R$, and hence $(\Ro e_d)_j = 0$ for $j = 2,\ldots,d-1$. Therefore if $\psi$ is symmetric with respect to the $(x_j,x_d)$-plane for $j = 1,2$, we can then reach $\Ro e_d = \pm e_d$, which implies that $\psi = \phi(\Ro^{-1}\,\cdot\,)$ coincides with either $\phi$ or $\phi(-\,\cdot\,)$. Since $\phi$ is distinct in $H^1_0(A_R)$ from $\phi(- \,\cdot\,)$, the claim readily follows.

Furthermore, thanks to Corollary \ref{C:rate}, we can take a constant $C$ independent of $\vep$ such that
\begin{align}
\|v_\vep(s)-\phi\|_{H^1_0(A_R)} 
&\leq C \left[ \left( J(v_{0,\vep}) - J(\phi) \right)^{1/2} + \|v_{0,\vep}^{q/2} - \phi^{q/2}\|_{L^2(A_R)} \right] \exp \left( - \frac{\nu_{d+1}}{q-1}s \right) \nonumber\\
&\leq C \vep \exp \left( - \frac{\nu_{d+1}}{q-1}s \right) \quad \mbox{ for } \ s \geq 0,\label{exp-stbl}
\end{align}
provided that $\vep > 0$ is small enough. Here we used the fact that
$$
0 \leq J(v_{0,\vep}) - J(\phi) \leq C\vep^2, \quad \|v_{0,\vep}^{q/2} - \phi^{q/2}\|_{L^2(A_R)} \leq C\vep,
$$
which can be derived from \eqref{v0e} and \eqref{c:cep} along with $q > 2$ (see \eqref{AAA}).

The rest of the proof runs as in~\cite[\S 5]{AM23+}. More precisely, we set
$$
w_\vep(s) := (-\Delta)^{-1} \left( v_\vep^{q-1}(s) - \phi^{q-1} \right) \quad \mbox{ for } \ s \geq 0,
$$
which solves
$$
\partial_s w_\vep(s) + v_\vep(s) - \phi = \mathsf{c}_q w_\vep(s) \ \mbox{ in } H^1_0(A_R) \quad \mbox{ for } \ s > 0.
$$
Set $W_\vep(s) := \mathbb{P}_{d+1}(w_\vep(s))$, where $\mathbb{P}_{d+1} : H^1_0(A_R) \to H^1_0(A_R)$ denotes the spectral projection onto the eigenspace corresponding to the eigenvalue $\mu_{d+1}$ of \eqref{ep}, for $s \geq 0$. Then we can derive the ODE,
$$
\partial_s W_\vep(s) + \mathbb{P}_{d+1}(v_\vep(s)-\phi) = \mathsf{c}_q W_\vep(s) \quad \mbox{ for } \ s > 0,
$$
which can also be rewritten\footnote{We used the relation that $\mathbb{P}_{d+1} = \mu_{d+1} \mathbb{P}_{d+1} \circ A$ with the operator $A : H^1_0(A_R) \to H^1_0(A_R)$ given by $A(w) = (-\Delta)^{-1}(\phi^{q-2} w)$ for $w \in H^1_0(A_R)$. Indeed, eigenpairs of $A$ are represented as $(1/\mu_j,e_j)$ for $j \in \N$ (see \eqref{ep}). Hence it follows from \eqref{R} that 
\begin{equation}\label{ooo}
\mathbb{P}_{d+1} \circ (-\Delta)^{-1} \mathcal{R}(v_\vep(s),\phi) = W_\vep(s) - \frac{q-1}{\mu_{d+1}} \mathbb{P}_{d+1} (v_\vep(s)-\phi).
\end{equation}} as
$$
\partial_s W_\vep(s) + \frac{\nu_{d+1}}{q-1} W_\vep(s) = \frac{\mu_{d+1}}{q-1} \mathbb{P}_{d+1} \circ (-\Delta)^{-1} \mathcal{R}(v_\vep(s),\phi) \ \mbox{ in } H^1_0(A_R) \quad \mbox{ for } \ s > 0,
$$
where $\mathcal{R}$ is given as
\begin{equation}\label{R}
\mathcal{R}(w,\phi) := |w|^{q-2}w - \phi^{q-1} - (q-1) \phi^{q-2} (w-\phi) \quad \mbox{ for } \ w \in H^1_0(A_R).
\end{equation}
Hence we have
\begin{align*}
 W_\vep(s) = \e^{-\frac{\nu_{d+1}}{q-1}s} W_\vep(0) + \frac{\mu_{d+1}}{q-1} \int^s_0 \e^{-\frac{\nu_{d+1}}{q-1}(s-\sigma)} \mathbb{P}_{d+1} \circ (-\Delta)^{-1} \mathcal{R}(v_\vep(\sigma),\phi) \, \d \sigma
\end{align*}
for $s \geq 0$. Moreover, we note that
\begin{equation}\label{R-est}
\|\mathcal{R}(w,\phi)\|_{L^{q'}(A_R)} \leq C \|w - \phi\|_{L^q(A_R)}^{1+\rho},
\end{equation}
(with $\rho = \min \{q-2,1\} > 0$) for some constant $C$ depending only on $R$, $q$ and being independent of $\vep \in (0,1)$. Noting that
\begin{align*}
 \left\| \mathbb{P}_{d+1} \circ (-\Delta)^{-1} \mathcal{R}(v_\vep(\sigma),\phi) \right\|_{H^1_0(A_R)} 
\leq C \|v_\vep(\sigma) - \phi \|_{H^1_0(A_R)}^{1+\rho}
\stackrel{\eqref{exp-stbl}}\leq C \vep^{1+\rho} \e^{-\frac{\nu_{d+1}}{q-1} (1+\rho) \sigma}
\end{align*}
and (by virtue of \eqref{v0e} and \eqref{c:cep} along with \eqref{R-est}) that
\begin{align*}
 \|W_\vep(0)\|_{H^1_0(A_R)} &\stackrel{\eqref{ooo}}\geq \frac{q-1}{\mu_{d+1}} \left\|\mathbb{P}_{d+1}(v_{0,\vep}-\phi)\right\|_{H^1_0(A_R)} - \left\|\mathcal{R}(v_{0,\vep},\phi)\right\|_{H^{-1}(A_R)} \\
&\geq \frac{q-1}{\mu_{d+1}} \left\|\vep c_\vep e_{d+1} \right\|_{H^1_0(A_R)} - C \left\|v_{0,\vep}-\phi\right\|_{H^1_0(A_R)}^{1+\rho}
\geq \frac{q-1}{2\mu_{d+1}} \vep c_\vep
\end{align*}
for $\vep > 0$ small enough, we deduce that
\begin{align*}
\MoveEqLeft{
 \|W_\vep(s)\|_{H^1_0(A_R)}
}\\
&\geq \e^{-\frac{\nu_{d+1}}{q-1} s} \|W_\vep(0)\|_{H^1_0(A_R)} - \frac{\mu_{d+1}}{q-1} \int^s_0 \e^{- \frac{\nu_{d+1}}{q-1}(s - \sigma)}\left\| \mathbb{P}_{d+1} \circ (-\Delta)^{-1} \mathcal{R}(v_\vep(\sigma),\phi) \right\|_{H^1_0(A_R)} \, \d \sigma\\
&\geq \frac{q-1}{2\mu_{d+1}} \vep c_\vep \e^{-\frac{\nu_{d+1}}{q-1} s} - \frac{\mu_{d+1}}{q-1} C \int^s_0 \e^{- \frac{\nu_{d+1}}{q-1}(s - \sigma)} \vep^{1+\rho} \e^{-\frac{\nu_{d+1}}{q-1}(1+\rho) \sigma} \, \d \sigma\\
&= \e^{-\frac{\nu_{d+1}}{q-1} s} \left( \frac{q-1}{2\mu_{d+1}} \vep c_\vep - \frac{\mu_{d+1}}{q-1} C \vep^{1+\rho} \int^s_0 \e^{\frac{\nu_{d+1}}{q-1}\sigma} \e^{-\frac{\nu_{d+1}}{q-1}(1+\rho) \sigma} \, \d \sigma \right)\\
&= \e^{-\frac{\nu_{d+1}}{q-1} s} \left[ \frac{q-1}{2\mu_{d+1}} \vep c_\vep + \frac{\mu_{d+1}}{\nu_{d+1} \rho} C \vep^{1+\rho} \left( \e^{-\frac{\nu_{d+1}}{q-1}\rho s} - 1\right) \right]\\
&\geq \e^{-\frac{\nu_{d+1}}{q-1} s} \left( \frac{q-1}{2\mu_{d+1}} \vep c_\vep - \frac{\mu_{d+1}}{\nu_{d+1} \rho} C \vep^{1+\rho} \right) > \frac{q-1}{4\mu_{d+1}} \vep c_\vep \e^{-\frac{\nu_{d+1}}{q-1} s} \quad \mbox{ for } \ s \geq 0
\end{align*}
for $\vep > 0$ small enough. Thus we have
\begin{align*}
\frac{q-1}{4\mu_{d+1}} \vep c_\vep \e^{-\frac{\nu_{d+1}}{q-1} s} 
&\leq \|W_\vep(s)\|_{H^1_0(A_R)} \leq \|w_\vep(s)\|_{H^1_0(A_R)}\\ 
&\leq \| v_\vep^{q-1}(s) - \phi^{q-1}(s) \|_{H^{-1}(A_R)}\\
&\stackrel{\eqref{R}}\leq \| \mathcal{R}(v_\vep(s), \phi) \|_{H^{-1}(A_R)} + (q-1) \left\| \phi^{q-2}(v_\vep(s)-\phi)\right\|_{H^{-1}(A_R)}\\
&\leq C \e^{-\frac{\nu_{d+1}}{q-1}(1+\rho)s} + C \left( \int_{A_R} |v_\vep(s)-\phi|^2 \phi^{q-2} \, \d x \right)^{1/2}.
\end{align*}
Here we also used the fact that
$$
\left\| w \right\|_{H^{-1}(A_R)}^2 \leq C \int_{A_R} |w(x)|^2 \phi^{2-q}(x) \, \d x \quad \mbox{ for } \ w \in L^2(A_R;\phi^{2-q}\, \d x),
$$
which follows from the duality of the embedding $H^1_0(A_R) \hookrightarrow L^2(A_R;\phi^{q-2}\, \d x)$. Thus the desired conclusion follows. This completes the proof.
\end{proof}

\section{A symmetry-induced jump in the convergence rate}\label{S:jump}

For each bounded domain $\Omega$ in $\R^d$ with boundary $\partial \Omega$, we consider the Dirichlet problem in $\Omega$,
\begin{equation}\label{eleq}
 - \Delta u = |u|^{q-2}u \ \mbox{ in } \Omega, \quad u = 0 \ \mbox{ on } \partial \Omega,
\end{equation}
under the assumption that $2 < q < 2^*$. Then the following lemma follows from~\cite[Theorem 1.2 and \S 4]{SaTe} (see also~\cite[Theorem 2.6]{BF21}).

\begin{proposition}\label{P:SaTe}
Let $\Omega_0$ be a bounded $C^{2,\alpha}$ domain in $\R^d$ and let $O \subset \R^d$ be an open bounded set such that $\overline{\Omega_0} \subset O$. Then there exists a sequence $(T_n)$ in $C^{2,\alpha}(O;O)$ such that
\begin{enumerate}
 \item[\rm (i)] there exists a constant $\vep_0 > 0$ such that $\det DT_n \geq \vep_0$ for every $n \in \N$\/{\rm ;}
 \item[\rm (ii)] $T_n$ converges to the identity $\mathrm{Id}$ in $C^{2,\alpha}(O;O)$\/{\rm ;} 
 \item[\rm (iii)] every positive solution to \eqref{eleq} in $\Omega_n := T_n(\Omega_0)$ is nondegenerate.
\end{enumerate}
\end{proposition}

In this section, we shall prove

\begin{proposition}\label{P:conv-ev}
Under the same assumptions as in Proposition {\rm \ref{P:SaTe}}, let $(T_n)$ be a sequence in $C^{2,\alpha}(\R^d,\R^d)$ such that {\rm (i)} and {\rm (ii)} of Proposition {\rm \ref{P:SaTe}} hold. Then the following {\rm (i)} and {\rm (ii)} hold true\/{\rm :}
\begin{enumerate}
 \item[(i)] Let $(u_n)$ be a sequence of least-energy solutions to \eqref{eleq} in $\Omega_n := T_n(\Omega_0)$. Then, up to a subsequence, $u_n$ converges to a least-energy solution $u_\infty$ of \eqref{eleq} in $H^1(\R^d)$. 
 \item[(ii)] Consider the following weighted eigenvalue problems\/{\rm :}
\begin{equation}\label{ev-0}
 - \Delta e_j = \mu_j |u_\infty|^{q-2} e_j \ \mbox{ in } \Omega_0, \quad e_j = 0 \ \mbox{ on } \partial \Omega_0
\end{equation}
and
\begin{equation}\label{ev-n}
 - \Delta e^n_j = \mu^n_j |u_n|^{q-2} e^n_j \ \mbox{ in } \Omega_n, \quad e^n_j = 0 \ \mbox{ on } \partial \Omega_n,
\end{equation}
where the eigenvalues are arranged in nondecreasing order, i.e., 
$$
0 < \mu_1 < \mu_2 \leq \cdots \leq \mu_j \leq \cdots \to \infty \quad \mbox{and} \quad 0 < \mu^n_1 < \mu^n_2 \leq \cdots \leq \mu^n_j \leq \cdots \to \infty,
$$
for $n \in \N$. Then for each $j \in \N$, the $j$-th eigenvalue $\mu_j^n$ of \eqref{ev-n} converges to the $j$-th eigenvalue $\mu_j$ of \eqref{ev-0} as $n \to \infty$.
\end{enumerate}
\end{proposition}

As a result, we can observe the following:

\begin{remark}[Jump in the convergence rate under certain domain perturbations]\label{R:jump}
Combining Theorem \ref{T:rate} with Propositions \ref{P:SaTe} and \ref{P:conv-ev} and noting that the constant $\mathsf{c}_q$ can be neglected by a simple scaling, we observe that the convergence rate of energy solutions to \eqref{v-pde}--\eqref{v-ic} may exhibit jumps under certain domain perturbations. Indeed, set $\Omega_0 = A_R$ for $R > 0$ large enough so that the conclusion of Theorem \ref{T:rate} holds. By (iii) of Proposition \ref{P:SaTe}, there exists a sequence of domains $\Omega_n := T_n(\Omega_0)$ such that all least-energy solutions to \eqref{elliptic} in $\Omega_n$ are nondegenerate for $n \in \N$. Let $\phi_n$ be a least-energy solution to \eqref{elliptic} in $\Omega_n$. In addition, let $v_n$ be an energy solution of \eqref{v-pde}--\eqref{v-ic} posed on $\Omega_n$ such that the initial datum $v_{0,n} \in \mathcal X$ is close enough to $\phi_n$ in $H^1_0(\Omega_n)$. Since $\phi_n$ is asymptotically stable on $\mathcal{X}$ (see~\cite[Theorem 2]{AK13}), $v_n(s)$ converges in $H^1_0(\Omega_n)$ to $\phi_n$ as $s \to \infty$. Denote by $\nu^n_j$ the $j$-th eigenvalue of the linearized operator $\mathcal{L}_{\phi_n}$ at $\phi_n$. Then noting that $\nu_1^n = 1-q < 0$ and $\nu_2^n > 0$ (indeed, the Morse index of least-energy solutions to \eqref{elliptic} is one), thanks to~\cite[Theorem 1.4]{A21}, we obtain
\begin{equation}\label{rate-Omn}
\|v_n(s)-\phi_n\|_{H^1_0(\Omega_n)} \leq C_n \e^{-\frac{\nu_2^n}{q-1}s} \quad \mbox{ for } \ s \geq 0,
\end{equation}
and actually, the above rate of convergence is also optimal in view of~\cite[Theorem 1.4]{AM23+}. Thanks to (i) of Proposition \ref{P:conv-ev}, it turns out that $\phi_n \to \phi$ in $H^1(\R^d)$, up to a subsequence, for some least-energy solution $\phi$ to \eqref{elliptic}. Then letting $\nu_j$ be the $j$-th eigenvalue of the linearized operator $\Lphi$ at $\phi$, thanks to Proposition \ref{P:nondeg}, we note that $\nu_1 = 1-q < 0 = \nu_2 = \cdots = \nu_d < \nu_{d+1}$. Moreover, by virtue of Proposition \ref{P:conv-ev}, we obtain 
$$
\nu_2^n = \mu^n_2 - \mathsf{c}_q (q-1) \to \nu_2 = \mu_2 - \mathsf{c}_q (q-1) = 0. 
$$
As a consequence, the rate of convergence in \eqref{rate-Omn} becomes slower as $n \to \infty$. On the other hand, Theorem \ref{T:rate} yields
$$
\|v(s)-\phi\|_{H^1_0(A_R)} \leq C \e^{-\frac{\nu_{d+1}}{q-1}s} \quad \mbox{ for } \ s \geq 0
$$
for all energy solutions $v = v(x,s)$ to \eqref{v-pde}--\eqref{v-ic}, provided that $v(s) \to \phi$ in $H^1_0(A_R)$ as $s \to \infty$ (according to the proof of Theorem \ref{T:opt}, such an energy solution $v$ actually exists). Therefore the rate of convergence appearing in \eqref{rate-Omn} exhibits a jump under the domain perturbations $\Omega_n \to \Omega_0 = A_R$ induced by the sequence of diffeomorphisms $(T_n)$.

In contrast, such phenomena occur neither in the linear diffusion case ($q = 2$) nor in the porous medium case ($1 < q < 2$). Indeed, in case $q = 2$, each solution to \eqref{u-pde}--\eqref{u-ic} can be represented as a Fourier series,
$$
u(x,t) = \sum_{j=1}^\infty c_j \e^{-\lambda_j t} e_j(x) \ \mbox{ with } \ c_j = (u_0,e_j)_{H^1_0(\Omega)},
$$
where $(\lambda_j,e_j)$ denotes the $j$-th eigenpair of the eigenvalue problem,
$$
- \Delta e_j = \lambda_j e_j \ \mbox{ in } \Omega, \quad e_j = 0 \ \mbox{ on } \partial \Omega,
$$
with $(e_j,e_k)_{H^1_0(\Omega)} = \delta_{jk}$ for $j,k \in \N$. Hence, from the fact that $\lambda_2 > \lambda_1$, the optimal convergence rate of the rescaled solution $\e^{\lambda_1 t} u(\cdot,t)$ to $e_1$ is $O(\e^{-(\lambda_2-\lambda_1)t})$ as $t \to \infty$. Therefore the convergence rate is continuous under domain perturbations, since the eigenvalues depend continuously on the domain. In case $1 < q < 2$, according to~\cite[Theorem 3]{ArPel81}, regardless of the geometry of $\Omega$, the optimal convergence rate of the rescaled function $(t+1)^{1/(2-q)} u(\cdot,t)$ to an asymptotic profile $\phi$ is $(t+1)^{-1}$ as $t \to \infty$. Finally, even for the fast diffusion case ($q > 2$), we remark that the optimal convergence rate obtained in~\cite{A21,AM23+} for \emph{nondegenerate} profiles is continuous under domain perturbations.
\end{remark}

In the rest of this section, we prove Proposition \ref{P:conv-ev}.

\begin{proof}[Proof of Proposition {\rm \ref{P:conv-ev}}]
We first prove the convergence of least-energy solutions. Denote by $d_{\Omega_n}$ and $d_{\Omega_0}$ the least energies of functionals $I_{\Omega_n}$ and $I_{{\Omega_0}}$, defined analogously to Section \ref{S:ker}, among nontrivial solutions to \eqref{eleq} with $\Omega$ replaced by $\Omega_n$ and ${\Omega_0}$, respectively. For each $u \in H^1_0(\Omega_n)$, we find that $w_u := u \circ T_n \in H^1_0({\Omega_0})$, i.e., $u = w_{u} \circ T_n^{-1}$. Then by Proposition \ref{P:SaTe}, it follows that
\begin{align*}
 d_{\Omega_n} &= \inf_{u \in H^1_0(\Omega_n) \setminus \{0\}} \left(\frac12-\frac1q\right) \left(\frac{\|\nabla u\|_{L^2(\Omega_n)}}{\|u\|_{L^q(\Omega_n)}}\right)^{2q/(q-2)}\\
&= \inf_{u \in H^1_0({\Omega_n}) \setminus \{0\}} \left(\frac12-\frac1q\right) \left(\frac{\| (DT_n^{-1} \circ T_n)^\top \nabla w_u |\mathrm{Jac} \, DT_n|^{1/2} \|_{L^2({\Omega_0})}}{\|w_u |\mathrm{Jac} \,DT_n|^{1/q} \|_{L^q(\Omega_0)}}\right)^{2q/(q-2)}\\
&\leq \left(\frac12-\frac1q\right) \left(\frac{\| (DT_n^{-1} \circ T_n)^\top \nabla w |\mathrm{Jac} \, DT_n|^{1/2} \|_{L^2({\Omega_0})}}{\|w |\mathrm{Jac} \,DT_n|^{1/q} \|_{L^q(\Omega_0)}}\right)^{2q/(q-2)} \\
&\to \left(\frac12-\frac1q\right) \left(\frac{\| \nabla w \|_{L^2({\Omega_0})}}{\|w \|_{L^q(\Omega_0)}}\right)^{2q/(q-2)} \ \mbox{ for } \ w \in H^1_0({\Omega_0}) \setminus \{0\},
\end{align*}
which implies that
\begin{equation}\label{rm:1}
\limsup_{n \to \infty} d_{\Omega_n} \leq d_{\Omega_0}.
\end{equation}
Let $(u_n)$ be a sequence of least-energy solutions to \eqref{eleq} in $\Omega_n$ for $n \in \N$. It then follows that
$$
\|\nabla u_n\|_{L^2(\R^d)}^2 = \|\nabla u_n\|_{L^2(\Omega_n)}^2 = \|u_n\|_{L^q(\Omega_n)}^q = \left(\frac12-\frac1q\right)^{-1} d_{\Omega_n},
$$
which along with \eqref{rm:1} yields, up to a subsequence,
$$
u_n \to u_\infty \quad \mbox{ weakly in } H^1(\R^d) \ \mbox{ and strongly in } L^r(\R^d) \ \mbox{ for } \ r \in [1,2^*),
$$
since $\Omega_n$ is a subset of the bounded set $O$ for any $n \in \N$. Moreover, as in the proof of Lemma \ref{L:d_AR-conv}, noting that, for any $\varphi \in C^\infty_c(\Omega_0)$ (resp., $\varphi \in C^\infty_c(\R^d \setminus \overline{\Omega_0})$), the support of $\varphi$ is included in $\Omega_n$ (resp., $\R^d \setminus \overline{\Omega_n}$) for $n \in \N$ large enough (see (ii) of Proposition \ref{P:SaTe}), we can verify that $u_\infty \in H^1_0(\Omega_0)$ is a weak solution to \eqref{eleq} in $\Omega_0$. Furthermore, by the Sobolev--Poincar\'e inequality $\|w\|_{L^q(O)} \leq C_{O} \|\nabla w\|_{L^2(O)}$ for $w \in H^1_0(O)$, we see that
$$
\|u_n\|_{L^q(O)}^2 \leq C_{O}^2 \|\nabla u_n\|_{L^2(O)}^2 = C_{O}^2 \|u_n\|_{L^q(O)}^q,
$$
which yields $\|u_n\|_{L^q(O)}^{q-2} \geq C_{O}^{-2} > 0$ for any $n \in \N$. Thus we infer that $u_\infty \neq 0$. Consequently, from \eqref{rm:1}, we obtain
\begin{align*}
 \limsup_{n \to \infty} \|\nabla u_n\|_{L^2(\R^d)}^2 \leq \left(\frac12 - \frac1q \right)^{-1} d_{\Omega_0} \leq \|\nabla u_\infty\|_{L^2(\R^d)}^2.
\end{align*}
By uniform convexity, we conclude that $u_n \to u_\infty$ strongly in $H^1(\R^d)$. Moreover, the limit $u_\infty$ fulfills $I_{\Omega_0}(u_\infty) = d_{\Omega_0}$, i.e., $u_\infty$ is a least-energy solution to \eqref{eleq} in $\Omega_0$. Furthermore, by Proposition \ref{P:SaTe}, we can also verify that
\begin{equation}\label{uT-conv}
 \|u_n \circ T_n - u_\infty\|_{H^1_0(O)} \to 0 
  \quad \mbox{ and } \quad \|u_n - u_\infty \circ T_n^{-1} \|_{H^1_0(O)} \to 0. 
\end{equation}

We now proceed to prove the convergence of the eigenvalue sequence $(\mu_j^n)$ as $n \to \infty$ for each $j \in \N$. To this end, we first recall that, for each $j \in \N$, the eigenvalues $\mu_j$ and $\mu^n_j$ of \eqref{ev-0} and \eqref{ev-n}, respectively, admit the following variational characterizations:
\begin{equation}\label{vc1}
\frac 1{\mu_j} = \sup_{\substack{Y \subset H^1_0({\Omega_0})\\\dim Y = j}} \inf_{\substack{w \in Y\\\|w\|_{H^1_0({\Omega_0})}=1}} \int_{{\Omega_0}} |u_\infty|^{q-2} w^2 \, \d x
\end{equation}
and
\begin{equation}\label{vc2}
\frac 1{\mu^n_j} = \sup_{\substack{Y \subset H^1_0({\Omega_n})\\\dim Y = j}} \inf_{\substack{w \in Y\\\|w\|_{H^1_0({\Omega_n})}=1}} \int_{{\Omega_n}} |u_n|^{q-2} w^2 \, \d x.
\end{equation}

We next establish the following lemma:

\begin{lemma}\label{L:ortho}
Let $j \in \N$ be fixed. Then there exists $n_0 \in \N$ such that the following hold true for all $n \geq n_0$\/{\rm :}
\begin{enumerate}
 \item[(i)] Let $Y$ be a $j$-dimensional subspace of $H^1_0({\Omega_0})$ spanned by an orthonormal system $\{e_i\}_{i=1}^j$. Then $\{e^n_i\}_{i=1}^j$, where $e^n_i := e_i \circ T_n^{-1}$, is linearly independent in $H^1_0(\Omega_n)$, that is, $Y_n := \mathrm{span} \{e^n_i\}_{i=1}^j$ is a $j$-dimensional subspace of $H^1_0(\Omega_n)$.
 \item[(ii)] Let $Y_n$ be a $j$-dimensional subspace of $H^1_0(\Omega_n)$ spanned by an orthonormal system $\{e_i^n\}_{i=1}^j$. Then $\{\tilde e^n_i\}_{i=1}^j$, where $\tilde e^n_i := e_i^n \circ T_n$, is linearly independent in $H^1_0({\Omega_0})$, that is, $\tilde Y_n := \mathrm{span} \{\tilde e^n_i\}_{i=1}^j$ is a $j$-dimensional subspace of $H^1_0({\Omega_0})$.
\end{enumerate}
\end{lemma} 

\begin{proof}
We prove (ii) only, as (i) can be proved similarly. For each $i,k \in \{1,\ldots,j\}$, we calculate
\begin{align*}
\lefteqn{
(\tilde e_i^n, \tilde e_k^n)_{H^1_0({\Omega_0})}
}\\
&= \int_{\Omega_0} ((D T_n)^\top (\nabla e_i^n \circ T_n)) \cdot ((D T_n)^\top (\nabla e_k^n \circ T_n)) \, \d x\\
&= \int_{\Omega_n} ( (D T_n \circ T_n^{-1})^\top \nabla e_i^n ) \cdot ( (D T_n \circ T_n^{-1})^\top \nabla e_k^n) |\det DT_n^{-1}| \, \d x\\
&= \int_{\Omega_n} ((D T_n \circ T_n^{-1})^\top \nabla e_i^n) \cdot ((D T_n \circ T_n^{-1})^\top \nabla e_k^n) (|\det DT_n^{-1}|-1) \, \d x\\
&\quad + \int_{\Omega_n} ((D T_n \circ T_n^{-1} - I)^\top \nabla e_i^n) \cdot ((D T_n \circ T_n^{-1} - I)^\top \nabla e_k^n) \, \d x\\
&\quad + \int_{\Omega_n} \nabla e_i^n \cdot ((D T_n \circ T_n^{-1} - I)^\top \nabla e_k^n) \, \d x + \int_{\Omega_n} ((D T_n \circ T_n^{-1} - I)^\top \nabla e_i^n) \cdot \nabla e_k^n \, \d x\\
&\quad + \int_{\Omega_n} \nabla e_i^n \cdot \nabla e_k^n \, \d x\\
&= \delta_{ik} + o_n(1)
\end{align*}
as $n \to \infty$. Here we denote by $I$ the unit $d \times d$ matrix and used the fact that $(e_i^n,e_k^n)_{H^1_0(\Omega_n)} = \delta_{ik}$. Hence we can take $n_0 \in \N$ large enough so that the matrix $\{(\tilde e_i^n,\tilde e_k^n)_{H^1_0(\Omega_n)}\}_{i,k=1,\ldots,d}$ is regular, and hence $\{\tilde e_i^n\}_{i=1}^j$ is linearly independent, provided that $n \geq n_0$. We also emphasize that $n_0$ can be chosen independently of the choice of $\{e_i^n\}_{i=1}^j$, since $\|e_i^n\|_{H^1_0(\Omega_n)} = 1$ for all $i \in \{ 1,\ldots,j \}$ and $n \in \N$.
\end{proof}

Let $Y_n$ be a $j$-dimensional subspace of $H^1_0(\Omega_n)$ spanned by $\{e^n_i\}_{i=1}^j$. Setting $\tilde Y_n$ as a $j$-dimensional subspace of $H^1_0({\Omega_0})$ spanned by $\{\tilde e_i^n\}_{i=1}^j$, where $\tilde e_i^n := e_i^n \circ T_n$ (see (ii) of Lemma \ref{L:ortho}), we derive from \eqref{vc1} that
\begin{align*}
\frac 1{\mu_j} 
&\geq \inf_{\substack{\tilde w \in \tilde Y_n\\ \|\tilde w\|_{H^1_0({\Omega_0})}=1}} \int_{{\Omega_0}} |u_\infty|^{q-2} \tilde w^2 \, \d x\\
&= \inf_{\substack{\tilde w \in \tilde Y_n\\ \|\tilde w\|_{H^1_0({\Omega_0})}=1}} \int_{\Omega_n} |u_\infty \circ T_n^{-1}|^{q-2} (\tilde w \circ T_n^{-1})^2 |\mathrm{Jac}\,DT_n^{-1}|\, \d x\\
&\geq \inf_{\substack{\tilde w \in \tilde Y_n\\ \|\tilde w\|_{H^1_0({\Omega_0})}=1}} \|\tilde w \circ T_n^{-1}\|_{H^1_0(\Omega_n)}^2 \inf_{\substack{w \in Y_n\\[0mm]\|w\|_{H^1_0(\Omega_n)}=1}} \int_{\Omega_n} |u_\infty \circ T_n^{-1}|^{q-2} w^2 |\mathrm{Jac}\,DT_n^{-1}|\, \d x.
\end{align*}
Since $\|w\|_{H^1_0(\Omega_n)} = 1$, we observe that
\begin{align*}
\MoveEqLeft{
\int_{\Omega_n} |u_\infty \circ T_n^{-1}|^{q-2} w^2 |\mathrm{Jac}\,DT_n^{-1}|\, \d x
}\\
&\geq \int_{\Omega_n} |u_n|^{q-2} w^2 \, \d x - \int_{\Omega_n} \left| |u_n|^{q-2} - |u_\infty \circ T_n^{-1}|^{q-2} \right| w^2 \, \d x\\
&\quad - \int_{\Omega_n} |u_\infty \circ T_n^{-1}|^{q-2} w^2 \left|1 - |\mathrm{Jac}\,DT_n^{-1}| \right|\, \d x\\
&\geq \int_{\Omega_n} |u_n|^{q-2} w^2 \, \d x - C \left\| u_n - u_\infty \circ T_n^{-1} \right\|_{L^{q}(O)}^\rho - C \left\| 1 - |\mathrm{Jac}\,DT_n^{-1}| \right\|_{L^\infty(O)},
\end{align*}
where $\rho := \min\{q-2,1\} > 0$, for some constant $C$ independent of $n$ and $w$. Moreover, recalling that $\|\tilde w\|_{H^1_0(\Omega_0)} = 1$, we can also prove that
\begin{align*}
\|\tilde w \circ T_n^{-1}\|_{H^1_0(\Omega_n)}^2
&= 
\int_{\Omega_0} |(D T_n^{-1} \circ T_n)^\top \nabla \tilde w|^2 |\det DT_n| \, \d x \\
&= 1 + \rho_n(\tilde w)
\end{align*}
for some reminder $\rho_n(\tilde w)$ satisfying
\begin{equation}\label{Rn}
|\rho_n(\tilde w)| \leq R_n := C\left( \left\| |\det D T_n| - 1 \right\|_{L^\infty(O)} + \left\|D T_n^{-1} \circ T_n - I \right\|_{L^\infty(O)} \right) \to 0,
\end{equation}
where $C$ is a constant independent of $\tilde w$ and $n$. Therefore combining these facts and taking the supremum with respect to $Y_n$, we obtain
\begin{align*}
\frac 1{\mu_j} \geq (1 - R_n) \left( \frac 1{\mu^n_j} - C \left\| u_n - u_\infty \circ T_n^{-1} \right\|_{L^{q}(O)}^\rho - C \left\| 1 - |\mathrm{Jac}\,DT_n^{-1}| \right\|_{L^\infty(O)} \right),
\end{align*}
which along with \eqref{uT-conv} yields
$$
\frac 1{\mu_j} \geq \limsup_{n \to \infty} \frac 1{\mu^n_j}.
$$
Exchanging the roles of $\mu_j$ and $\mu_j^n$ and using \eqref{vc2} and (i) of Lemma \ref{L:ortho}, we can similarly derive that
\begin{align*}
\frac 1{\mu_j^n} 
&\geq (1 - \tilde R_n) \left( \frac 1{\mu_j} - C \left\| u_\infty - u_n \circ T_n \right\|_{L^{q}(O)}^\rho - C \left\| 1 - |\mathrm{Jac}\,DT_n| \right\|_{L^\infty(O)} \right)
\end{align*}
for some sequence $(\tilde R_n)$ defined analogously to $(R_n)$ of \eqref{Rn}. It then follows from \eqref{uT-conv} that
$$
\liminf_{n \to \infty} \frac 1{\mu^n_j} \geq \frac 1{\mu_j}.
$$
Thus we conclude that
$$
\lim_{n \to \infty} \mu^n_j = \mu_j \quad \mbox{ for } \ j \in \N,
$$
which completes the proof.
\end{proof}


\appendix

\section{Proof of Lemma \ref{L:LapEV}}\label{apdx:comm}

We first give a proof of Lemma \ref{L:LapEV}.

\begin{proof}[Proof of Lemma {\rm \ref{L:LapEV}}]
Recall that $\lambda_R$ is characterized by 
$$
\lambda_R = \inf_{ w \in H^1_{0,\mathrm{rad}} (A_R) \setminus \{0\} } 
	\frac{ \left\| \nabla w \right\|_{L^2(A_R)}^2 }{ \left\| w \right\|_{L^2 (A_R)}^2 } 
	= \inf_{ w \in H^1_{0,\mathrm{rad}} (A_R) \setminus \{0\} } 
	\frac{ \int_R^{R+1} |w'(r)|^2 r^{d-1}\, \d r  }{ \int_R^{R+1} |w(r)|^2 r^{d-1}\, \d r  }.
$$
We compute the Rayleigh quotient for $w (r) = \sin (  \pi (r - R )  )$. We note that
$$
\begin{aligned}
\int_{R}^{R+1} |w(r)|^2 r^{d-1} \, \d r 
= \int_0^1 \sin^2 \left( \pi t \right) (R+t)^{d-1}\, \d t 
&= R^{d-1} \int_0^1 \sin^2( \pi t ) \left( 1 + \frac{t}{R} \right)^{d-1} \, \d t\\
&= R^{d-1} \left( 1 + O (R^{-1}) \right) \int_0^1 \sin^2 (\pi t) \, \d t \\
&= \frac{R^{d-1}}{2} \left( 1 + O (R^{-1}) \right) 
\end{aligned}
$$
and 
$$
\begin{aligned}
\int_{R}^{R+1} |w'(r)|^2 r^{d-1}\, \d r
= R^{d-1} \left( 1 + O(R^{-1}) \right) \pi^2 \int_0^1 \cos^2 (\pi t) \, \d t 
= \frac{R^{d-1}}{2} \left( 1 + O(R^{-1}) \right) \pi^2.
\end{aligned}
$$
It then follows that 
$$
 \lambda_R \leq \left( 1 + O(R^{-1}) \right) \pi^2. 
$$
	
As for the lower bound, the same change of variables leads to 
\begin{align*}
\int_R^{R+1} |w(r)|^2 r^{d-1}\, \d r 
&= \int_0^1 |w(R+t)|^2 (R+t)^{d-1}\, \d r\\
&= R^{d-1}\int_0^1 |w(R+t)|^2 \left( 1 + R^{-1} t \right)^{d-1}\, \d t
\end{align*}
and
$$
\int_R^{R+1} |w'(r)|^2 r^{d-1}\, \d r 
 = R^{d-1} \int_0^1 |w'(R+t)|^2 \left( 1 + R^{-1} t \right)^{d-1}\, \d t.
$$
Therefore it follows that
$$
\begin{aligned}
\lambda_R &= \inf_{ w \in H^1_{0,\mathrm{rad}} (A_R) \setminus \{0\} } 
 \frac{ \int_0^1 |w'(R+t)|^2 \left( 1 + R^{-1} t \right)^{d-1}\, \d t }
 { \int_0^1 |w(R+t)|^2 \left( 1 + R^{-1} t \right)^{d-1}\, \d t }.
 \end{aligned}
$$
In the one-dimensional case, it is known that 
$$
\pi^2 = \inf_{v\in H^1_0(0,1) \setminus \{0\}} \frac{\int_0^1 |v'(t)|^2 \, \d t}{\int_0^1 |v(t)|^2 \, \d t}.
$$
Hence, letting $R \to \infty$, we obtain 
$$
\begin{aligned}
\liminf_{R \to \infty} \lambda_R 
= \liminf_{R \to \infty} \inf_{ w \in H^1_{0,\mathrm{rad}} (A_R) \setminus \{0\} } \frac{ \int_0^1 |w'(R+t)|^2 \left( 1 + R^{-1} t \right)^{d-1} \, \d t }{ \int_0^1 |w(R+t)|^2 \left( 1 + R^{-1} t \right)^{d-1}\, \d t }
\geq \pi^2. 
\end{aligned}
$$
Here we also used the fact that $1 \leq (1 + R^{-1}t)^{d-1} \leq (1 + R^{-1})^{d-1}$ for $t \in (0,1)$. Therefore we conclude that $\lambda_R \to \pi^2$ as $R \to \infty$. The proof is completed. 
\end{proof}

\section{Proof of \eqref{quan-re}~\cite{BF21}}\label{A:S:smoothing}

In this section, we give a proof of \eqref{quan-re}, which is almost same as the original one in~\cite{BF21} but also with a very slight modification for adapting it to the requirements in the present paper. Let $\Omega$ be a bounded $C^{1,1}$ domain in $\R^{d}$ with boundary $\partial \Omega$. We consider the Cauchy--Dirichlet problem,
\begin{alignat}{3}
\partial_s (|v|^{q-2}v) &= \Delta v + \mathsf{c}_q |v|^{q-2}v \quad &&\mbox{ in } \Omega \times \R_+,\label{v-pde-a}\\
v &= 0 \quad &&\mbox{ on } \partial \Omega \times \R_+,\label{v-bc-a}\\
v(\cdot,0) &= v_0 \in \mathcal{X} \quad &&\mbox{ in } \Omega\label{v-ic-a}
\end{alignat}
as well as the stationary problem,
\begin{equation}\label{elliptic-a}
- \Delta \phi = \mathsf{c}_q |\phi|^{q-2}\phi \ \mbox{ in } \Omega, \quad \phi = 0 \ \mbox{ on } \partial \Omega.
\end{equation}

The following proposition is a very slightly modified version of~\cite[Theorem 4.1]{BF21}. In contrast to~\cite{BF21}, no additional relation between $v$ and $\phi$ is assumed here; this distinction will play a role in the proof of Lemma \ref{L:re-exp}. 

\begin{proposition}[Uniform estimates for the relative error]\label{P:rel-err}
Let $\Omega$ be a bounded $C^{1,1}$ domain in $\R^{d}$ with boundary $\partial \Omega$. Let $v = v(x,s)$ be a nonnegative energy solution to the Cauchy--Dirichlet problem \eqref{v-pde-a}--\eqref{v-ic-a} for some $v_0 \in \mathcal{X}$ such that $\sup_{s \geq {s_0}} \|v(s)\|_{L^\infty(\Omega)} \leq a$ with ${s_0} := \log 2$ for some constant $a < \infty$. Let $\phi = \phi(x)$ be a positive solution to \eqref{elliptic-a}. Set
$$
h(x,s) := \frac{v(x,s)}{\phi(x)} - 1 \quad \mbox{ for } \ x \in \Omega \ \mbox{ and } \ s > 0.
$$ 
Then there exists a constant $C \geq 0$ depending only on ${d}$, $\Omega$, $q$, $a$, $\|\phi\|_{H^1_0(\Omega)}$ such that
\begin{align*}
\|h(s)\|_{L^\infty(\Omega)} &\leq \frac{\e^{\frac{2}{q-2}\delta}}{\delta} \left[ C ( 1 + \delta ) \sup_{\sigma \geq s - \delta} \left\| v^{q-1}(\sigma) - \phi^{q-1} \right\|_{L^{1}(\Omega)}^{1/{d}} + \frac{\delta^2}{q-2} \right]
\end{align*}
for any $s > 0$ and $\delta > 0$ with $s - \delta \geq {s_0}$. In particular, if $v^{q-1}(s)$ converges to $\phi^{q-1}$ strongly in $L^1(\Omega)$ as $s \to \infty$, then
$$
h(s) \to 0 \quad\mbox{ uniformly in } \Omega \ \mbox{ as } \ s \to \infty.
$$
\end{proposition}

\begin{remark}[Bound for the supremum]
It is also noteworthy that, for any $q \in (2,2^*)$, the supremum $\sup_{s \geq {s_0}} \|v(s)\|_{L^\infty(\Omega)}$ with ${s_0} := \log 2$ can be estimate from above in terms of ${d}$, $q$ and the Rayleigh quotient $\mathcal{Q}(v_0) = \|\nabla v_0\|_{L^2(\Omega)}/\|v_0\|_{L^q(\Omega)}$. Indeed, recalling~\cite[Lemma 1]{A16}, we note that
\begin{align*}
 \|v(s)\|_{L^\infty(\Omega)} \leq C_0 \mathcal{Q}(v_0)^{\frac{4q}{\kappa(q-2)}}
\quad \mbox{ for all } \ s \geq {s_0} = \log 2,
\end{align*}
where $\kappa := 2{d} - {d} q + 2q > 0$ (by $q < 2^*$) and $C_0$ is a constant depending only on $d$, $q$. Furthermore, due to the elliptic regularity, for any solutions $\phi$ to \eqref{elliptic}, there exists a constant $C_1$ depending only on $d$, $\Omega$, $q$, $\|\phi\|_{H^1_0(\Omega)}$ such that
$$
\|\phi\|_{L^\infty(\Omega)} \leq C_1.
$$
\end{remark}


We establish several lemmas for proving Proposition \ref{P:rel-err}. 

\begin{lemma}\label{L:hint}
There exists a constant $C \geq 0$ depending only on ${d}$, $\Omega$, $q$, $a$, $\|\phi\|_{H^1_0(\Omega)}$ such that
\begin{align*}
\left| \int^{s_2}_{s_1} h(x,s) \, \d s\right|
&\leq C ( 1 + s_2-s_1 ) \sup_{\sigma \geq s_1} \left\| v^{q-1}(\sigma) - \phi^{q-1} \right\|_{L^{1}(\Omega)}^{1/{d}}
\end{align*}
for $s_2 \geq s_1 \geq {s_0} = \log 2$ and $x \in \Omega$.
\end{lemma}

\begin{proof}
The following identity plays a crucial role:
\begin{align*}
 - \Delta (h(s) \phi) = - \Delta (v(s) - \phi)
= \mathsf{c}_q v^{q-1}(s) - \partial_s \left(v^{q-1} \right)(s) - \mathsf{c}_q \phi^{q-1},
\end{align*}
which yields
$$
h(s) \phi = (-\Delta)^{-1} \left[ \mathsf{c}_q \left(v^{q-1}(s) - \phi^{q-1} \right) - \partial_s \left(v^{q-1}\right)(s) \right].
$$
Integrate both sides in $s$ over $(s_1,s_2)$ to see that
\begin{align*}
\int^{s_2}_{s_1} h(s) \phi \, \d s
&= \mathsf{c}_q \int^{s_2}_{s_1} (-\Delta)^{-1} \left(v^{q-1}(s)-\phi^{q-1}\right) \, \d s - (-\Delta)^{-1} \left( v^{q-1}(s_2) - v^{q-1}(s_1) \right)\\
&=: \mathsf{c}_q I_1 + I_2
\end{align*}
for $s_2 \geq s_1 \geq 0$. In what follows, we shall estimate the terms $I_1$, $I_2$ in terms of the Green function $G_\Omega = G_\Omega(x,y)$ associated with the Dirichlet problem,
\begin{equation}\label{Dir}
-\Delta w = g \ \mbox{ in } \Omega, \quad w = 0 \ \mbox{ on } \partial \Omega.
\end{equation}
By abuse of notation, we write $G_\Omega * g$ for the integral operator associated with the Green kernel $G_\Omega$, that is, we simply denote
$$
G_\Omega * g(x) := \int_\Omega G_\Omega(x,y) g(y) \, \d y \quad \mbox{ for } \ x \in \Omega,
$$
which solves \eqref{Dir} for $g \in L^p(\Omega)$.

We first estimate $I_1$ as follows:
\begin{align*}
 |I_1| \leq  \int^{s_2}_{s_1} \left| (-\Delta)^{-1} \left(v^{q-1}(s)-\phi^{q-1}\right)(x) \right| \, \d s =  \int^{s_2}_{s_1} \left| G_\Omega * \left(v^{q-1}(s)-\phi^{q-1}\right)(x) \right| \, \d s.
\end{align*}
Here we note that
\begin{align*}
\left| G_\Omega * \left(v^{q-1}(s)-\phi^{q-1}\right) (x)\right|
&\leq \left( \int_{B_r(x)} + \int_{\Omega \setminus B_r(x)} \right)\left| G_\Omega(x,y) \left( v^{q-1}(y,s) - \phi^{q-1}(y) \right) \right| \, \d y\\
&\leq \left\| v^{q-1}(s) - \phi^{q-1} \right\|_{L^\infty(\Omega)} \int_{B_r(x)} G_\Omega(x,y) \, \d y\\
&\quad + \left( \sup_{y \in \Omega \setminus B_r(x)} G_\Omega(x,y) \right) \int_\Omega \left| v^{q-1}(y,s) - \phi^{q-1}(y) \right| \, \d y,
\end{align*}
where $s \geq {s_0}$ and $B_r(x)$ denotes the ball centered at $x$ of radius $r > 0$. Moreover, we recall that $G_\Omega(x,y)$ complies with the following estimates: for ${d} \geq 3$, there exists a constant $K = K({d},\Omega) \geq 0$ such that
\begin{align*}
0 \leq G_\Omega(x,y) \leq \frac K {|x-y|^{{d}-2}} \left( \frac{\phi(x)}{|x-y|} \wedge 1\right) \left( \frac{\phi(y)}{|x-y|} \wedge 1\right)
\end{align*}
(see~\cite[Theorem 3.3, p.333]{GrWi82}). Here we used the fact that $d_{\partial \Omega}(x) := \mathrm{dist}(x,\partial \Omega) \asymp \phi(x)$. We shall treat only the case where ${d} \geq 3$.\footnote{As for ${d} = 2$, we have $G_\Omega(x,y) \lesssim \log (1 + d_{\partial \Omega}(x)d_{\partial \Omega}(y)/|x-y|^2)$ with $d_{\partial \Omega}(x) := \mathrm{dist}(x,\partial \Omega)$ (see~\cite[Theorem 6.23]{ChZha95} and~\cite[Lemma 2.1]{AACh09}). Hence for each $r > 0$ and $x,y \in \Omega$, dividing the computation into two cases, $|x-y| \leq d_{\partial \Omega}(y)/2$ (i.e., $d_{\partial \Omega}(y) \leq 2 d_{\partial \Omega}(x)$) and $|x-y| > d_{\partial \Omega}(y)/2$, one can derive that $G_\Omega(x,y) \lesssim \log (1 + 2 d_{\partial \Omega}(x)/|x-y|)$. Therefore it follows from the fact that $d_{\partial \Omega}(x) \asymp \phi(x)$ that
$$
\int_{B_r(x)} G_\Omega(x,y) \, \d y \lesssim \phi(x) r, \quad \sup_{y \in \Omega \setminus B_r(x)} G_\Omega(x,y) \lesssim \phi(x) r^{-1}.
$$
The rest of the proof proceeds as before. The case where ${d} = 1$ is simpler since the Green function is bounded.} We observe that
$$
\int_{B_r(x)} G_\Omega(x,y) \, \d y \leq K \phi(x) \int_{B_r(0)} \frac 1 {|\zeta|^{{d}-1}} \, \d \zeta \leq C \phi(x) r, \quad
\sup_{y \in \Omega \setminus B_r(x)} G_\Omega(x,y) \leq K \phi(x) r^{1-{d}}.
$$
Therefore we have
\begin{align*}
\lefteqn{
\left| G_\Omega * \left(v^{q-1}(s)-\phi^{q-1}\right)(x) \right|
}\\
&\leq C \phi(x) r \left( \sup_{\sigma \geq {s_0}} \|v(\sigma)\|_{L^\infty(\Omega)}^{q-1} + \|\phi\|_{L^\infty(\Omega)}^{q-1} \right) + K \phi(x) r^{1-{d}} \left\| v^{q-1}(s) - \phi^{q-1} \right\|_{L^{1}(\Omega)}.
\end{align*}
Optimizing the right-hand side in $r$, we obtain
$$
\left| G_\Omega * \left(v^{q-1}(s)-\phi^{q-1}\right)(x) \right|
\leq C \phi(x) \left\| v^{q-1}(s) - \phi^{q-1} \right\|_{L^{1}(\Omega)}^{1/{d}},
$$
where $C$ depends only on ${d}$, $\Omega$, $q$, $a$ and $\|\phi\|_{H^1_0(\Omega)}$. Hence
\begin{align*}
 |I_1| \leq C  \phi(x) \left( \sup_{s \geq s_1} \left\| v^{q-1}(s) - \phi^{q-1} \right\|_{L^{1}(\Omega)}^{1/{d}} \right) (s_2-s_1)
\end{align*}
for $s_2 \geq s_1 \geq {s_0}$ and $x \in \Omega$.

We move on to the estimate for $I_2$. We similarly observe that
\begin{align*}
 |I_2| &\leq \int_\Omega G_\Omega(x,y) \left| v^{q-1}(y,s_2) - v^{q-1}(y,s_1) \right| \, \d y\\
&\leq 2 \left( \sup_{s \geq {s_0}} \|v(s)\|_{L^\infty(\Omega)}^{q-1} \right) C \phi(x) r + K \phi(x) r^{1-{d}} \left\|v^{q-1}(s_2) - v^{q-1}(s_1) \right\|_{L^{1}(\Omega)}
\end{align*}
for $s_2 \geq s_1 \geq {s_0}$ and $x \in \Omega$. Optimizing the right-hand side in $r$, we have
$$
|I_2| \leq C \phi(x) \left\|v^{q-1}(s_2) - v^{q-1}(s_1) \right\|_{L^{1}(\Omega)}^{1/{d}}
$$
for $s_2 \geq s_1 \geq {s_0}$ and $x \in \Omega$. Consequently, dividing both sides by $\phi(x) > 0$, we obtain
\begin{align*}
\left| \int^{s_2}_{s_1} h(x,s) \, \d s\right|
&\leq C  \left( \sup_{s \geq s_1} \left\| v^{q-1}(s) - \phi^{q-1} \right\|_{L^{1}(\Omega)}^{1/{d}} \right) (s_2-s_1) + C \left\|v^{q-1}(s_2) - v^{q-1}(s_1)\right\|_{L^{1}(\Omega)}^{1/{d}}\\
&\leq C (1 + s_2-s_1) \sup_{\sigma \geq s_1} \left\| v^{q-1}(\sigma) - \phi^{q-1} \right\|_{L^{1}(\Omega)}^{1/{d}}
\end{align*}
for $s_2 \geq s_1 \geq {s_0}$ and $x \in \Omega$. Thus the proof has been completed.
\end{proof}

We next set up the following

\begin{lemma}[B\'enilan--Crandall estimates]
It holds that
\begin{equation}\label{h-ineq}
h(x,s) \leq h(x,\sigma) \e^{\frac{2}{q-2}(s-\sigma)} +  \frac{2}{q-2}(s-\sigma) \e^{\frac{2}{q-2}(s-\sigma)}
\end{equation}
for $s \geq \sigma \geq {s_0} = \log 2$ and $x \in \Omega$.
\end{lemma}

\begin{proof}
Recalling the so-called B\'enilan--Crandall estimate (see, e.g.,~\cite[Lemma 8.1]{Vazquez},~\cite[Theorem 3.7]{AS25}),
$$
\partial_t \left(u^{q-1}\right)(t) \leq \frac{\mathsf{c}_q}t u^{q-1}(t) \quad \mbox{ for } \ t > 0,
$$
where $u$ is an energy solution to \eqref{u-pde}--\eqref{u-ic} and $u$ and $v$ are related by \eqref{v} (see also \eqref{X}), we have
\begin{align*}
\partial_s \left( v^{q-1} \right) (s) \leq \frac{\mathsf{c}_q}{1-\e^{-s}} v^{q-1}(s) \leq 2 \mathsf{c}_q v^{q-1}(s) \quad \mbox{ for } \ s \geq {s_0} = \log 2,
\end{align*}
which yields
$$
v(s) \leq v(\sigma) \e^{\frac{2}{q-2}(s-\sigma)} \quad \mbox{ for } \ s \geq \sigma \geq {s_0}.
$$
Thus we obtain
\begin{align*}
h(s) + 1 = \frac{v(s)}{\phi} \leq \frac{v(\sigma)}{\phi} \e^{\frac{2}{q-2}(s-\sigma)} = \left( h(\sigma)+1 \right) \e^{\frac{2}{q-2}(s-\sigma)} \ \mbox{ for } \ s \geq \sigma \geq {s_0}.
\end{align*}
Noting that $\e^a - 1 \leq a \e^a$ for $a \geq 0$, we conclude that
\begin{align*}
h(s) \leq h(\sigma) \e^{\frac{2}{q-2}(s-\sigma)} +  \frac{2}{q-2}(s-\sigma) \e^{\frac{2}{q-2}(s-\sigma)}  \ \mbox{ for } \ s \geq \sigma \geq {s_0}.
\end{align*}
This completes the proof.
\end{proof}

Using the lemma above, we can prove

\begin{lemma}\label{L:h-hint}
It holds that
\begin{align*}
 |h(x,s)| &\leq \frac{\e^{\frac{2}{q-2}\delta}}{\delta} \left( \left| \int^{s+\delta}_s h(x,\sigma) \, \d \sigma \right| + \left| \int^{s}_{s-\delta} h(x,\sigma) \, \d \sigma \right| + \frac{\delta^2}{q-2} \right)
\end{align*}
for $\delta > 0$, $s \geq {s_0} + \delta$ and $x \in \Omega$.
\end{lemma}

\begin{proof}
Let $s_2 \geq s_1 \geq {s_0}$. Substituting $\sigma = s_1$ to \eqref{h-ineq} and integrate both sides in $s$ over $(s_1,s_2)$, we find that
\begin{align*}
\int^{s_2}_{s_1} h(s) \, \d s
 \leq \frac{q-2}2 h(s_1) \left[ \e^{\frac{2}{q-2}(s_2-s_1)} - 1 \right] + \frac{1}{q-2}(s_2-s_1)^2 \e^{\frac{2}{q-2}(s_2-s_1)}.
\end{align*}
Since $\e^{a}-1 \geq a$ for $a \geq 0$, we have
\begin{align*}
 h(s_1) &\geq \frac 2 {q-2} \left[ \e^{\frac{2}{q-2}(s_2-s_1)} - 1 \right]^{-1} \left( - \left| \int^{s_2}_{s_1} h(s) \, \d s \right| - \frac{1}{q-2}(s_2-s_1)^2 \e^{\frac{2}{q-2}(s_2-s_1)} \right)\\
&\geq - \frac{1}{s_2-s_1}\left( \left| \int^{s_2}_{s_1} h(s) \, \d s \right| + \frac{1}{q-2}(s_2-s_1)^2 \e^{\frac{2}{q-2}(s_2-s_1)} \right)
\end{align*}
for $s_2 \geq s_1 \geq {s_0}$. On the other hand, it follows from \eqref{h-ineq} that
\begin{equation}\label{wxwx}
h(\sigma) \geq h(s) \e^{-\frac{2}{q-2}(s-\sigma)} - \frac{2}{q-2}(s-\sigma)
\end{equation}
for $s \geq \sigma \geq {s_0}$. Now setting $s = s_2$ in \eqref{wxwx} and integrate both sides in $\sigma$ over $(s_1,s_2)$, we deduce that
\begin{align*}
\int^{s_2}_{s_1} h(\sigma) \, \d \sigma
 \geq \frac{q-2}2 h(s_2) \left[ 1 - \e^{-\frac{2}{q-2}(s_2-s_1)} \right] - \frac{1}{q-2}(s_2-s_1)^2.
\end{align*}
Therefore using the relation that $1 - \e^{-a} \geq a \e^{-a}$ for $a \geq 0$, we obtain
\begin{align*}
 h(s_2) &\leq \frac 2 {q-2} \left[ 1 - \e^{-\frac{2}{q-2}(s_2-s_1)} \right]^{-1}
 \left( \left| \int^{s_2}_{s_1} h(\sigma) \, \d \sigma \right| + \frac{1}{q-2}(s_2-s_1)^2 \right)\\
&\leq \frac{\e^{\frac{2}{q-2}(s_2-s_1)}}{s_2-s_1} \left( \left| \int^{s_2}_{s_1} h(\sigma) \, \d \sigma \right| + \frac{1}{q-2}(s_2-s_1)^2 \right)
\end{align*}
for $s_2 \geq s_1 \geq {s_0}$. Now, let $s > 0$ and $\delta > 0$ satisfy $s - \delta \geq {s_0}$. Then we conclude that
\begin{align*}
- \frac{1}{\delta}\left( \left| \int^{s+\delta}_{s} h(\sigma) \, \d \sigma \right| + \frac{\delta^2}{q-2} \e^{\frac{2}{q-2}\delta} \right)
 \leq  h(s) \leq \frac{\e^{\frac{2}{q-2}\delta}}{\delta} \left( \left| \int^{s}_{s-\delta} h(\sigma) \, \d \sigma \right| + \frac{\delta^2}{q-2} \right),
\end{align*}
whence the desired conclusion follows.
\end{proof}

Combining all these facts, we are in a position to verify Proposition \ref{P:rel-err}.

\begin{proof}[Proof of Proposition {\rm \ref{P:rel-err}}]
Let $s > 0$ and $\delta > 0$ satisfy $s-\delta \geq {s_0}$. Combining Lemma \ref{L:h-hint} with Lemma \ref{L:hint}, we immediately obtain
\begin{align*}
\|h(s)\|_{L^\infty(\Omega)} &\leq \frac{\e^{\frac{2}{q-2}\delta}}{\delta} \left[ C ( 1 + \delta ) \sup_{\sigma \geq s - \delta} \left\| v^{q-1}(\sigma) - \phi^{q-1} \right\|_{L^{1}(\Omega)}^{1/d} + \frac{\delta^2}{q-2} \right].
\end{align*}
We next prove the uniform convergence of the relative error $h(s)$. From the assumption that
$$
\rho(s) := \sup_{\sigma \geq s - 1} \left\| v^{q-1}(\sigma) - \phi^{q-1} \right\|_{L^{1}(\Omega)} \to 0 \quad \mbox{ as } \ s \to \infty,
$$
for any $s \geq {s_0} + 1$ large enough, it follows that $\delta := \rho(s)^{1/(2d)} \in (0,1)$. Then we have
\begin{align*}
\|h(s)\|_{L^\infty(\Omega)}
\leq \frac{\e^{\frac{2}{q-2}\rho(s)^{1/(2d)}}}{\rho(s)^{1/(2d)}} \left[C \left( 1 +  \rho(s)^{1/(2d)} \right) \rho(s)^{1/d} + \frac{\rho(s)^{1/d}}{q-2} \right]\to 0 \quad \mbox{ as } \ s \to \infty,
\end{align*}
which completes the proof.
\end{proof}

\section{Morse--Bott lemma}\label{A:MB}

In this section, we give a proof of Lemma \ref{L:MB}.

\begin{proof}[Proof of Lemma {\rm \ref{L:MB}}]
Suppose to the contrary that there exists a sequence $(m_n)$ in $\M_\phi$ such that 
\begin{equation}\label{a01}
\|m_n-\phi\|_{H^1_0(\Omega)} \to 0 \quad \mbox{ and } \quad \ker \mathcal{L}_{m_n} \supsetneq T_{m_n} \M_\phi.
\end{equation}
Then for each $n \in \N$, there exists $\psi_n \in \ker \mathcal{L}_{m_n}$ such that 
\begin{equation}\label{a02}
\|\psi_n\|_{L^q(\Omega)} = 1 \quad \mbox{ and } \quad (\psi_n,v)_{H^1_0(\Omega)} = 0 \ \mbox{ for all } \ v \in T_{m_n} \M_\phi.
\end{equation}
Since $m_n \in \M_\phi$, we find from (MB) that
$$
- \Delta m_n = \mathsf{c}_q m_n^{q-1} \ \mbox{ in } \Omega, \quad m_n = 0 \ \mbox{ on } \partial \Omega, \quad \|m_n\|_{L^q(\Omega)}^q = \|\phi\|_{L^q(\Omega)}^q.
$$
Hence as $\mathcal{L}_{m_n}\psi_n = 0$ in $H^{-1}(\Omega)$, we have
$$
\|\nabla \psi_n\|_{L^2(\Omega)}^2 = \mathsf{c}_q (q-1) \int_\Omega m_n^{q-2} \psi_n^2 \, \d x \leq \mathsf{c}_q (q-1) \|m_n\|_{L^q(\Omega)}^{q-2} = \mathsf{c}_q (q-1) \|\phi\|_{L^q(\Omega)}^{q-2}
$$
for all $n \in \N$. Then there exists $\psi_\infty \in H^1_0(\Omega)$ such that, up to a subsequence, $\psi_n \to \psi$ weakly in $H^1_0(\Omega)$ and strongly in $L^q(\Omega)$. Thus recalling \eqref{a01} and \eqref{a02}, we obtain 
$$
\mathcal{L}_\phi \psi_\infty = 0 \ \mbox{ in } H^{-1}(\Omega) \quad \mbox{ and } \quad \|\psi_\infty\|_{L^q(\Omega)} = 1.
$$
Moreover, since $\M_\phi$ is a submanifold of $H^1_0(\Omega)$, we can deduce that\footnote{Indeed, set $p := \dim \M_\phi < \infty$. We can take $\delta_0 > 0$ and $p$ continuous vector fields $X_j : \M_\phi \cap B_{H^1_0(\Omega)}(\phi;\delta_0) \to H^1_0(\Omega)$ ($j=1,\ldots,p$) such that $T_u \M_\phi$ is spanned by $\{X_j(u)\}_{j=1,\ldots,p}$ for each $u \in \M_\phi \cap B_{H^1_0(\Omega)}(\phi;\delta_0)$. Since \eqref{a02} yields $(\psi_n,X_j(m_n))_{H^1_0(\Omega)} = 0$ for $n \in \N$ and $j = 1,\ldots,p$, we have
$$
(\psi_\infty,X_i(\phi))_{H^1_0(\Omega)} = 0 \quad \mbox{ for } \ j = 1,\ldots,p.
$$}
$$
(\psi_\infty,v)_{H^1_0(\Omega)} = 0 \ \mbox{ for all } \ v \in T_\phi \M_\phi,
$$
whence it follows that $\ker \mathcal{L}_\phi \supset \mathrm{span}\,\{T_\phi \M_\phi,\psi_\infty\}$. However, it contradicts the fact that $\ker \mathcal{L}_\phi = T_\phi \M_\phi$ along with $\psi_\infty \perp T_\phi \M_\phi$. This completes the proof.
\end{proof}

\section*{Declaration on the Use of AI Tools}

The authors occasionally consulted ChatGPT (OpenAI) for minor language editing and background information. All mathematical results and proofs in this manuscript were developed independently by the authors without the use of generative AI tools. The authors take full responsibility for the content of the manuscript.

\end{document}